\documentclass{amsart}
\usepackage{amssymb,mathtools}
\usepackage[british]{babel}
\usepackage{enumitem}
\usepackage[latin1]{inputenc}

\usepackage{url}

\usepackage{stmaryrd}
\usepackage{tikz-cd}

\newtheorem{theorem}{Theorem}
\numberwithin{theorem}{section}
\newtheorem{corollary}[theorem]{Corollary}
\newtheorem{lemma}[theorem]{Lemma}
\newtheorem{proposition}[theorem]{Proposition}
\theoremstyle{definition}
\newtheorem{definition}[theorem]{Definition}
\newtheorem{remark}[theorem]{Remark}
\newtheorem{example}[theorem]{Example}

\newtheorem{assumption}[theorem]{Standing Assumption}
\newtheorem*{claim}{Claim}

\newcommand\pto{\mathrel{\twoheadrightarrow_p}} 

\newcommand{\rca}{\mathsf{RCA}}
\newcommand{\atrs}{\mathsf{ATR_0^{set}}}
\newcommand{\lo}{\mathsf{LO}}
\newcommand{\tl}{\vartriangleleft}

\newcommand{\tr}{\operatorname{Tr}}
\newcommand{\nf}{\mathrel{=_{\operatorname{NF}}}}
\newcommand{\scr}{\mathsf{SC}}
\newcommand{\hau}{\mathsf{H}}
\newcommand{\sub}{\operatorname{sub}}

\newcommand{\rng}{\operatorname{rng}}
\newcommand{\dom}{\operatorname{dom}}
\newcommand{\To}{\Rightarrow}
\newcommand{\supp}{\operatorname{supp}}

\newcommand{\ord}{\operatorname{Ord}}
\newcommand{\en}{\operatorname{en}}
\newcommand{\len}{\operatorname{len}}
\newcommand{\ax}{\operatorname{Ax}}
\newcommand{\ot}{\mathbf O}

\newcommand{\rk}{\operatorname{rk}}
\newcommand{\ad}{\operatorname{Ad}}

\title[Well ordering principles for iterated $\Pi^1_1$-comprehension]{Well ordering principles\\ for iterated $\Pi^1_1$-comprehension}
\author{Anton Freund and Michael Rathjen}

\address{Anton Freund, Department of Mathematics, Technical University of Darmstadt, Schloss\-garten\-str.~7, 64289~Darmstadt, Germany}
\email{freund@mathematik.tu-darmstadt.de}

\address{Michael Rathjen, Department of Pure Mathematics, University of Leeds, Leeds LS2\,9JT, United Kingdom}
\email{rathjen@maths.leeds.ac.uk}

\thanks{The work of Anton Freund has been funded by the Deutsche Forschungsgemeinschaft (DFG, German Research Foundation) -- Project number 460597863. Michael Rathjen was partially supported by JTF grant 60842}

\begin{document}

\begin{abstract}
We introduce ordinal collapsing principles that are inspired by proof theory but have a set theoretic flavor. These principles are shown to be equivalent to iterated $\Pi^1_1$-comprehension and the existence of admissible sets, over weak base theories. Our work extends a previous result on the non-iterated case, which had been conjectured in Montalb\'an's ``Open questions in reverse mathematics" (Bull.\ Symb.\ Log.\ 17(3)2011). This previous result has already been applied to the reverse mathematics of combinatorial and set theoretic principles. The present paper is a significant contribution to a general approach that connects these~fields.
\end{abstract}


\maketitle

\section{Introduction}\label{sect:intro}

Well ordering principles assert that certain (computable) transformations of linear orders preserve well foundedness. Historically, the first example concerns the transformation of a linear order~$X$ into the set
\begin{equation*}
\omega(X):=\{\langle x_0,\ldots,x_{n-1}\rangle\,|\,x_0,\ldots,x_{n-1}\in X\text{ and }x_{n-1}\leq_X\ldots\leq_Xx_0\}
\end{equation*}
of finite non-increasing sequences in~$X$, ordered lexicographically. As shown by \mbox{J.-Y.}~Gi\-rard \cite[Theorem~5.4.1]{girard87} and J.~Hirst~\cite{hirst94}, the statement that `$\omega(X)$ is well founded whenever the same holds for~$X$' is equivalent to a set existence principle known as arithmetical comprehension. The latter is, in turn, equivalent to important mathematical results such as the Arzel\`a-Ascoli theorem or the infinite Ramsey theorem (for each fixed number of at least three colours). To make clear that these equivalences are informative, we point out that they are established in a weak base system $\rca_0$ (`recursive comprehension axiom'). They are part of a research programme known as `reverse mathematics', developed by H.~Friedman~\cite{friedman-rm} and S.~Simpson (see his textbook~\cite{simpson09} for a comprehensive introduction).

The literature contains many more equivalences between well ordering principles, statements about set existence, and mathematical theorems~\cite{rathjen-afshari,friedman-mw,marcone-montalban,rathjen-atr,rathjen-model-bi,rathjen-weiermann-atr,thomson-thesis,thomson-rathjen-Pi-1-1}. At the same time, there is a fundamental limitation: The statement that `$X$ is well founded' has complexity~$\Pi^1_1$ (one universal quantification over infinite sets). Given a computable transformation~$D$ of linear orders, the principle that `$D(X)$ is well founded whenever the same holds for~$X$' will thus be $\Pi^1_2$ (`for all -- exists'). It is known that principles of this form cannot be equivalent to more abstract set existence statements, such as the principle of $\Pi^1_1$-comprehension from reverse mathematics or the `minimal bad sequence lemma' of C.~Nash-Williams~\cite{nash-williams63} (see the analysis by A.~Marcone~\cite{marcone-bad-sequence}).

To overcome this limitation, one can consider order transformations of higher type, which have other transformations as arguments or values. More precisely, the latter should be dilators in the sense of J.-Y.~Girard~\cite{girard-pi2}, i.\,e., particularly uniform transformations $X\mapsto D(X)$ of well orders (see below for details). In the prime example from the literature, a given dilator~$D$ is transformed into a linear order~$\vartheta(D)$ that represents a relativized Bachmann-Howard ordinal (details below). The statement that `$\vartheta(D)$ is well founded for every dilator~$D$' is equivalent to the principle of $\Pi^1_1$-comprehension, as shown by the first author~\cite{freund-thesis,freund-equivalence,freund-categorical,freund-computable}. For related work by the second author we refer to~\cite{rathjen-axiomatic-thinking} and to Section~6 of the earlier paper~\cite{rathjen-atr}. The equivalence with $\Pi^1_1$-comprehension had been conjectured in A.~Montalb\'an's list of `Open questions in reverse mathematics'~\cite{montalban-open-problems}). 

The cited result on~$\vartheta(D)$ has become the basis for a very satisfactory analysis of the minimal bad sequence lemma in terms of a uniform Kruskal theorem~\cite{frw-kruskal}, for a new approach to Friedman's gap condition~\cite{freund-kruskal-gap,freund-BH-derivative}, and for another equivalence that involves patterns of resemblance~\cite{freund-patterns} (which resolves a further open question from Montalb\'an's list~\cite{montalban-open-problems}). These applications show why well ordering principles are relevant: they connect very intricate constructions from proof theory to reverse mathematics, set theory, and core mathematics. The present paper shows that these connections extend far beyond the existing literature. Specifically, we will study iterated $\Pi^1_1$-comprehension or, equivalently, hierarchies of admissible sets. In particular, we will obtain a characterization of $\Pi^1_1$-transfinite recursion, which is equivalent to mathematical results such as the Galvin-Prikry theorem from Ramsey theory (as shown by K.~Tanaka~\cite{tanaka89}). We will also characterize the statement that `every set is contained in a countable $\beta$-model of $\Pi^1_1$-comprehension', which solves an important case of the general Conjecture~6.1 from~\cite{rathjen-atr}.

Let us recall some terminology that is needed to state our result. We write $\lo$ for the category with linear orders as objects and embeddings (strictly increasing functions) as morphisms. By $[\cdot]^{<\omega}$ we denote the finite subset functor on the category of sets, with
\begin{align*}
[X]^{<\omega}&:=\text{`the set of finite subsets of~$X$'},\\
[f]^{<\omega}(a)&:=\{f(x)\,|\,x\in a\}\quad\text{(for $f:X\to Y$ and $a\in[X]^{<\omega}$)}.
\end{align*}
We will suppress the forgetful functor from linear orders to sets. In the following definition, this allows us to view both~$D$ and~$[\cdot]^{<\omega}$ as functors from linear orders to sets, so that we can consider a natural transformation between them. By $\rng(f)$ we denote the range (in the sense of `image') of a function~$f$.

\begin{definition}\label{def:dilator}
A predilator consists of a functor $D:\lo\to\lo$ and a natural transformation $\supp:D\Rightarrow[\cdot]^{<\omega}$ such that the `support condition'
\begin{equation*}
\rng(D(f))=\{\sigma\in D(Y)\,|\,\supp_Y(\sigma)\subseteq\rng(f)\}
\end{equation*}
is satisfied for every embedding~$f:X\to Y$ of linear orders. If $D(X)$ is well founded for any well order~$X$, then $D$ (together with $\supp$) is a dilator.
\end{definition}

Girard additionally demands that $D(f)\leq D(g)$ follows from $f\leq g$ (pointwise inequalities between morphisms), which is automatic for dilators but not for predilators (see~\cite[Proposition~2.3.10]{girard-pi2} or also~\cite[Lemma~5.3]{frw-kruskal}). Apart from this, our definition is equivalent to Girard's, which does not mention supports but demands that $D$ preserves direct limits and pullbacks (see~\cite[Remark~2.2.2]{freund-thesis}). Predilators are determined by their restrictions to the category of finite orders, essentially because any linear order is the union of its finite suborders. As observed by Girard, this allows us to treat predilators as sets (rather than proper classes) and to represent them in reverse mathematics (assuming their values on finite orders are countable). To make the present paper more readable, we will not work with representations explicitly. The reader who desires a detailed formalization of our considerations in reverse mathematics will find a blueprint in~\cite[Section~2]{freund-computable}.

The aforementioned characterization of $\Pi^1_1$-comprehension can now be made more precise. For a subset $a$ and an element~$y$ of a linear order~$X$, we write
\begin{equation*}
a\subseteq_X y\quad:\Leftrightarrow\quad x<_X y\text{ for all $x\in a$}.
\end{equation*}
This fits with the usual identification of ordinals with their sets of predecessors. The following notion -- first defined in~\cite{freund-equivalence} -- is inspired by Rathjen's notation system for the Bachmann-Howard ordinal (see~\cite{rathjen-weiermann-kruskal}).

\begin{definition}\label{def:Bachmann-Howard-fp}
A Bachmann-Howard collapse for a predilator~$D$ consists of a linear order~$X$ and a function $\vartheta:D(X)\to X$ such that
\begin{enumerate}[label=(\roman*)]
\item $\sigma<_{D(X)}\tau$ and $\supp_X(\sigma)\subseteq_X\vartheta(\tau)$ entail $\vartheta(\sigma)<_X\vartheta(\tau)$,
\item we have $\supp_X(\sigma)\subseteq_X\vartheta(\sigma)$ for all~$\sigma\in D(X)$.
\end{enumerate}
If such a $\vartheta$ exists, we call~$X$ a Bachmann-Howard fixed point of~$D$.
\end{definition} 

In~\cite[Section~4]{freund-computable} it is shown that any predilator~$D$ has a minimal Bachmann-Howard fixed point~$\vartheta(D)$, which is computable with a representation of~$D$ as oracle. We can now give a precise formulation of the result that was mentioned above.

\begin{theorem}[{\cite{freund-equivalence,freund-computable}}]\label{thm:Pi^1_1-CA}
The following are equivalent over~$\rca_0$:
\begin{enumerate}[label=(\roman*)]
\item $\Pi^1_1$-comprehension,
\item any dilator has a well founded Bachmann-Howard fixed point,
\item if $D$ is a dilator, then $\vartheta(D)$ is well founded.
\end{enumerate}
\end{theorem}

Let us point out that (ii) and (iii) have different virtues. Since $D\mapsto\vartheta(D)$ is a computable transformation, statement~(iii) is a well ordering principle of higher type, as discussed above. The explicit construction of~$\vartheta(D)$ reveals that the strength of~(ii) lies in well foundedness, not in the existence of Bachmann-Howard fixed points as linear orders. On the other hand, statement~(ii) has the advantage that it is very easy to formulate. This demonstrates another advantage of well ordering principles: they allow us to condense central ideas of ordinal analysis into elegant set theoretic principles. With a grain of salt, we suggest to view these principles as `large cardinal axioms' in the computable realm.

We now describe how Theorem~\ref{thm:Pi^1_1-CA} will be generalized in the present paper. The product $X\times Y$ of linear orders is defined as usual, namely by
\begin{equation*}
(x,y)<_{X\times Y}(x',y')\quad:\Leftrightarrow\quad\text{$x<_X x'$ or ($x=x'$ and $y<_Y y'$)}.
\end{equation*}
Given functions $f:X\to X'$ and $g:Y\to Y'$, we define $f\times g:X\times Y\to X'\times Y'$ by $(f\times g)(x,y):=(f(x),g(y))$. Note that we omit one pair of parentheses to improve readability. If $f$ or $g$ is the identity on $X=X'$ or $Y=Y'$, respectively, we write $X\times g$ or $f\times Y$ rather than $f\times g$. By Example~\ref{ex:Buchholz-psi}, the following generalizes the $\psi$-functions of W.~Buchholz~\cite{buchholz-new-system}.

\begin{definition}\label{def:nu-collapse}
Given a well order~$\nu$ and a predilator~$D$, a $\nu$-collapse for~$D$ consists of a linear order~$X$ and an embedding~$\pi:X\to\nu\times D(X)$ with the following two properties: First, we demand that the relation $\tl$ on~$X$ that is given by
\begin{equation*}
s\tl t\quad:\Leftrightarrow\quad s\in\supp_X(\tau)\text{ for }\pi(t)=(\alpha,\tau)
\end{equation*}
is well founded (think of $s$ as a subterm of~$t$). For $\gamma<\nu$, we use recursion along~$\tl$ to define $G^D_\gamma:X\to[D(X)]^{<\omega}$ and simultaneously $G_\gamma:D(X)\to[D(X)]^{<\omega}$ by
\begin{align*}
G^D_\gamma(t)&:=\begin{cases}
\{\tau\}\cup G_\gamma(\tau) & \text{if }\pi(t)=(\alpha,\tau)\text{ with }\alpha\geq\gamma,\\
\emptyset & \text{if }\pi(t)=(\alpha,\tau)\text{ with }\alpha<\gamma,
\end{cases}\\
G_\gamma(\tau)&:=\bigcup\{G^D_\gamma(s)\,|\,s\in\supp_X(\tau)\}.
\end{align*}
Secondly, we now demand that $\pi$ has range
\begin{equation*}
\rng(\pi)=\{(\alpha,\tau)\in\nu\times D(X)\,|\,G_\alpha(\tau)\subseteq_{D(X)}\tau\}.
\end{equation*}
If such a $\pi$ exists, we say that $X$ is a $\nu$-fixed point of~$D$.
\end{definition}

Concerning the formalization in~$\rca_0$, we note that all ranks with respect to the order~$\tl$ are finite (by the finiteness of supports), and that the functions $G_\gamma$ and $G_\gamma^D$ are computable. Instead of the embedding~$\pi$, we will often consider its partial inverse~$\psi:\nu\times D(X)\to_p X$, which can be seen as a collapsing function in the sense of impredicative ordinal analysis (see the following example). While some readers may prefer to reformulate the definition in terms of~$\psi$, we feel that the use of~$\pi$ has notational advantages. Note that we cannot expect $\psi$ to be total, because the order type of~$\nu\times D(X)$ will typically exceed the one of~$X$. Very roughly, the condition on $\rng(\pi)$ ensures that $\psi$ has a large domain of definition. Given that $\pi$ and hence $\psi$ is order preserving, this means that~$X$ must have large order type.

\begin{example}\label{ex:Buchholz-psi}
To turn the transformation $X\mapsto\omega(X)$ into a dilator, we declare
\begin{align*}
\omega(f)(\langle x_0,\ldots,x_{n-1}\rangle)&:=\langle f(x_0),\ldots,f(x_{n-1})\rangle,\\
\supp^\omega_X(\langle x_0,\ldots,x_{n-1}\rangle)&:=\{x_0,\ldots,x_{n-1}\}.
\end{align*}
Consider Buchholz' order~$\mathsf{OT}$ from~\cite[Section~2]{buchholz-new-system}, and let $\mathsf{P}\subseteq\mathsf{OT}$ be the suborder of principal terms, which have the form $D_\alpha t$ with $\alpha<\omega+1$ and $t\in\mathsf{OT}$. We want to show that $\mathsf P$ is an $(\omega+1)$-fixed point of the dilator~$\omega(\cdot)$. Up to the obvious isomorphism $\mathsf{OT}\cong\omega(\mathsf{P})$, we can define $\pi:\mathsf{P}\to(\omega+1)\times\omega(\mathsf P)$ by $\pi(D_\alpha t):=(\alpha,t)$. Clause (${\prec}2$) from the cited paper by Buchholz ensures that $\pi$ is an embedding. Given $s\tl D_\alpha t$ with $t=\langle t_0,\ldots,t_{n-1}\rangle$, we invoke the definition of $\tl$ to get
\begin{equation*}
s\in\supp^\omega_{\mathsf P}(t)=\{t_0,\ldots,t_{n-1}\}.
\end{equation*}
The latter entails that $s$ is a subterm of $D_\alpha t$ (in the usual sense), which ensures that $\tl$ is well founded. The isomorphism $\mathsf{OT}\cong\omega(\mathsf{P})$ identifies $t\in\mathsf{P}\subseteq\mathsf{OT}$ with the element $\langle t\rangle\in\omega(\mathsf{P})$. Up to this identification, the function $G_\gamma:\omega(\mathsf{P})\to[\omega(\mathsf{P})]^{<\omega}$ from Definition~\ref{def:nu-collapse} is an extension of $G^\omega_\gamma:\mathsf{P}\to[\omega(\mathsf{P})]^{<\omega}$. Based on this observation, one readily checks that our function $G_\gamma$ coincides with $G_\gamma:\mathsf{OT}\to[\mathsf{OT}]^{<\omega}$ as defined by Buchholz, still modulo $\mathsf{OT}\cong\omega(\mathsf{P})$. In view of Buchholz' clause~($\mathsf{OT}3$), it follows that $\pi$ has range as required by Definition~\ref{def:nu-collapse}.
\end{example}

In Section~\ref{sect:fixed-points-exist}, we explicitly construct a $\nu$-fixed point~$\psi_\nu(D)$ of a given predilator~$D$. More precisely, the order~$\psi_\nu(D)$ will be given as a term system that is computable relative to~$\nu$ and~$D$, so that its existence is known in the axiom system~$\rca_0$. We will also show that $\psi_\nu(D)$ is isomorphic to any other $\nu$-fixed point of~$D$, so that $\nu$-fixed points are essentially unique. This confirms the significance of Example~\ref{ex:Buchholz-psi}. Let us now state our main result, which is further explained below. The proof spans most of our paper and will be completed in Section~\ref{sect:conclude}.

\begin{theorem}\label{thm:main}
Provably in $\rca_0$, the following principles are equivalent for any infinite well order $\nu$:
\begin{enumerate}[label=(\roman*)]
\item $\Pi^1_1$-recursion along~$\nu$,
\item any dilator has a well founded $\nu$-fixed point,
\item if $D$ is a dilator, then $\psi_\nu(D)$ is well founded.
\end{enumerate}
Over $\atrs$, statements (i) to~(iii) are also equivalent to the following:
\begin{enumerate}[label=(\roman*)]\setcounter{enumi}{3}
\item for any set~$u$, there is a sequence of admissible sets $\mathsf{Ad}_\alpha\ni u$ for $\alpha<\nu$, such that $\alpha<\beta<\nu$ entails $\mathsf{Ad}_\alpha\in\mathsf{Ad}_\beta$ (where we consider~$\nu$ as an ordinal).
\end{enumerate}
\end{theorem}

The restriction to infinite~$\nu$ is convenient, because it will allow us to reduce to the case where~$\nu$ is of limit type. In $\rca_0$ one can also prove the equivalence  for $\nu=1$ and hence for each finite~$\nu$ that is fixed externally, as we shall see in Corollary~\ref{cor:1-fixed-points} (based on Theorem~\ref{thm:Pi^1_1-CA}). What we will not show is that $\rca_0$ proves the equivalence uniformly for all finite~$\nu$. We believe that this could be establish by our methods, but this would seem to require a separate treatment of the successor case, which we were keen to avoid.

Let us now explain statement~(i) from Theorem~\ref{thm:main}. Given $Y\subseteq\mathbb N$ and $\alpha<\nu$, we write $Y_\alpha$ for the set of all $x\in\mathbb N$ such that (the Cantor code of) the pair~$\langle\alpha,x\rangle$ is contained in~$Y$. In other words, we view $Y$ as a representation of the sequence of sets $Y_\alpha\subseteq\mathbb N$ with $\alpha<\nu$. Its initial segments are represented by the sets
\begin{equation*}
Y^\alpha:=\{\langle\gamma,x\rangle\in Y\,|\,\gamma<\alpha\}=\{\langle\gamma,x\rangle\in\alpha\times\mathbb N\,|\,x\in Y_\gamma\}\subseteq\mathbb N.
\end{equation*}
For a formula~$\varphi(x,\alpha,X)$, possibly with further parameters, let $H_\varphi(Y)$ be (the obvious formalization of) the statement
\begin{equation*}
H_\varphi(Y)\quad:\Leftrightarrow\quad Y_\alpha=\{x\in\mathbb N\,|\,\varphi(x,\alpha,Y^\alpha)\}\text{ for all }\alpha<\nu.
\end{equation*}
More intuitively, this expresses that the sets $Y_\alpha\subseteq\mathbb N$ are built by recursion along~$\nu$, where~$\varphi$ determines the recursion step. Let us recall that $\Pi^1_1$-formulas have the form $\forall X\subseteq\mathbb N.\,\theta$ for a formula $\theta$ that contains quantifiers $\forall n\in\mathbb N$ and $\exists n\in\mathbb N$ only. Statement~(i) from Theorem~\ref{thm:main} is the axiom schema that consists of all statements
\begin{equation*}
\forall x_1,\ldots,x_m\in\mathbb N\,\forall X_1,\ldots,X_n\subseteq\mathbb N\,\exists Y\subseteq\mathbb N.\, H_\varphi(Y)
\end{equation*}
for a $\Pi^1_1$-formula~$\varphi$ with number and set parameters $x_1,\ldots,x_m$ and~$X_1,\ldots,X_n$.

Before we discuss the axiom system $\atrs$ and statement~(iv) from Theorem~\ref{thm:main}, we consider some instances that are relevant in their own right (see Section~\ref{sect:conclude} for proofs). First, the following result was promised in~\cite{rathjen-axiomatic-thinking}, for a projected article with the title `A~proof-theoretic characterization of $\beta$-models of $\Pi^1_1$-comprehension', which we have incorporated into the present more general paper.

\begin{corollary}\label{cor:beta-model-Pi11-CA}
The following are equivalent over $\rca_0$:
\begin{enumerate}[label=(\roman*)]
\item every subset of~$\mathbb N$ is contained in a countable $\beta$-model of $\Pi^1_1$-comprehension,
\item any dilator has a well founded $\omega$-fixed point,
\item if $D$ is a dilator, then $\psi_\omega(D)$ is well founded.
\end{enumerate}
\end{corollary}

Secondly, the axiom schema and rule of $\Delta^1_2$-comprehension are closely connected to iterations of $\Pi^1_1$-recursion along fixed~$\nu<\varepsilon_0$ and $\nu<\omega^\omega$, respectively, as shown by H.~Friedman~\cite{friedman70} and S.~Feferman~\cite{feferman70} (see also the presentation by W.~Pohlers~\cite[Section~3.2]{pohlers98}). Our Theorem~\ref{thm:main} yields analogous connections with the well foundedness of $\nu$-fixed points. Finally, we obtain the following corollary when we quantify over~$\nu$. To confirm the significance of this result, we recall that $\Pi^1_1$-transfinite recursion is equivalent to the Galvin-Prikry theorem and to the principle of $\Delta^0_2$-determinacy, due to Tanaka~\cite{tanaka89,tanaka90}.

\begin{corollary}\label{cor:Pi11-TR}
The following are equivalent over $\rca_0$:
\begin{enumerate}[label=(\roman*)]
\item $\Pi^1_1$-transfinite recursion, i.\,e., the principle that $\Pi^1_1$-recursion is available along any well order~$\nu$,
\item any dilator has a well founded $\nu$-fixed point for every well order~$\nu$,
\item if $D$ is a dilator and $\nu$ is any well order, then $\psi_\nu(D)$ is well founded.
\end{enumerate}
\end{corollary}

Let us now complete our explanation of Theorem~\ref{thm:main}. The axiom system $\atrs$ is a set theory due to Simpson~\cite{simpson82,simpson09}, who showed that it is conservative over the axiom system~$\mathsf{ATR_0}$ (`arithmetical transfinite recursion') from reverse mathematics. Its axioms ensure that all primitive recursive set functions (in the sense of R.~Jensen and C.~Karp~\cite{jensen-karp}) are total and that every well order is isomorphic to an ordinal (`axiom beta'). We also include the axiom that all sets are countable, as in~\cite{simpson09} (while~\cite{simpson82} marks this axiom as `optional'). 

We also recall that an admissible set is a transitive model of Kripke-Platek set theory. For~$\nu=1$, the equivalence between (i) and~(iv) has been shown by G.~J\"ager~\cite{jaeger-admissibles} (see also~\cite[Section~1.4]{freund-thesis}). The extension to general~$\nu$ can probably be considered as known, but we will also obtain a new -- if rather indirect -- proof in the present paper. Indeed, we will work in $\atrs$ to prove the circle of implications
\begin{equation*}
(i)\quad\Rightarrow\quad(ii)\quad\Leftrightarrow\quad(iii)\quad\Rightarrow\quad(iv)\quad\Rightarrow\quad(i)
\end{equation*}
between the statements from Theorem~\ref{thm:main}. In order to obtain the equivalence of~(i),~(ii) and~(iii) over $\rca_0$, we will argue that each of these statements entails arithmetical transfinite recursion (consider Theorem~\ref{thm:1fp-to-Bachmann-Howard} together with Theorem~\ref{thm:Pi^1_1-CA} above). Note that~(iv) cannot be (directly) considered over~$\rca_0$, as it is a statement of set theory rather than reverse mathematics.

Statements~(ii) and~(iii) of Theorem~\ref{thm:main} are equivalent because~$\psi_\nu(D)$ is the unique $\nu$-fixed point of~$D$ (up to isomorphism), as mentioned above and proved in Section~\ref{sect:fixed-points-exist}. The implication from (i) to~(ii) is established in Section~\ref{sect:wo-proof}, where we relativize Buchholz'~\cite{buchholz-distinguished} method of `distinguished sets' to a given dilator (cf.~the relativization to a single order in~\cite[Section~12.3.1]{thomson-rathjen-Pi-1-1}). In Section~\ref{sect:conclude} we recall the standard proof that~(iv) implies~(i).

To prove the crucial implication from~(ii) to~(iv), we will generalize the argument that was given for~$\nu=1$ in~\cite{freund-equivalence}. There we developed a notion of \mbox{$\beta$-}proof (cf.~\cite{girard-intro}) that is sound and complete for the class of models $\mathbb L^u_\alpha$, i.\,e., the stages of the constructible hierarchy over a transitive~$u=:\mathbb L^u_0$. By completeness, the existence of an admissible set~$\mathbb L^u_\alpha$ (which implies~(i) of Theorem~\ref{thm:Pi^1_1-CA}) was reduced to the claim that there is no $\beta$-proof of contradiction in Kripke-Platek set theory. This claim is a natural target for ordinal analysis, which is specialized in consistency proofs based on large well orders. Specifically, one argues that the height of a given $\beta$-proof can be bounded by some dilator~$D$. Based on the well order~$\vartheta(D)$ from~(ii) of Theorem~\ref{thm:Pi^1_1-CA}, one can employ J\"ager's ordinal analysis of Kripke-Platek set theory~\cite{jaeger-kripke-platek}, to conclude that the given \mbox{$\beta$-}proof does not derive a contradiction.

In the argument from~\cite{freund-equivalence} that we have sketched in the previous paragraph, the relevant $\beta$-proofs consist of a tree $S_X$ for each linear order~$X$ (see~\cite[Section~4]{freund-equivalence}). The aforementioned dilator~$D$ is essentially given by $D(X)=S_X$ with the Kleene-Brouwer order. In the present paper, we obtain corresponding trees $S^R_X$ that depend not only on a linear order~$X$ but also on a given embedding~$R:\nu\to X$, which corresponds to the sequence of admissible sets in~(iv) of Theorem~\ref{thm:main} (see Section~\ref{sect:search-trees}). However, we cannot allow~$D(X)$ to depend on~$R$, because~(ii) of Theorem~\ref{thm:main} requires a dilator, i.\,e., a transformation whose arguments are linear orders without additional structure. This new obstacle is resolved in Section~\ref{sect:search-to-notations}, which can be seen as the main technical contribution of the present paper. To complete the proof that (ii) implies~(iv) in Theorem~\ref{thm:main}, we then adapt the classical ordinal analysis for iterated admissible sets, developed by J\"ager and Pohlers~\cite{jaeger-pohlers82} and streamlined by Buchholz~\cite{buchholz-local-predicativity} (see also the earlier work on inductive definitions~\cite{bfps-inductive} and the detailed results in~\cite{Rathjen_PhD_2013}). Our `abstract' version of this ordinal analysis is worked out in Sections~\ref{sect:operators} and~\ref{sect:ordinal-analysis}. In the final Section~\ref{sect:conclude}, we combine all previous work into official proofs of Theorem~\ref{thm:main} and Corollaries~\ref{cor:beta-model-Pi11-CA} and~\ref{cor:Pi11-TR}.

\section{Existence and uniqueness of $\nu$-fixed points}\label{sect:fixed-points-exist}

In the present section, we construct a $\nu$-fixed point~$\psi_\nu(D)$ of a given predilator~$D$ for an arbitrary well order~$\nu$. Before, we show that all $\nu$-fixed points of~$D$ are isomorphic, which will entail that~$\psi_\nu(D)$ is essentially unique. The following result is central for our uniqueness proof.

\begin{proposition}\label{prop:FP-initial}
For well orders~$\mu$ and $\nu$, consider a $\mu$-collapse $\pi:X\to\mu\times D(X)$ and a $\nu$-collapse $\kappa:Y\to\nu\times D(Y)$ of a predilator~$D$. Given an embedding~$I:\mu\to\nu$, there is a unique embedding $f:X\to Y$ such that
\begin{equation*}
\begin{tikzcd}
X\arrow[r,"\pi"]\arrow[d,swap,"f"] & \mu\times D(X)\arrow[d,"I\times D(f)"]\\
Y\arrow[r,"\kappa"] & \nu\times D(Y)
\end{tikzcd}
\end{equation*}
is a commutative diagram.
\end{proposition}
\begin{proof}
Write $\tl$ for the well founded relation on~$X$ that is given by Definition~\ref{def:nu-collapse}. To prepare the proof of existence, we establish a more general form of uniqueness. For the purpose of this proof, let us say that a (finite or infinite) set $a\subseteq X$ is closed if $s\tl t\in a$ implies $s\in a$. We write $\iota_a:a\hookrightarrow X$ for the inclusion. By the definition of $\tl$ and the support condition from Definition~\ref{def:dilator}, any closed~$a$ validates
\begin{equation*}
t\in a\,\Rightarrow\,\supp_X(\tau)\subseteq a=\rng(\iota_a)\,\Rightarrow\,\tau\in\rng(D(\iota_a))\qquad\text{for}\qquad\pi(t)=(\alpha,\tau).
\end{equation*}
Given that $D(\iota_a)$ is an embedding, we get a unique embedding $\pi_a$ such that
\begin{equation*}
\begin{tikzcd}
a\arrow[r,"\pi_a"]\arrow[d,swap,"\iota_a"] & \mu\times D(a)\arrow[d,"\operatorname{Id}\times D(\iota_a)"]\\
X\arrow[r,"\pi"] & \mu\times D(X)
\end{tikzcd}
\end{equation*}
commutes. By an $a$-approximation, we shall mean an embedding $f_a:a\to Y$ such that the diagram from the proposition commutes if we replace~$X,\pi,f$ by $a,\pi_a,f_a$. When $a$ is the entire order $X$, then the functions $\iota_a$ and $D(\iota_a)$ are the identity on $a=X$ and $D(a)=D(X)$, respectively, since $D$ is a functor. In this case, the functions $\pi_a$ and $\pi$ will thus coincide, which means that an $X$-approximation is a function~$f$ as in the proposition. Our strong form of uniqueness reads as follows.

\begin{claim}
Given any $a$-approximation~$f_a$ and $b$-approximation~$f_b$ for closed $a,b\subseteq X$, we have $f_a(t)=f_b(t)$ for all $t\in a\cap b$.
\end{claim}

\noindent To prove the claim, one checks that $c:=a\cap b$ is closed and that $f_a\!\restriction\!c$ and $f_b\!\restriction\!c$ are $c$-approximations (write $f_a\!\restriction\!c=f_a\circ\iota$ with $\iota:c\hookrightarrow a$). To conclude, we consider an arbitrary $c$-approximation~$f$ and show that its values are uniquely determined. Given $t\in c$, write $\pi_c(t)=(\alpha,\tau)$ and consider the inclusion $\iota:\supp_c(\tau)\hookrightarrow c$. By the support condition, we can write $\tau=D(\iota)(\tau_0)$, where $\tau_0$ is unique since $D(\iota)$ is an embedding. As $f$ is a $c$-approximation, we obtain
\begin{equation*}
\kappa\circ f(t)=(I\times D(f))\circ\pi_c(t)=(I(\alpha),D(f\circ\iota)(\tau_0)).
\end{equation*}
Given that $\kappa$ is an embedding, this means that $f(t)$ is determined by $f\circ\iota$. We can deduce uniqueness by induction over~$\tl$, as $s\in\rng(\iota)$ implies~$s\tl t$. To see the latter, note that we have
\begin{equation*}
\pi(t)=\pi\circ\iota_c(t)=(\operatorname{Id}\times D(\iota_c))\circ\pi_c(t)=(\alpha,D(\iota_c)(\tau)),
\end{equation*}
and that the naturality of $\supp:D\Rightarrow[\cdot]^{<\omega}$ yields
\begin{equation*}
\rng(\iota)=\supp_c(\tau)=[\iota_c]^{<\omega}\circ\supp_c(\tau)=\supp_X(D(\iota_c)(\tau)).
\end{equation*}
As a next step towards existence, we show that approximations can be combined:

\begin{claim}
Consider a family $\langle f_i\,|\,i\in I\rangle$ of $a_i$-approximations $f_i$ for closed $a_i\subseteq X$. The function $f:a=\bigcup_{i\in I}a_i\to Y$ with $f(t)=f_i(t)$ for $t\in a_i$ is an $a$-approximation.
\end{claim}

\noindent Note that $a$ is closed and that $f$ is well defined by the previous claim. To show that $f$ is an $a$-approximation, we need to consider at most two indices at a time, namely, when we check that $f$ is an order embedding. This means that the claim for general~$I$ reduces to the one for~$I=\{0,1\}$. We establish the latter by induction on the cardinality $|a_0\cup a_1|\in\mathbb N\cup\{\infty\}$. The crucial step is to show
\begin{equation*}
t_0<_X t_1\,\Rightarrow\,f_0(t_0)<_Y f_1(t_1)\qquad\text{for}\qquad t_i\in a_i.
\end{equation*}
Let $a_i'\subseteq a_i$ consist of the predecessors of~$t_i$ in the transitive closure of~$\tl$. Then the set $c:=a_0'\cup a_1'$ is finite and cannot contain both $t_0$ and $t_1$, as $\tl$ is well founded. Due to the induction hypothesis, the restrictions $f_i\!\restriction\!a_i'$ can thus be combined into a $c$-approximation~$f'$. Put $\pi_i:=\pi_d$ with $d=a_i$. As in the proof of uniqueness, we can write $\pi_i(t_i)=(\alpha_i,D(\iota_i')(\tau_i))$ with $\iota_i':a_i'\hookrightarrow a_i$. For $\iota_i:a_i\hookrightarrow X$ we get
\begin{equation*}
\pi(t_i)=\pi\circ\iota_i(t_i)=(\operatorname{Id}\times D(\iota_i))\circ\pi_i(t_i)=(\alpha_i,D(\iota_i\circ\iota_i')(\tau_i)).
\end{equation*}
Let us also consider the inclusions $\iota_i'':a_i'\hookrightarrow c$ and $\iota_c:c\hookrightarrow X$. Clearly,
\begin{equation*}
\begin{tikzcd}
a_i'\arrow[r,,hook,"{\iota_i'}"]\arrow[d,hook,swap,"{\iota_i''}"] & a_i\arrow[d,hook,"\iota_i"]\\
c\arrow[r,hook,"\iota_c"] & X
\end{tikzcd}
\end{equation*}
is a commutative diagram. Aiming at the implication above, we now assume $t_0<t_1$. As $\pi$ is an embedding, we get either $\alpha_0<\alpha_1$ or $\alpha_0=\alpha_1$ and
\begin{equation*}
D(\iota_c)\circ D(\iota_0'')(\tau_0)=D(\iota_0\circ\iota_0')(\tau_0)<D(\iota_1\circ\iota_1')(\tau_1)=D(\iota_c)\circ D(\iota_1'')(\tau_1),
\end{equation*}
which entails $D(\iota_0'')(\tau_0)<D(\iota_1'')(\tau_1)$. By the choice of~$f'$ we have $f'\!\restriction\!a_i'=f_i\!\restriction\!a_i'$, or equivalently $f'\circ\iota_i''=f_i\circ\iota_i'$. Hence the last inequality entails
\begin{equation*}
D(f_0\circ\iota_0')(\tau_0)=D(f')\circ D(\iota_0'')(\tau_0)<D(f')\circ D(\iota_1'')(\tau_1)=D(f_1\circ\iota_1')(\tau_1).
\end{equation*}
To conclude $f_0(t_0)<f_1(t_1)$, it is thus enough to observe
\begin{equation*}
\kappa\circ f_i(t_i)=(I\times D(f_i))\circ\pi_i(t_i)=(I(\alpha_i),D(f_i\circ\iota_i')(\tau_i)).
\end{equation*}
Now that this second claim is proved, the proposition is reduced to the following:

\begin{claim}
Given any $t\in X$, there is an $a$-approximation for some finite closed $a\ni t$.
\end{claim}

\noindent Arguing by induction on~$\tl$, we can use the previous claim to produce a $b$-approxi\-mation $f$ for some finite closed $b\subseteq X$ that contains all $s\tl t$. As before, we can write $\pi(t)=(\alpha,D(\iota_b)(\tau))$ with $\iota_b:b\to X$. To extend $f$ into a function $f':a\to Y$ on the closed set $a:=b\cup\{t\}$, we would like to stipulate $\kappa\circ f'(t)=(I(\alpha),D(f)(\tau))$. For this purpose, we need to show that the right side lies in the range of~$\kappa$. Let us write $G^{D,Z}_\gamma:Z\to[D(Z)]^{<\omega}$ and $G^{Z}_\gamma:D(Z)\to[D(Z)]^{<\omega}$ for the functions from Definition~\ref{def:nu-collapse}, where $Z$ can be $X$ or~$Y$. Analogous functions for $Z=b$ arise by
\begin{align*}
G^{D,b}_\gamma(s)&:=\begin{cases}
\{\sigma\}\cup G^b_\gamma(\sigma) & \text{if }\pi_b(s)=(\alpha,\sigma)\text{ with }\alpha\geq\gamma,\\
\emptyset & \text{if }\pi_b(s)=(\alpha,\sigma)\text{ with }\alpha<\gamma,
\end{cases}\\
G^b_\gamma(\sigma)&:=\bigcup\{G^{D,b}_\gamma(r)\,|\,r\in\supp_b(\sigma)\}.
\end{align*}
To see that this recursion is well founded, note that $\pi_b(s)=(\alpha,\sigma)$ and $r\in\supp_b(\sigma)$ entail $r\tl s$, as in the proof of the first claim. By induction along~$\tl$ we get
\begin{equation*}
\begin{aligned}
[D(\iota_b)]^{<\omega}\circ G^{D,b}_\gamma&=G^{D,X}_\gamma\circ\iota_b,\\
[D(\iota_b)]^{<\omega}\circ G^{b}_\gamma&=G^{X}_\gamma\circ D(\iota_b),
\end{aligned}\qquad
\begin{aligned}
[D(f)]^{<\omega}\circ G^{D,b}_\gamma&=G^{D,Y}_{I(\gamma)}\circ f,\\
[D(f)]^{<\omega}\circ G^{b}_\gamma&=G^{Y}_{I(\gamma)}\circ D(f).
\end{aligned}
\end{equation*}
For $t\in X$ with $\pi(t)=(\alpha,D(\iota_b)(\tau))$ as above, we can invoke Definition~\ref{def:nu-collapse} to get
\begin{equation*}
[D(\iota_b)]^{<\omega}\circ G^{b}_\alpha(\tau)=G^{X}_\alpha\circ D(\iota_b)(\tau)\subseteq_{D(X)}D(\iota_b)(\tau).
\end{equation*}
The latter entails $G^{b}_\alpha(\tau)\subseteq_{D(b)}\tau$ and then
\begin{equation*}
G^{Y}_{I(\alpha)}\circ D(f)(\tau)=[D(f)]^{<\omega}\circ G^{b}_\alpha(\tau)\subseteq_{D(Y)}D(f)(\tau).
\end{equation*}
Again by Definition~\ref{def:nu-collapse}, it follows that $(I(\alpha),D(f)(\tau))$ lies in the range of~$\kappa$. As~indicated above we can thus define $f':a=b\cup\{t\}\to Y$ by stipulating
\begin{equation*}
\kappa\circ f'(t)=(I(\alpha),D(f)(\tau))
\end{equation*}
and $f'\!\restriction\!b=f$. The fact that~$f'$ is order preserving is readily deduced from the following observation: For $s\in b$ with $\pi_b(s)=(\beta,\sigma)$ we have $\pi(s)=(\beta,D(\iota_b)(\sigma))$, and since $f$ is a $b$-approximation we get
\begin{equation*}
\kappa\circ f'(s)=\kappa\circ f(s)=(I\times D(f))\circ\pi_b(s)=(I(\beta),D(f)(\sigma)).
\end{equation*}
To see that the diagram from the proposition commutes with $c,\pi_c,f'$ at the place of~$X,\pi,f$, we note that $f'\!\restriction\!b=f$ amounts to $f=f'\circ\iota$ with $\iota:b\hookrightarrow c$. For~$s\in b$ or~$s=t$, we see that $\pi(s)=(\beta,D(\iota_b)(\sigma))$ yields $\pi_c(s)=(\beta,D(\iota)(\sigma))$ and hence
\begin{equation*}
(I\times D(f'))\circ\pi_c(s)=(I(\beta),D(f'\circ\iota)(\sigma))=(I(\beta),D(f)(\sigma)),
\end{equation*}
which coincides with $\kappa\circ f'(s)$ as computed above.
\end{proof}

In terminology from category theory, the proposition shows that any $\nu$-fixed point satisfies the universal property of an initial object. As the following proof makes explicit, this entails that $\nu$-fixed points are essentially unique. For an application of Proposition~\ref{prop:FP-initial} with $\mu<\nu$, we refer to Corollary~\ref{cor:no-fp-monotone} below.

\begin{corollary}\label{cor:nu-fp-unique}
All $\nu$-fixed points of a given predilator are order isomorphic.
\end{corollary}
\begin{proof}
Consider $\nu$-fixed points $\pi:X\to\nu\times D(X)$ and $\kappa:Y\to\nu\times D(Y)$, and write $I:\nu\to\nu$ for the identity. Two applications of the previous proposition (one with $X$ and $Y$ interchanged) yield embeddings $f:X\to Y$ and $g:Y\to X$ with
\begin{equation*}
\pi\circ g\circ f=(I\times D(g))\circ\kappa\circ f=(I\times D(g))\circ(I\times D(f))\circ\pi=(I\times D(g\circ f))\circ\pi.
\end{equation*}
If $\operatorname{Id}_X$ is identity on~$X$, then $D(\operatorname{Id}_X)$ is the identity on~$D(X)$, as~$D$ is a functor. Hence we also have $\pi\circ\operatorname{Id}_X=(I\times D(\operatorname{Id}_X))\circ\pi$. We can conclude $g\circ f=\operatorname{Id}_X$ by the uniqueness part of the previous proposition. The analogous argument shows that $f\circ g$ is the identity on~$Y$, so that $f$ is indeed an isomorphism.
\end{proof}

To prepare the construction of $\nu$-fixed points, we recall a notion of normal form that is due to Girard~\cite{girard-pi2}. Where the context suggests it, we identify $n\in\mathbb N$ and the finite order~$\{0,\ldots,n-1\}$ (with the usual order between natural numbers). We also agree to write $|a|=\{0,\ldots,|a|-1\}$ for the cardinality of a finite set~$a$.

\begin{definition}\label{def:normal-form}
The trace of a predilator~$D$ is defined as
\begin{equation*}
\tr(D):=\{(n,\sigma)\,|\,n\in\mathbb N\text{ and }\sigma\in D(n)\text{ with }\supp_n(\sigma)=n\}.
\end{equation*}
We say that $\sigma\in D(X)$ has normal form $\sigma\nf D(e)(\sigma_0)$ with $e:n\to X$ for some~$n\in\mathbb N$ if we have $(n,\sigma_0)\in\tr(D)$ and $\sigma$ is indeed equal to $D(e)(\sigma_0)$.
\end{definition}

Let us recall a standard observation:

\begin{lemma}\label{lem:NF-unique}
Any $\sigma\in D(X)$ has a unique normal form $\sigma\nf D(e)(\sigma_0)$.
\end{lemma}
\begin{proof}
If $\sigma$ has normal form as given, then $e$ is determined as the unique embedding with domain $n:=|\supp_X(\sigma)|$ and range~$\supp_X(\sigma)\subseteq X$, as naturality yields
\begin{equation*}
\supp_X(\sigma)=\supp_X\circ D(e)(\sigma_0)=[e]^{<\omega}\circ\supp_n(\sigma_0)=[e]^{<\omega}(n)=\rng(e).
\end{equation*}
For existence, consider~$e$ as determined. The support condition from Definition~\ref{def:dilator} ensures that $\sigma=D(e)(\sigma_0)$ holds for some $\sigma_0\in D(n)$. By the equations above, we see that $\supp_X(\sigma)=\rng(e)$ entails $\supp_n(\sigma_0)=n$ and hence $(n,\sigma_0)\in\tr(D)$.
\end{proof}

In order to construct a $\nu$-fixed point~$\psi_\nu(D)$ of a given predilator~$D$, we shall first build an order $\psi_\nu^+(D)\supseteq\psi_\nu(D)$ that admits an order isomorphism
\begin{equation*}
\psi_\nu^+(D)\cong\nu\times D\left(\psi_\nu^+(D)\right).
\end{equation*}
We will later show that $\psi_\nu(D)$ is well founded when $D$ is a dilator (cf.~Theorem~\ref{thm:main}). The same cannot hold for~$\psi_\nu^+(D)$, which explains the auxiliary status of this order. Indeed, when we have $\nu>1$ and~$D$ admits embeddings $X\hookrightarrow D(X)$, then the order type of~$\nu\times D(X)$ will always exceed the one of~$X$.

\begin{definition}\label{def:psi-plus}
Consider an ordinal~$\nu$ and a predilator~$D$. The set $\psi_\nu^+(D)$ of terms is generated by the following recursive clause: Given a finite set $a\subseteq\psi_\nu^+(D)$, we add a term $\psi_\alpha(a,\sigma)\in\psi_\nu^+(D)$ for each $\alpha<\nu$ and each $\sigma\in D(|a|)$ with $(|a|,\sigma)\in\tr(D)$.
\end{definition}

Note that $\psi^+_\nu(D)$ is non-empty if the same holds for~$D(0)$. Let us consider
\begin{equation*}
l:\psi_\nu^+(D)\to\mathbb N\quad\text{with}\quad l\left(\psi_\alpha(a,\sigma)\right):=1+\textstyle\sum_{t\in a}2\cdot l(t).
\end{equation*}
The following definition determines $s\leq t$ by recursion on $l(s)+l(t)$. In particular, the factor~$2$ in the definition of~$l$ allows us to determine the restriction of $\preceq$ to~$a\cup b$. We demand that this restriction is linear, to ensure that $D(a\cup b)$ is defined.

\begin{definition}\label{def:order-psi-plus}
In order to define a binary relation $\preceq$ on $\psi_\nu^+(D)$ by recursion, we declare that $\psi_\alpha(a,\sigma)\preceq\psi_\beta(b,\tau)$ holds precisely if $a\cup b$ is linearly ordered by $\preceq$ and
\begin{enumerate}[label=(\roman*)]
\item we have $\alpha<\beta$,
\item or we have $\alpha=\beta$ and $D(e_a)(\sigma)\leq_{D(a\cup b)}D(e_b)(\tau)$ for the strictly increasing functions $e_a:|a|\to a\cup b$ and $e_b:|b|\to a\cup b$ with range $a$ and $b$, respectively.
\end{enumerate}
\end{definition}

The condition that $a\cup b$ is linearly ordered is made redundant by the following.

\begin{lemma}\label{lem:psi-plus-linear}
The relation $\preceq$ is a linear order on~$\psi_\nu^+(D)$.
\end{lemma}
\begin{proof}
By induction on~$n\in\mathbb N$, one can simultaneously show
\begin{alignat*}{3}
t&\preceq t\quad && \text{for }l(t)<n,\\
r\preceq s\text{ and }s&\preceq t\text{ imply }r\preceq t\quad && \text{for }l(r)+l(s)+l(t)<n,\\
s\preceq t\text{ and }t&\preceq s\text{ imply }s=t\quad && \text{for }l(s)+l(t)<n,\\
s\preceq t\text{ }&\text{or }t\preceq s\quad && \text{for }l(s)+l(t)<n.
\end{alignat*}
Let us establish transitivity for $r=\psi_\alpha(a,\rho)$, $s=\psi_\beta(b,\sigma)$ and $t=\psi_\gamma(c,\tau)$. The induction hypothesis ensures that $\preceq$ is linear on $d:=a\cup b\cup c$ (due to the factor~$2$ in the definition of~$l$ and since transitivity is trivial when all three relevant terms are equal). Given $r\preceq s$ and $s\preceq t$, the conclusion $r\preceq t$ is immediate unless we have $\alpha=\beta=\gamma$ as well as
\begin{equation*}
D(e_a^{a\cup b})(\rho)\leq_{D(a\cup b)} D(e_b^{a\cup b})(\sigma)\quad\text{and}\quad D(e_b^{b\cup c})(\sigma)\leq_{D(b\cup c)} D(e_c^{b\cup c})(\tau),
\end{equation*}
where $e_u^v:|u|\to v$ is strictly increasing with range~$u\subseteq v$. Note that $\iota_v^w\circ e_u^v=e_u^w$ holds for the inclusion $\iota_v^w:v\hookrightarrow w$. After composing the previous inequalities with $D(\iota_{a\cup b}^d)$ and $D(\iota_{b\cup c}^d)$, respectively, we can invoke transitivity in~$D(d)$ to get
\begin{equation*}
D(\iota_{a\cup c}^d)\circ D(e_a^{a\cup c})(\rho)=D(e_a^d)(\rho)\leq_{D(d)} D(e_c^d)(\tau)=D(\iota_{a\cup c}^d)\circ D(e_c^{a\cup c})(\tau).
\end{equation*}
We obtain $D(e_a^{a\cup c})(\rho)\leq_{D(a\cup c)}D(e_c^{a\cup c})(\tau)$, so that clause~(ii) of Definition~\ref{def:order-psi-plus} yields the desired inequality~$r\preceq t$. By similar but easier arguments, we can reduce the reflexivity and linearity of~$\preceq$ to the corresponding properties of orders~$D(d)$. To~establish antisymmetry, we must show that $s=t$ follows from
\begin{equation*}
D(e_b^{b\cup c})(\sigma)=D(e_c^{b\cup c})(\tau).
\end{equation*}
The expressions on both sides of this equation are normal forms in the sense of Definition~\ref{def:normal-form}, as Definition~\ref{def:psi-plus} ensures that $(|b|,\sigma)$ and $(|c|,\tau)$ lie in~$\tr(D)$. Hence Lemma~\ref{lem:NF-unique} allows us to conclude.
\end{proof}

To obtain an order isomorphism $\psi_\nu^+(D)\cong\nu\times D\left(\psi_\nu^+(D)\right)$ as promised above, it suffices to map $\psi_\alpha(a,\sigma)$ to $(\alpha,D(e_a)(\sigma))$, where $e_a:|a|\to\psi_\nu^+(D)$ is strictly increasing with range~$a$. This fact will not be used, but a very similar result is shown in the proof of Theorem~\ref{thm:psi-is-fp} below. We now single out the desired suborder.

\begin{definition}\label{def:psi}
In the following, let $e_a:|a|\to\psi^+_\nu(\sigma)$ denote the strictly increasing function with range~$a$ and the indicated codomain. For each ordinal $\gamma<\nu$ we define a function $G^+_\gamma:\psi_\nu^+(D)\to[D(\psi_\nu^+(D))]^{<\omega}$ by recursion over terms, stipulating
\begin{equation*}
G^+_\gamma(\psi_\alpha(a,\sigma)):=\begin{cases}
\{D(e_a)(\sigma)\}\cup\bigcup\{G^+_\gamma(r)\,|\,r\in a\} & \text{if $\alpha\geq\gamma$},\\
\emptyset & \text{if $\alpha<\gamma$}.
\end{cases}
\end{equation*}
The suborder $\psi_\nu(D)\subseteq\psi_\nu^+(D)$ is determined by the recursive clause
\begin{equation*}
\psi_\alpha(a,\sigma)\in\psi_\nu(D)\quad:\Leftrightarrow\quad a\subseteq\psi_\nu(D)\text{ and }\bigcup\{G^+_\alpha(r)\,|\,r\in a\}\subseteq_{D(\psi_\nu^+(D))}D(e_a)(\sigma).
\end{equation*}
\end{definition}

Let us now establish the main result of this section.

\begin{theorem}\label{thm:psi-is-fp}
The order~$\psi_\nu(D)$ is a $\nu$-fixed point of a given predilator~$D$.
\end{theorem}
\begin{proof}
Write $\iota:\psi_\nu(D)\hookrightarrow\psi^+_\nu(D)$ for the inclusion and $e'_a:|a|\to\psi_\nu(D)$ for the strictly increasing function with range~$a$, so that $e_a=\iota\circ e'_a$ is the same function as in Definition~\ref{def:psi}. Now consider the function
\begin{equation*}
\pi:\psi_\nu(D)\to\nu\times D(\psi_\nu(D))\quad\text{with}\quad\pi(\psi_\alpha(a,\sigma)):=(\alpha,D(e'_a)(\sigma)).
\end{equation*}
One readily shows that $\pi(s)\leq\pi(t)$ entails $s\preceq t$ (factorize $e'_a=\iota_{a\cup b}\circ e_a^{a\cup b}$ with $\iota_{a\cup b}:a\cup b\hookrightarrow\psi_\nu(D)$ as in the proof of Lemma~\ref{lem:psi-plus-linear}). Since the codomain of~$\pi$ is a linear order, it follows that $\pi$ is an embedding. With $X:=\psi_\nu(D)$ we compute
\begin{equation*}
\supp_X\left(D(e'_a)(\sigma)\right)=[e'_a]^{<\omega}\left(\supp_{|a|}(\sigma)\right)=[e'_a]^{<\omega}(|a|)=a.
\end{equation*}
Here the first equality holds since $\supp:D\Rightarrow[\cdot]^{<\omega}$ is natural, while the second one relies on $(|a|,\sigma)\in\tr(D)$ according to Definition~\ref{def:psi-plus}. The binary relation~$\tl$ that is determined in Definition~\ref{def:nu-collapse} can thus be characterized by
\begin{equation*}
s\tl\psi_\alpha(a,\sigma)\quad\Leftrightarrow\quad s\in a,
\end{equation*}
which entails that it is well founded. Let the functions $G^D_\gamma:\psi_\nu(D)\to[D(\psi_\nu(D))]^{<\omega}$ and $G_\gamma:D(\psi_\nu(D))\to[D(\psi_\nu(D))]^{<\omega}$ be given as in Definition~\ref{def:nu-collapse}. By induction along~$\tl$ one readily shows
\begin{equation*}
G^+_\gamma\circ\iota=[D(\iota)]^{<\omega}\circ G^D_\gamma.
\end{equation*}
In view of Definition~\ref{def:psi}, we can deduce that $\psi_\alpha(a,\sigma)\in\psi_\nu(D)$ entails
\begin{equation*}
[D(\iota)]^{<\omega}\circ G_\alpha(D(e'_a)(\sigma))=\bigcup\{G^+_\alpha(r)\,|\,r\in a\}\subseteq_{D(\psi_\nu^+(D))} D(\iota)\circ D(e'_a)(\sigma).
\end{equation*}
It follows that we have
\begin{equation*}
\rng(\pi)\subseteq\{(\alpha,\tau)\in\nu\times D(\psi_\nu(D))\,|\,G_\alpha(\tau)\subseteq_{D(\psi_\nu(D))}\tau\},
\end{equation*}
as Definition~\ref{def:nu-collapse} demands. To show that the converse of this inclusion holds as well, we consider an arbitrary element $(\alpha,\tau)$ of the right side. Writing $X=\psi_\nu(D)$, we put $a:=\supp_X(\tau)$. The support condition from Definition~\ref{def:dilator} yields a $\sigma\in D(|a|)$ with $\tau=D(e'_a)(\sigma)$. As in the proof of Lemma~\ref{lem:NF-unique} we get $(|a|,\sigma)\in\tr(D)$, which allows us to form the term $\psi_\alpha(a,\sigma)\in\psi_\nu^+(D)$. Given $G_\alpha(\tau)\subseteq\tau$, we get
\begin{equation*}
\bigcup\{G^+_\alpha(r)\,|\,r\in a\}=[D(\iota)]^{<\omega}\circ G_\alpha(\tau)\subseteq_{D(\psi^+_\nu(D))} D(\iota)(\tau)=D(e_a)(\sigma).
\end{equation*}
This entails that $\psi_\alpha(a,\sigma)$ does even lie in~$\psi_\nu(D)$. By construction, we can now conclude that $(\alpha,\tau)=\pi(\psi_\alpha(a,\sigma))$ is contained in the range of~$\pi$.
\end{proof}

By Corollary~\ref{cor:nu-fp-unique}, any $\nu$-fixed point of~$D$ is isomorphic to~$\psi_\nu(D)$, which confirms that statements~(ii) and~(iii) from Theorem~\ref{thm:main} are equivalent. If the equivalence with~(i) is to hold, then~(iii) must become stronger as $\nu$ grows. We conclude the section with a direct proof that this is the case.

\begin{corollary}\label{cor:no-fp-monotone}
If $\psi_\nu(D)$ is well founded, then so is $\psi_\mu(D)$ for any $\mu<\nu$.
\end{corollary}
\begin{proof}
Given $\mu<\nu$, there is an embedding of~$\mu$ into~$\nu$. By Proposition~\ref{prop:FP-initial} (which applies due to Theorem~\ref{thm:psi-is-fp}), we get an embedding of $\psi_\mu(D)$ into $\psi_\nu(D)$.
\end{proof}

\section{A proof of well foundedness}\label{sect:wo-proof}

In this section, we prove that~(i) implies~(ii) in Theorem~\ref{thm:main}, i.\,e., we use iterated $\Pi^1_1$-comprehension to show that $\nu$-fixed points of dilators are well \mbox{founded}. To make the general case more transparent, we provide an argument for $\nu=1$ first.

\begin{remark}\label{rmk:nu-one}
We show that any $1$-fixed point~$X$ of a dilator~$D$ is well founded. Consider a $1$-collapse $\pi:X\to D(X)$, where $D(X)$ is identified with $1\times D(X)$. Up to this identification, Definition~\ref{def:nu-collapse} yields
\begin{equation*}
s\tl t\quad\Leftrightarrow\quad s\in{\supp_X}\circ\pi(t),
\end{equation*}
and the definitions of $G^D_0:X\to[D(X)]^{<\omega}$ and $G_0:D(X)\to[D(X)]^{<\omega}$ become
\begin{equation*}
G_0^D(t)=\{\pi(t)\}\cup G_0(\pi(t))\quad\text{and}\quad G_0(\tau)=\bigcup\{G_0^D(s)\,|\,s\in\supp_X(\tau)\}.
\end{equation*}
Furthermore, the condition on the range of~$\pi$ does now read
\begin{equation*}
\rng(\pi)=\{\tau\in D(X)\,|\,G_0(\tau)\subseteq_{D(X)}\tau\}.
\end{equation*}
As a special feature of the case $\nu=1$, we get
\begin{equation*}
s\tl t\quad\Rightarrow\quad \pi(s)\in G_0^D(s)\subseteq G_0(\pi(t))\subseteq_{D(X)}\pi(t)\quad\Rightarrow\quad s<t.
\end{equation*}
Assuming $\Pi^1_1$-comprehension, we may form the well founded part~$W$ of~$X$, which can be given as the intersection of all sets $Z\subseteq X$ such that we have $t\in Z$ whenever $s\in Z$ holds for all~$s<_X t$. One readily shows that~$W$ is well founded with
\begin{equation*}
t\in W\quad\Leftrightarrow\quad s\in W\text{ for all }s\in X\text{ with }s<_Xt.
\end{equation*}
Write $\iota:W\hookrightarrow X$ for the inclusion. By the previous observations and the support condition from Definition~\ref{def:dilator}, we get
\begin{equation*}
t\in W\quad\Rightarrow\quad{\supp_X}\circ\pi(t)\subseteq W=\rng(\iota)\quad\Rightarrow\quad\pi(t)\in\rng(D(\iota)).
\end{equation*}
It follows that there is a function
\begin{equation*}
\kappa:W\to D(W)\quad\text{with}\quad D(\iota)\circ\kappa=\pi\circ\iota.
\end{equation*}
We will show that~$\kappa$ is a $1$-collapse of~$D$. Once this has been achieved, we can invoke Corollary~\ref{cor:nu-fp-unique} to learn that $X\cong W$ is well founded, as desired. In fact, the existence part of Proposition~\ref{prop:FP-initial} yields an embedding $f:X\to W$ with $\kappa\circ f=(I\times D(f))\circ\pi$, where $I:\nu\to\nu$ is the identity. By the uniqueness part of the same proposition, the composition $\iota\circ f$ must be the identity on~$X=W$. It remains to show that $\kappa$ satisfies the conditions from Definition~\ref{def:nu-collapse}. The latter ensures that~$\pi$ is an order embedding, so that the same holds for~$\kappa$. Given $s,t\in W$, we observe that the naturality of~$\supp:D\Rightarrow[\cdot]^{<\omega}$ yields
\begin{equation*}
[\iota]^{<\omega}\circ{\supp_W}\circ\kappa(t)={\supp_X}\circ D(\iota)\circ\kappa(t)={\supp_X}\circ\pi\circ\iota(t),
\end{equation*}
so that $\iota(s)\tl\iota(t)$ is equivalent to $s\in{\supp_W}\circ\kappa(t)$. This shows that the restriction of~$\tl$ to~$W$ coincides with the relation that~$\kappa$ induces according to Definition~\ref{def:nu-collapse}. The latter also yields functions $G_W^D:W\to[D(W)]^{<\omega}$ and $G_W:D(W)\to[D(W)]^{<\omega}$, which are given by
\begin{equation*}
G_W^D(t)=\{\kappa(t)\}\cup G_W(\kappa(t))\quad\text{and}\quad G_W(\tau)=\bigcup\{G_W^D(s)\,|\,s\in\supp_W(\tau)\}.
\end{equation*}
A straightforward induction along~$\tl$ shows that we have
\begin{equation*}
[D(\iota)]^{<\omega}\circ G^D_W=G_0^D\circ\iota\quad\text{and}\quad [D(\iota)]^{<\omega}\circ G_W=G_0\circ D(\iota).
\end{equation*}
By the aforementioned condition on the range of~$\pi$, we obtain
\begin{equation*}
[D(\iota)]^{<\omega}\circ G_W\circ\kappa(t)=G_0\circ D(\iota)\circ\kappa(t)=G_0\circ\pi\circ\iota(t)\subseteq_{D(X)}\pi\circ\iota(t)=D(\iota)\circ\kappa(t)
\end{equation*}
for any~$t\in W$. Since $D(\iota)$ is an embedding, we can conclude
\begin{equation*}
\rng(\kappa)\subseteq\{\tau\in D(W)\,|\,G_W(\tau)\subseteq_{D(W)}\tau\}.
\end{equation*}
It remains to establish the converse inclusion. Note that~$D(W)$ is well founded, as~$D$ is a dilator and~$W$ is a well order. We argue by (main) induction on~$\tau\in D(W)$ to prove the crucial implication
\begin{equation*}
G_W(\tau)\subseteq_{D(W)}\tau\quad\Rightarrow\quad\tau\in\rng(\kappa).
\end{equation*}
Assuming the premise, we get $G_0(D(\iota)(\tau))\subseteq_{D(X)} D(\iota)(\tau)$ as above, which allows us to write $D(\iota)(\tau)=\pi(t)$ with~$t\in X$. We will show $t\in W$, so that we obtain
\begin{equation*}
D(\iota)\circ\kappa(t)=\pi\circ\iota(t)=\pi(t)=D(\iota)(\tau).
\end{equation*}
Since~$D(\iota)$ is an embedding, we can conclude $\tau=\kappa(t)\in\rng(\kappa)$ as desired. In order to get~$t\in W$, we establish
\begin{equation*}
s\in X\text{ and }s<_Xt\quad\Rightarrow\quad s\in W
\end{equation*}
by (side) induction on~$s$ in the order~$\tl$. For~$r\tl s<t$ we get~$r<t$, so that the induction hypothesis yields~$r\in W$. This shows that we have ${\supp_X}\circ\pi(s)\subseteq\rng(\iota)$. We can thus write $\pi(s)=D(\iota)(\sigma)$, due to the support condition. As above, the condition on the range of~$\pi$ entails $G_W(\sigma)\subseteq_{D(W)}\sigma$. Since $s<t$ implies $\sigma<\tau$, the main induction hypothesis yields~$\sigma=\kappa(s')$ for some~$s'\in W$. In view of
\begin{equation*}
\pi(s)=D(\iota)(\sigma)=D(\iota)\circ\kappa(s')=\pi\circ\iota(s')
\end{equation*}
we get $s=\iota(s')\in W$, as needed to complete the side induction step.
\end{remark}

The previous remark is loosely inspired by~\cite[Section~10]{rathjen-weiermann-kruskal}. Similarly, the following generalization to~$\nu>1$ can be seen as an `abstract' version of~\cite[Section~12]{thomson-rathjen-Pi-1-1}. For all result up to Theorem~\ref{thm:Pi11-rec-to-wf}, we fix a $\nu$-collapse $\pi:X\to\nu\times D(X)$ of a dilator~$D$ (note that $D$ preserves well foundedness).

\begin{definition}
For each $\alpha<\nu$ we put
\begin{equation*}
X_\alpha:=\{t\in X\,|\,\pi(t)=(\gamma,\tau)\text{ with }\gamma\leq\alpha\}.
\end{equation*}
Furthermore, we define $E^D_\alpha:X\to[X_\alpha]^{<\omega}$ and $E_\alpha:D(X)\to[X_\alpha]^{<\omega}$ by
\begin{align*}
E^D_\alpha(t)&:=\begin{cases}
\{t\} & \text{if $t\in X_\alpha$},\\
E_\alpha(\tau) & \text{if $\pi(t)=(\gamma,\tau)$ with $\gamma>\alpha$},
\end{cases}\\
E_\alpha(\tau)&:=\bigcup\{E^D_\alpha(s)\,|\,s\in\supp_X(\tau)\}.
\end{align*}
This amounts to a recursion along the well founded relation~$\tl$ from Definition~\ref{def:nu-collapse}.
\end{definition}

Note that each set~$X_\alpha$ is an initial segment of~$X$, since~$\pi$ is an embedding.

\begin{definition}\label{def:W-alpha}
By $\Pi^1_1$-recursion on~$\alpha<\nu$, define $W_\alpha$ as the well founded part of
\begin{equation*}
M_\alpha:=\{t\in X_\alpha\,|\,E^D_\gamma(t)\subseteq W_\gamma\text{ for all }\gamma<\alpha\}.
\end{equation*}
Let us also set $W:=\bigcup\{W_\alpha\,|\,\alpha<\nu\}$.
\end{definition}

We point out that the sets $W_\alpha$ are distinguished (`ausgezeichnet') in the sense of Buchholz~\cite{buchholz-distinguished}, modulo the fact that we are in a somewhat more abstract setting.

\begin{lemma}
For $\alpha\leq\beta$ we have $W_\alpha=W_\beta\cap X_\alpha=W\cap X_\alpha$.
\end{lemma}
\begin{proof}
For $\alpha<\beta$ and $t\in W_\beta\cap X_\alpha$ we get $t\in E^D_\alpha(t)\subseteq W_\alpha$ by the definition of~$M_\beta$. To establish $W_\alpha\subseteq W_\beta$, we argue by induction on~$\beta$. For $\alpha\leq\gamma<\beta$, the induction hypothesis ensures that $t\in W_\alpha\subseteq X_\gamma$ entails $E^D_\gamma(t)=\{t\}\subseteq W_\gamma$, so that we get
\begin{equation*}
W_\alpha\subseteq M_\beta\cap X_\alpha\subseteq M_\alpha.
\end{equation*}
By definition of the well founded part, $W_\beta$ is the largest initial segment of~$M_\beta$ that is well founded. The given inclusions entail that $W_\alpha$ is such a segment and hence contained in~$W_\beta$. More explicitly, induction on~$t\in W_\alpha$ yields~$t\in W_\beta$.
\end{proof}

As~$W$ is the union of well founded initial segments, we get the following.

\begin{corollary}\label{cor:W-wf}
The suborder $W\subseteq X$ is well founded.
\end{corollary}

In the next lemma, we collect some basic facts for later use.

\begin{lemma}\label{lem:E-basic}
The following holds for any $\alpha,\beta<\nu$, any $s,t\in X$ and any $\tau\in D(X)$:
\begin{enumerate}[label=(\alph*)]
\item Given $s\in E^D_\beta(t)$ and $\alpha\leq\beta$, we get $E^D_\alpha(s)\subseteq E^D_\alpha(t)$. The same holds when~$E^D_\beta(t)$ is replaced by~$E_\beta(\tau)$.
\item If $\pi(t)=(\alpha,\tau)$, then we have $E_\alpha(\tau)\subseteq_X t$.
\item From $(\alpha,\tau)\in\rng(\pi)$ we get $(\beta,\tau)\in\rng(\pi)$ for any $\beta\geq\alpha$.
\end{enumerate}
\end{lemma}
\begin{proof}
(a) We argue by induction on~$t$ in the order~$\tl$. For $s=t$, the claim is trivial. In the remaining case, we have $\pi(t)=(\delta,\tau)$ with $\delta>\beta\geq\alpha$. We get $s\in E^D_\beta(r)$ for some $r\tl t$, so that the induction hypothesis yields
\begin{equation*}
E^D_\alpha(s)\subseteq E^D_\alpha(r)\subseteq E_\alpha(\tau)=E^D_\alpha(t).
\end{equation*}
(b) By induction on~$s$ in the order~$\tl$, we prove the auxiliary claim
\begin{equation*}
r\in E^D_\alpha(s)\text{ and }\pi(r)=(\alpha,\rho)\quad\Rightarrow\quad\rho\in G^D_\alpha(s).
\end{equation*}
Assuming the antecedent, we must have $\pi(s)=(\gamma,\sigma)$ with $\gamma\geq\alpha$, so that
\begin{equation*}
G^D_\alpha(s)=\{\sigma\}\cup\bigcup\{G^D_\alpha(s')\,|\,s'\in\supp_X(\sigma)\}.
\end{equation*}
For $r=s$ we obtain $\rho=\sigma\in G^D_\alpha(s)$. In the remaining case we have $r\in E^D_\alpha(s')$ for some $s'\tl s$, so that the induction hypothesis yields $\rho\in G^D_\alpha(s')\subseteq G^D_\alpha(s)$. To deduce the lemma, consider an arbitrary $r\in E_\alpha(\tau)$. Write $\pi(r)=(\delta,\rho)$, necessarily with~$\delta\leq\alpha$. If we have $\delta<\alpha$, then we immediately get $\pi(r)<\pi(t)$ and hence~$r<t$. Now assume~$\delta=\alpha$, and note that we have $r\in E^D_\alpha(s)$ for some $s\in\supp_X(\tau)$. By the auxiliary claim and the condition on~$\rng(\pi)$ in Definition~\ref{def:nu-collapse}, we get
\begin{equation*}
\rho\in G^D_\alpha(s)\subseteq G_\alpha(\tau)\subseteq_{D(X)}\tau.
\end{equation*}
Once again this yields $\pi(r)<\pi(t)$ and hence~$r<t$, as required for $E_\alpha(\tau)\subseteq_X t$.

(c) Given~$\alpha\leq\beta$, one checks $G^D_\beta(t)\subseteq G^D_\alpha(t)$ by a straightforward induction on~$t$ in the order~$\tl$. The same inclusion then holds with $G_\gamma(\tau)$ at the place of~$G^D_\gamma(t)$. Now it suffices to recall the condition on~$\rng(\pi)$ from Definition~\ref{def:nu-collapse}.
\end{proof}

Inspired by~\cite[Definition~12.64]{thomson-rathjen-Pi-1-1}, we introduce the following crucial sets.

\begin{definition}
Let us put
\begin{align*}
B&:=\{\tau\in D(X)\,|\,\text{we have $t\in W$ whenever~$\pi(t)=(\alpha,\tau)$ for some~$\alpha<\nu$}\},\\
M&:=\{\tau\in D(X)\,|\,\text{we have $E_\gamma(\tau)\subseteq W$ for all~$\gamma<\nu$}\}.
\end{align*}
\end{definition}

All of the following results rely on the standing assumption that~$D$ is a dilator. Note that we only use this assumption once, namely in the following proof.

\begin{lemma}
The suborder~$M\subseteq D(X)$ is well founded.
\end{lemma}
\begin{proof}
Given any $\tau\in M$, pick a $\gamma<\nu$ such that the finite set~$\supp_X(\tau)$ is fully contained in~$X_\gamma$. By the definition of~$M$, we obtain
\begin{equation*}
W\supseteq E_\gamma(\tau)=\bigcup\{E^D_\gamma(s)\,|\,s\in\supp_X(\tau)\}=\supp_X(\tau).
\end{equation*}
For the inclusion $\iota:W\hookrightarrow X$, we get $\tau\in\rng(D(\iota))$ by the support condition from Definition~\ref{def:dilator}. Hence~$M$ lies in the range of the embedding~$D(\iota):D(W)\to D(X)$. To conclude, note that~$D(W)$ is well founded as~$D$ is a dilator.
\end{proof}

The next result is the technical core of this section.

\begin{proposition}\label{prop:M-in-B}
We have $M\subseteq B$.
\end{proposition}
\begin{proof}
We argue by (main) induction over the well order~$M$, i.\,e., we assume $\tau\in M$ and $\{\sigma\in M\,|\,\sigma<\tau\}\subseteq B$ to derive~$\tau\in B$. Aiming at the latter, consider an arbitrary~$t\in X$ such that $\pi(t)=(\alpha,\tau)$ holds for some~$\alpha$. We need to prove~$t\in W$. Given $\tau\in M$, we get $t\in M_\alpha$ via
\begin{equation*}
E^D_\gamma(t)=E_\gamma(\tau)\subseteq W\cap X_\gamma=W_\gamma\quad\text{for}\quad\gamma<\alpha.
\end{equation*}
Since $W_\alpha$ is the accessible part of~$M_\alpha$, we can conclude~$t\in W_\alpha\subseteq W$ once the following is established (cf.~\cite[Lemma~12.65]{thomson-rathjen-Pi-1-1}):
\begin{claim}
Given any $\alpha<\nu$ and $t\in X$ with $\pi(t)=(\alpha,\tau)$, we obtain $s\in W_\alpha$ for all elements $s\in M_\alpha$ with $s<t$.
\end{claim}
To prove this claim, we argue by (side) induction on~$s$ in the transitive closure of the well founded relation~$\tl$, or alternatively on~$h(s)$ for
\begin{equation*}
h:X\to\mathbb N\quad\text{with}\quad h(s):=\max(\{0\}\cup\{h(r)+1\,|\,r\tl s\}).
\end{equation*}
It will be important that the induction hypothesis is available for all $\alpha$ and hence for various~$t$, while $\tau$ remains fixed as above. In the side induction step, we first assume that $s\in X_\gamma$ holds for some $\gamma<\alpha$. Given $s\in M_\alpha$, we then get
\begin{equation*}
s\in E^D_\gamma(s)\subseteq W_\gamma\subseteq W_\alpha.
\end{equation*}
In the remaining case, we have $\pi(s)=(\alpha,\sigma)$ with $\sigma<\tau$, as~$s<t$ entails $\pi(s)<\pi(t)$. To use the main induction hypothesis, we want to show~$\sigma\in M$, which amounts to
\begin{equation*}
E_\gamma(\sigma)\subseteq W\quad\text{for all }\gamma<\nu.
\end{equation*}
We prove the latter by (auxiliary) induction on~$\gamma$. For $\gamma<\alpha$ we can invoke~$s\in M_\alpha$ to get $E_\gamma(\sigma)=E^D_\gamma(s)\subseteq W_\gamma$. In the case of~$\gamma=\alpha$, we use Lemma~\ref{lem:E-basic}(a) to obtain
\begin{equation*}
E^D_\delta(r)\subseteq E_\delta(\sigma)\subseteq W\cap X_\delta=W_\delta\quad\text{for all }r\in E_\alpha(\sigma)\text{ and }\delta<\alpha,
\end{equation*}
which yields $E_\alpha(\sigma)\subseteq M_\alpha$. By Lemma~\ref{lem:E-basic}(b), we have $E_\alpha(\sigma)\subseteq_X s<t$. Furthermore, it is not hard to see that the elements of~$E_\alpha(\sigma)$ lie below~$s$ in the transitive closure of~$\tl$ (alternatively check $h(r')\leq h(r)$ for $r'\in E_\alpha^D(r)$ by induction over~$\tl$). We can thus use the side induction hypothesis to get $E_\alpha(\sigma)\subseteq W_\alpha$. Finally, we consider the case of $\gamma>\alpha$. The auxiliary induction hypothesis entails $E_\gamma(\sigma)\subseteq M_\gamma$ as before. By Lemma~\ref{lem:E-basic}(c) we find $s',t'\in X$ with $\pi(s')=(\gamma,\sigma)$ and $\pi(t')=(\gamma,\tau)$. In view of~$\sigma<\tau$ we get
\begin{equation*}
E_\gamma(\sigma)\subseteq_X s'<t'.
\end{equation*}
Thus the desired inclusion $E_\gamma(\sigma)\subseteq W_\gamma$ follows from the side induction hypothesis (now with $\gamma$ and $t'$ at the place of~$\alpha$ and~$t$). This completes the auxiliary induction and hence the proof of~$\sigma\in M$, as noted above. We can now invoke the main induction hypothesis to get~$\sigma\in B$. Given $\pi(s)=(\alpha,\sigma)$, this yields \mbox{$s\in W\cap X_\alpha=W_\alpha$}, which concludes the steps of side induction (claim) and main induction.
\end{proof}

In Remark~\ref{rmk:nu-one}, we have exploited the fact that $t\in W$ and $s\tl t$ entail $s\in W$. The proof that we have given breaks down for~$\nu>1$. However, we get the desired closure property for an inductively generated suborder:

\begin{definition}\label{def:V-sub-W}
Let $V\subseteq W$ be given by the recursive clause
\begin{equation*}
t\in V\quad:\Leftrightarrow\quad t\in W\text{ and }s\in V\text{ for all }s\tl t.
\end{equation*}
\end{definition}

In the following result, the implication $\Rightarrow$ is the closure property mentioned above. The converse implication encapsulates most previous work of this section.

\begin{corollary}\label{cor:V-pi-closed}
For $t\in X$ with $\pi(t)=(\alpha,\tau)$ we have
\begin{equation*}
t\in V\quad\Leftrightarrow\quad\supp_X(\tau)\subseteq V.
\end{equation*}
\end{corollary}
\begin{proof}
Since $s\tl t$ amounts to~$s\in\supp_X(\tau)$, it suffices to show that $\supp_X(\tau)\subseteq V$ implies $t\in W$. For~$\gamma<\nu$, a straightforward induction over~$\tl$ shows that $s\in V$ entails $E^D_\gamma(s)\subseteq V$. Given $\supp_X(\tau)\subseteq V$, we thus get
\begin{equation*}
E_\gamma(\tau)=\bigcup\{E^D_\gamma(s)\,|\,s\in\supp_X(\tau)\}\subseteq V\subseteq W.
\end{equation*}
This shows~$\tau\in M$, so that Proposition~\ref{prop:M-in-B} yields~$\tau\in B$, which entails~$t\in W$.
\end{proof}

Finally, we deduce the main result of this section, which shows that (i) implies~(ii) in Theorem~\ref{thm:main}. To justify the formulation of the following theorem, we recall that $\nu$-fixed points exist and are essentially unique, by Theorem~\ref{thm:psi-is-fp} and Corollary~\ref{cor:nu-fp-unique}.

\begin{theorem}\label{thm:Pi11-rec-to-wf}
If $\Pi^1_1$-recursion along~$\nu$ is available, then the $\nu$-fixed point of any dilator is well founded.
\end{theorem}
\begin{proof}
Consider a dilator~$D$ and a $\nu$-fixed point~$X$ with collapse~$\pi:X\to\nu\times D(X)$. Using $\Pi^1_1$-recursion along~$\nu$, we can construct sets $W_\alpha$ as in Definition~\ref{def:W-alpha}, to obtain suborders $V\subseteq W\subseteq X$ as in Definition~\ref{def:V-sub-W}. Note that $V$ is well founded by Corollary~\ref{cor:W-wf}. We shall show that $V$ is a $\nu$-fixed point of~$D$. Once this is achieved, we can use Corollary~\ref{cor:nu-fp-unique} to conclude that $X\cong V$ is well founded. In fact, we could derive $X=V$ via Proposition~\ref{prop:FP-initial} (as in Remark~\ref{rmk:nu-one}). Write $\iota:V\to X$ for the inclusion. By the previous corollary and the support condition from Definition~\ref{def:dilator}, we get $\tau\in\rng(D(\iota))$ whenever we have $\pi(t)=(\alpha,\tau)$ with $t\in V$. We thus obtain an embedding $\kappa$ so that
\begin{equation*}
\begin{tikzcd}
V\arrow[r,"\kappa"]\arrow[d,swap,"\iota"] & \nu\times D(V)\arrow[d,"I\times D(\iota)"]\\
X\arrow[r,"\pi"] & \nu\times D(X)
\end{tikzcd}
\end{equation*}
commutes. Concerning the constructions from Definition~\ref{def:nu-collapse}, we note that $\kappa$ and~$\pi$ induce the same relation~$\tl$ on~$V\subseteq X$, as in Remark~\ref{rmk:nu-one}. The cited definition also yields functions $G^{D,Z}_\gamma:Z\to[D(Z)]^{<\omega}$ and $G^Z_\gamma:D(Z)\to[D(Z)]^{<\omega}$ for $Z=X$ and for $Z=V$, which are defined with respect to~$\kappa$ and~$\pi$. As in Remark~\ref{rmk:nu-one}, a straightforward induction over~$\tl$ shows
\begin{equation*}
[D(\iota)]^{<\omega}\circ G^{D,V}_\gamma=G^{D,X}_\gamma\circ\iota\quad\text{and}\quad [D(\iota)]^{<\omega}\circ G^V_\gamma=G^X_\gamma\circ D(\iota).
\end{equation*}
It remains to establish the crucial condition from Definition~\ref{def:nu-collapse}, i.\,e., the equation
\begin{equation*}
\rng(\kappa)=\{(\alpha,\tau)\in\nu\times D(V)\,|\,G^V_\alpha(\tau)\subseteq_{D(V)}\tau\}.
\end{equation*}
We point out that the analogous condition is given for~$\pi$, as the latter is a $\nu$-collapse. As in Remark~\ref{rmk:nu-one}, one derives the inclusion~$\subseteq$ and shows that $G^V_\alpha(\tau)\subseteq_{D(V)}\tau$ entails $(\alpha,D(\iota)(\tau))=\pi(t)$ for some~$t\in X$. Note that we have
\begin{equation*}
{\supp_X}\circ D(\iota)(\tau)=[\iota]^{<\omega}\circ\supp_V(\tau)\subseteq V,
\end{equation*}
as $\supp:D\Rightarrow[\cdot]^{<\omega}$ is a natural transformation. Crucially, we can now infer~$t\in V$ by the non-trivial direction of Corollary~\ref{cor:V-pi-closed}. In view of
\begin{equation*}
(I\times D(\iota))(\alpha,\tau)=(\alpha,D(\iota)(\tau))=\pi(t)=\pi\circ\iota(t)=(I\times D(\iota))\circ\kappa(t),
\end{equation*}
we get $(\alpha,\tau)=\kappa(t)\in\rng(\kappa)$ as desired.
\end{proof}

\section{Booting up: Bachmann-Howard fixed points and Veblen hierarchy}\label{sect:booting-up}

In the first part of this section, we establish a connection between Bachmann-Howard fixed points and \mbox{$1$-fixed} points (cf.~Definitions~\ref{def:Bachmann-Howard-fp} and~\ref{def:nu-collapse}). This will allow us to use $\Pi^1_1$-comprehension whenever the well foundedness of $1$-fixed points is given, due to Theorem~\ref{thm:Pi^1_1-CA} (proved in~\cite{freund-equivalence,freund-computable}). Amongst others, $\Pi^1_1$-comprehension secures the Veblen hierarchy of normal functions. In the second part of this section, we discuss a functor~$\Gamma$ that represents this hierarchy. It will be used in our proof that~(iii) implies~(iv) in Theorem~\ref{thm:main}.

We begin with the easier part of the connection, which will not be needed in this paper but completes the picture in a satisfactory way:

\begin{proposition}\label{prop:theta-to-psi}
Assume that $Z$ is a Bachmann-Howard fixed point of a given pre\-dilator~$D$. Then some suborder~$X\subseteq Z$ is a $1$-fixed point of~$D$.
\end{proposition}
\begin{proof}
By assumption, we have a Bachmann-Howard collapse $\vartheta:D(Z)\to Z$. To~see that $\vartheta$ is injective, consider an inequality $\sigma<\tau$ in the linear order~$D(Z)$. If we have $\supp_Z(\sigma)\subseteq_{D(Z)}\vartheta(\tau)$, then clause~(i) of Definition~\ref{def:Bachmann-Howard-fp} yields $\vartheta(\sigma)<\vartheta(\tau)$. Otherwise, there is an $r\in\supp_Z(\sigma)$ with $\vartheta(\tau)\leq r<\vartheta(\sigma)$, where the second inequality relies on clause~(ii) of the cited definition. We shall assume that $\vartheta$ is also surjective and that
\begin{equation*}
s\tl\vartheta(\tau)\quad:\Leftrightarrow\quad s\in\supp_Z(\tau)
\end{equation*}
defines a well founded relation on~$Z$. To justify these assumptions, we point out that they hold when~$Z$ is the minimal Bachmann-Howard fixed point~$\vartheta(D)$ that was constructed in~\cite[Section~4]{freund-computable}. In other words, we can replace~$Z$ by $\vartheta(D)\subseteq Z$ to satisfy the additional assumptions. Let us now define $G^D:Z\to[D(Z)]^{<\omega}$ and simultaneously $G:D(Z)\to[D(Z)]^{<\omega}$ by the recursive clauses
\begin{equation*}
G^D(\vartheta(\tau)):=\{\tau\}\cup G(\tau)\quad\text{and}\quad G(\tau):=\bigcup\{G^D(s)\,|\,s\in\supp_Z(\tau)\}.
\end{equation*}
By induction on~$s$ in the order~$\tl$, we can show
\begin{equation*}
G^D(s)\subseteq_{D(Z)}\tau\quad\Rightarrow\quad s<\vartheta(\tau).
\end{equation*}
Indeed, assume that the premise holds for $s=\vartheta(\sigma)$. We then have $\sigma\in G^D(s)$ and hence $\sigma<\tau$. To conclude by clause~(i) of Definition~\ref{def:Bachmann-Howard-fp}, we note that $r\in\supp_Z(\sigma)$ entails $G^D(r)\subseteq G^D(s)$, so that $r<\vartheta(\tau)$ follows by induction hypothesis. Now set
\begin{equation*}
Y:=\{s\in Z\,|\,s=\vartheta(\sigma)\text{ with }G(\sigma)\subseteq_{D(Z)}\sigma\}.
\end{equation*}
To generate $X\subseteq Y$, we inductively declare
\begin{equation*}
t\in X\quad:\Leftrightarrow\quad t\in Y\text{ and }s\in X\text{ for all }s\tl t.
\end{equation*}
Write $\iota:X\hookrightarrow Z$ for the inclusion. For $\vartheta(\tau)\in X$ we get $\supp_Z(\tau)\subseteq X=\rng(\iota)$. Hence we have $\tau=D(\iota)(\sigma)$ for a (necessarily unique) element~$\sigma\in D(X)$, by the support condition from Definition~\ref{def:dilator}. We thus find a function~$\pi$ such that
\begin{equation*}
\begin{tikzcd}
X\arrow[r,"\pi"]\arrow[d,swap,hook,"\iota"] & D(X)\arrow[d,"D(\iota)"]\\
Z & D(Z)\arrow[l,swap,"\vartheta"]
\end{tikzcd}
\end{equation*}
is a commutative diagram. Clearly $\pi$ is injective. To conclude that it is an order embedding, we assume $\pi(s)<\pi(t)$ and deduce~$s<t$. Given $s\in X$, we get
\begin{equation*}
\vartheta\circ D(\iota)\circ\pi(s)=\iota(s)\in X\subseteq Y.
\end{equation*}
By the definition of~$Y$, this yields $G(D(\iota)\circ\pi(s))\subseteq_{D(Z)}D(\iota)\circ\pi(s)$ and hence
\begin{equation*}
G^D(\vartheta\circ D(\iota)\circ\pi(s))=\{D(\iota)\circ\pi(s)\}\cup G(D(\iota)\circ\pi(s))\subseteq_{D(Z)}D(\iota)\circ\pi(t).
\end{equation*}
Due to the implication that was shown above, one can infer $s<t$ via
\begin{equation*}
\iota(s)=\vartheta\circ D(\iota)\circ\pi(s)<\vartheta\circ D(\iota)\circ\pi(t)=\iota(t).
\end{equation*}
After some straightforward verifications, we can conclude that $\pi$ is a $1$-collapse of the predilator $D$ (where we identify $D(X)$ and $1\times D(X)$ as in Remark~\ref{rmk:nu-one}).
\end{proof}

Let $D$ and $E$ be predilators with associated transformations $\supp^D:D\Rightarrow[\cdot]^{<\omega}$ and $\supp^E:E\Rightarrow[\cdot]^{<\omega}$. The predilator $E\circ D$ consists of the usual composition as functors and the transformation $\supp^{E\circ D}:E\circ D\Rightarrow[\cdot]^{<\omega}$ that is given by
\begin{equation*}
\supp^{E\circ D}_X(\sigma):=\bigcup\{\supp^D_X(\rho)\,|\,\rho\in\supp^E_{D(X)}(\sigma)\}.
\end{equation*}
It is straightforward to check that the conditions from Definition~\ref{def:dilator} are satisfied. In the following theorem, we write $\omega$ for the predilator from Example~\ref{ex:Buchholz-psi} (see also the beginning of Section~\ref{sect:intro}). The result is an abstract version of~\cite[Corollary~3.1]{rathjen-weiermann-kruskal}, which provides a similar connection between concrete ordinal notation systems.

\begin{theorem}\label{thm:1fp-to-Bachmann-Howard}
Any $1$-fixed point of $\omega\circ D$ is a Bachmann-Howard fixed point of~$D$, where $D$ can be any predilator.
\end{theorem}
\begin{proof}
Consider a $1$-collapse $\pi:X\to\omega\circ D(X)=:E(X)$, where we identify $E(X)$ and $1\times E(X)$ as before. Let $G^E_0:X\to[E(X)]^{<\omega}$ and $G_0:E(X)\to[E(X)]^{<\omega}$ be given as in Definition~\ref{def:nu-collapse} (see also Remark~\ref{rmk:nu-one}), so that we have
\begin{equation*}
\rng(\pi):=\{\tau\in\omega\circ D(X)\,|\,G_0(\tau)\subseteq_{\omega\circ D(X)}\tau\}.
\end{equation*}
We need to define a function $\vartheta:D(X)\to X$ that satisfies clauses~(i) and~(ii) from Definition \ref{def:Bachmann-Howard-fp}. As in the first paragraph of Section~\ref{sect:intro}, we write elements of $\omega\circ D(X)$ in the form $\langle\sigma_0,\ldots,\sigma_{n-1}\rangle$, for elements $\sigma_0\geq\ldots\geq\sigma_{n-1}$ of~$D(X)$. In particular, a given~$\sigma\in D(X)$ gives rise to an element $\langle\sigma\rangle\in\omega\circ D(X)$, which allows us to form
\begin{equation*}
\sigma^\star:=\max\big(\{\langle\rangle\}\cup G_0(\langle\sigma\rangle)\big)\in\omega\circ D(X).
\end{equation*}
Writing $\sigma^\star=\langle\sigma_0,\ldots,\sigma_{n-1}\rangle$, we now set
\begin{equation*}
\sigma^+:=\langle\sigma_0,\ldots,\sigma_{i(\sigma)-1},\sigma\rangle\quad\text{with}\quad i(\sigma):=\min\big(\{i<n\,|\,\sigma_i<\sigma\}\cup\{n\}\big).
\end{equation*}
Note that we have $\sigma^+\in\omega\circ D(X)$, as the definition of $i(\sigma)$ ensures $\sigma_{i(\sigma)-1}\geq\sigma$. Informally, we point out that the given construction corresponds to $\sigma^+=\sigma^\star+\omega^\sigma$ in terms of ordinal arithmetic. Let us now show
\begin{equation*}
G_0(\sigma^+)\subseteq G_0(\sigma^\star)\cup G_0(\langle\sigma\rangle)\subseteq G_0(\langle\sigma\rangle)\subseteq_{\omega\circ D(X)}\sigma^+.
\end{equation*}
The first inclusion reduces to the analogous inclusions for $\supp^{\omega\circ D}_X$ and $\supp^\omega_{D(X)}$, which we get by the definition of supports in Example~\ref{ex:Buchholz-psi}. Concerning the second inclusion, we note that $G_0(\langle\rangle)$ is empty, since the same holds for $\supp^\omega_{D(X)}(\langle\rangle)$ and hence for $\supp^{\omega\circ D}_X(\langle\rangle)$. In the remaining case we have $\sigma^\star\in G_0(\langle\sigma\rangle)$. Here we can infer $G_0(\sigma^\star)\subseteq G_0(\langle\sigma\rangle)$ from the general fact that $\rho\in G^E_0(s)$ entails $G_0(\rho)\subseteq G_0^E(s)$, which is readily verified by induction on~$s$ in the order~$\tl$ from Definition~\ref{def:nu-collapse}. Finally, we see that $r\in G_0(\langle\sigma\rangle)$ entails $r\leq\sigma^\star<\sigma^+$, by the definition of~$\sigma^\star$ and as we have $\sigma_{i(\sigma)}<\sigma$ or $i(\sigma)=n$ (recall that $\omega\circ D(X)$ is ordered lexicographically). For any $\sigma\in D(X)$, we have shown $G_0(\sigma^+)\subseteq_{\omega\circ D(X)}\sigma^+$, which entails $\sigma^+\in\rng(\pi)$. This allows us to form the function
\begin{equation*}
\vartheta:D(X)\to X\quad\text{with}\quad\pi\circ\vartheta(\sigma)=\sigma^+,
\end{equation*}
which is unique since~$\pi$ is an embedding. To verify clause~(ii) of Definition~\ref{def:Bachmann-Howard-fp}, we show $r<\vartheta(\sigma)$ for a given $r$ in the set $\supp^D_X(\sigma)$. The latter is equal to $\supp^{\omega\circ D}_X(\langle\sigma\rangle)$, as we have $\supp^\omega_X(\langle\sigma\rangle)=\{\sigma\}$. We thus get
\begin{equation*}
\pi(r)\in G^{\omega\circ D}_0(r)\subseteq G_0(\langle\sigma\rangle)\quad\text{and hence}\quad\pi(r)\leq\sigma^\star<\sigma^+=\pi\circ\vartheta(\sigma),
\end{equation*}
which yields $r<\vartheta(\sigma)$ as desired. In order to prepare the remaining verification, we recall that $s\tl t$ entails $\pi(s)<\pi(t)$, as observed in Remark~\ref{rmk:nu-one}. One can derive that $\rho\in G^{\omega\circ D}_0(t)$ entails $\rho\leq\pi(t)$, by a straightforward induction on~$t$ in the order~$\tl$. Aiming at clause~(i) of Definition~\ref{def:Bachmann-Howard-fp}, we now assume
\begin{equation*}
\sigma<_{D(X)}\tau\quad\text{and}\quad\supp^D_X(\sigma)\subseteq_X\vartheta(\tau).
\end{equation*}
For an arbitrary $r\in\supp^D_X(\sigma)$ and any $\rho\in G^{\omega\circ D}_0(r)$, we get
\begin{equation*}
\rho\leq\pi(r)<\pi\circ\vartheta(\tau)=\tau^+.
\end{equation*}
In view of $\supp^{\omega\circ D}_X(\langle\sigma\rangle)=\supp^D_X(\sigma)$ from above, this yields
\begin{equation*}
G_0(\langle\sigma\rangle)=\bigcup\left\{\left.G^{\omega\circ D}_0(r)\,\right|\,r\in\supp^{\omega\circ D}_X(\langle\sigma\rangle)\right\}\subseteq_{\omega\circ D(X)}\tau^+.
\end{equation*}
Together with $0\leq\tau^\star<\tau^+$, we get $\sigma^\star<\tau^+$. The latter and $\sigma<\tau$ entail
\begin{equation*}
\pi\circ\vartheta(\sigma)=\sigma^+<\tau^+=\pi\circ\vartheta(\tau)
\end{equation*}
and hence $\vartheta(\sigma)<\vartheta(\tau)$, by basic considerations about the lexicographic order.
\end{proof}

As noted at the beginning of Section~\ref{sect:intro}, the statement that ``$\omega(X)$ is well founded for any well order~$X$" is equivalent to arithmetical comprehension and hence unprovable in the theory~$\rca_0$. The latter can prove that $\omega$ is a predilator but not that it is a dilator. To prepare the use of Theorem~\ref{thm:1fp-to-Bachmann-Howard} over~$\rca_0$, we show the following proposition. It is interesting to compare the result with~\cite[Theorem~2.2]{freund-predicative-collapsing}.

\begin{proposition}\label{prop:omega-1fp}
For any linear order~$Y$, the order $\omega(Y)$ is a $1$-fixed point of a predilator~$D$ with $D(X)=1+Y\times X$ (see the proof for a detailed definition of~$D$).
\end{proposition}
\begin{proof}
Recall the notation for products from the paragraph before Definition~\ref{def:nu-collapse}. To complete the definition of~$1+Y\times X$, we introduce general notation for the sum of linear orders~$Z_0$ and $Z_1$, which will also be needed later. The underlying set of our sum is the disjoint union
\begin{equation*}
Z_0+Z_1:=\{z_0\,|\,z_0\in Z_0\}\cup\{Z_0+z_1\,|\,z_1\in Z_1\}.
\end{equation*}
To determine the order, we declare that $z_0\mapsto z_0$ and $z_1\mapsto Z_0+z_1$ are embeddings of $Z_0$ and $Z_1$ into $Z_0+Z_1$, while $z_0<Z_0+z_1$ holds for any $z_i\in Z_i$. Given embeddings $f_i:Z_i\to Z_i'$, we define $f_0+f_1:Z_0+Z_1\to Z'_0+Z'_1$ by
\begin{equation*}
(f_0+f_1)(z_0):=f_0(z_0)\quad\text{and}\quad (f_0+f_1)(Z_0+z_1):=Z'_0+f_1(z_1).
\end{equation*}
If $f_0$ or $f_1$ is the identity on $Z_0=Z'_0$ or $Z_1=Z'_1$, respectively, we write $Z_0+f_1$ or $f_0+Z_1$ rather than $f_0+f_1$. Let us agree that $\times$ binds stronger than $+$ and that $1=\{0\}$ denotes the singleton order. For our fixed order~$Y$, this explains the transformations $X\mapsto D(X):=1+Y\times X$ and $f\mapsto D(f):=1+Y\times f$ of orders and embeddings. To turn $D$ into a dilator, we define $\supp_X:D(X)\to[X]^{<\omega}$ by
\begin{equation*}
\supp_X(0):=\emptyset\quad\text{and}\quad\supp_X(1+(y,x)):=\{x\}.
\end{equation*}
Let us now consider the embedding $\pi:\omega(Y)\to 1+Y\times\omega(Y)$ with
\begin{equation*}
\pi(\langle\rangle):=0\quad\text{and}\quad\pi(\langle y_0,\ldots,y_n\rangle):=1+(y_0,\langle y_1,\ldots,y_n\rangle).
\end{equation*}
To see that $\pi$ is a $1$-collapse of~$D$, we need to show
\begin{equation*}
\rng(\pi)=\{\tau\in D\circ\omega(Y)\,|\,G_0(\tau)\subseteq_{D\circ\omega(Y)}\tau\},
\end{equation*}
with $G_0:D\circ\omega(Y)\to[D\circ \omega(Y)]^{<\omega}$ as in Definition~\ref{def:nu-collapse} (see also Remark~\ref{rmk:nu-one}). First note that we have $0\in\rng(\pi)$ while $\supp_X(0)$ and hence $G_0(0)$ is empty. Let us now consider $\tau=1+(y_0,\langle y_1,\ldots,y_n\rangle)$. We then have $G_0(\tau)=G^D_0(\langle y_1,\ldots,y_n\rangle)$, where $G^D_0:\omega(Y)\to[D\circ\omega(Y)]^{<\omega}$ is recursively given by $G^D_0(\langle\rangle)=\{0\}$ and
\begin{equation*}
G^D_0(\langle z_0,\ldots,z_m\rangle)=\{1+(z_0,\langle z_1,\ldots,z_m\rangle)\}\cup G^D_0(\langle z_1,\ldots,z_m\rangle).
\end{equation*}
Let us observe that $1+(z_0,\langle z_1,\ldots,z_m\rangle)$ is the largest element of this set, by a straightforward induction on~$m$ (note $z_1\leq z_0$ and $\langle z_2,\ldots,z_m\rangle<\langle z_1,\ldots,z_m\rangle$). If we have $n=0$ and hence $\tau=1+(y_0,\langle\rangle)$, then we get $G_0(\tau)=\{0\}\subseteq_{D\circ\omega(Y)}\tau$ as well as $\tau=\pi(\langle y_0\rangle)\in\rng(\pi)$. In the case of $n>0$, we need to show
\begin{equation*}
\tau=1+(y_0,\langle y_1,\ldots,y_n\rangle)\in\rng(\pi)\quad\Leftrightarrow\quad 1+(y_1,\langle y_2,\ldots,y_n\rangle)<\tau.
\end{equation*}
Given $\langle y_1,\ldots,y_n\rangle\in\omega(Y)$, we see that both sides are equivalent to $y_1\leq y_0$.
\end{proof}

Based on Theorem~\ref{thm:Pi^1_1-CA}, we can now derive that the equivalence from Theorem~\ref{thm:main} holds for~$\nu=1$. This allows us to use $\Pi^1_1$-comprehension whenever the well foundedness of $\nu$-fixed points is given. In view of Proposition~\ref{prop:theta-to-psi}, the following can be seen as a strengthening of Theorem~\ref{thm:Pi^1_1-CA}.

\begin{corollary}\label{cor:1-fixed-points}
For each fixed~$\nu\in\mathbb N\backslash\{0\}$, the following are equivalent over~$\rca_0$:
\begin{enumerate}[label=(\roman*)]
\item $\Pi^1_1$-comprehension,
\item the $\nu$-fixed point of any dilator is well founded,
\item any dilator has a well founded $\mu$-fixed point for some well order~$\mu\neq\emptyset$.
\end{enumerate}
\end{corollary}
\begin{proof}
By iterated applications of~(i), we obtain $\Pi^1_1$-recursion along~$\nu$, as the latter is fixed externally. We can then invoke Theorem~\ref{thm:Pi11-rec-to-wf} to get~(ii), which does clearly imply~(iii). Assuming the latter, we argue that any given dilator~$D$ has a well founded Bachmann-Howard fixed point, to infer~(i) via Theorem~\ref{thm:Pi^1_1-CA}. In any application of~(iii) we may assume $\mu=1$, due to Corollary~\ref{cor:no-fp-monotone}. If~$Y$ is a well order, then the predilator from the previous proposition is a dilator, provably in~$\rca_0$. In the presence of~(iii), we can conclude that~$\omega(Y)$ is well founded. So we know that $\omega\circ D$ is a dilator. Using~(iii) again, we get a well founded $1$-fixed point of~$\omega\circ D$. By Theorem~\ref{thm:1fp-to-Bachmann-Howard}, this is the desired Bachmann-Howard fixed point of~$D$.
\end{proof}

In the rest of this section, we discuss a dilator~$\Gamma$ such that $\Gamma(X)$ represents the Veblen function~$\varphi$ up to the $X$-th ordinal~$\alpha$ with $\varphi(\alpha,0)=\alpha$ (such $\alpha$ are called `strongly critical'). The Veblen function plays an important role in ordinal analysis (see e.\,g.~\cite[Chapters~V and~VII]{schuette77}) and can also be analysed in terms of computability theory (as done by Marcone and Montalb\'an~\cite{marcone-montalban}). We will use the dilator~$\Gamma$ in our proof that~(iii) implies~(iv) in Theorem~\ref{thm:main}, where we mimic traditional ordinal analysis in a more abstract setting. To understand the following, it is not indispensable but certainly helpful to know the set theoretic approach to the Veblen function, for which we refer to~\cite[Section~3]{pohlers-proof-theory}.

The next definition is equivalent to~\cite[Definition~2.5]{rathjen-atr}, despite a small difference in clause~(ii'). A detailed justification of the recursion is given after the definition. The abbreviations $\scr$ and $\hau$ stand for `strongly critical' ordinals and `Haupt\-zahlen'. The latter is German for (additively) `principal numbers'. We write $\overline\varphi$ in order to save the symbol $\varphi$ for Definition~\ref{def:total-Veblen} below.

\begin{definition}\label{def:Gamma(X)}
Given a linear order~$X$, we define sets $\scr\subseteq\hau\subseteq\Gamma(X)$ of terms, a binary relation $<_{\Gamma(X)}$ on~$\Gamma(X)$ and a function $h:\Gamma(X)\to\Gamma(X)$ by simultaneous recursion. We write $s\leq_{\Gamma(X)}t$ to abbreviate the disjunction of $s<_{\Gamma(X)}t$ and $s=t$, where the latter denotes equality as terms. The terms are generated as follows:
\begin{enumerate}[label=(\roman*)]
\item We have terms $0\in\Gamma(X)\backslash\hau$ and $\Gamma_x\in\scr\subseteq\hau\subseteq\Gamma(X)$ for all $x\in X$.
\item Assume that we are given terms $s,t\in\Gamma(X)$ with $h(t)\leq_{\Gamma(X)}s$, such that we have $t\neq 0$ or $s\notin\scr$. We then add a term~$\overline\varphi st\in\hau\backslash\scr\subseteq\Gamma(X)$.
\item Given $n>1$ terms $t_0,\ldots,t_{n-1}\in\hau$ with $t_{i+1}\leq_{\Gamma(X)} t_i$ for~$i<n-1$, we add a term $\langle t_0,\ldots,t_{n-1}\rangle\in\Gamma(X)\backslash\hau$.
\end{enumerate}
To determine~$h$, we put $h(\Gamma_x):=\Gamma_x$ and $h(\overline\varphi st):=s$ as well as $h(t):=0$ in the remaining cases. Let us abbreviate $\langle\rangle:=0$ and $\langle t\rangle:=t$ for $t\in\hau$, so that any element of $\Gamma(X)$ can be uniquely written in the form~$\langle t_0,\ldots,t_{n-1}\rangle$ with $n\in\mathbb N$. We declare that $<_{\Gamma(X)}$ is the minimal relation with the following closure properties:
\begin{enumerate}[label=(\roman*')]
\item We have $r<_{\Gamma(X)}\Gamma_y$ for $r=\Gamma_x$ with $x<_X y$, for $r=\overline\varphi st$ with $s,t<_{\Gamma(X)}\Gamma_x$, and for $r=\langle r_0,\ldots,r_{n-1}\rangle$ with $n=0$ or $r_0<_{\Gamma(X)}\Gamma_y$.
\item We have $r<_{\Gamma(X)}\overline\varphi st$ for $r=\Gamma_x$ with $r\leq_{\Gamma(X)}s$ or $r\leq_{\Gamma(X)}t$, for a term $r=\langle r_0,\ldots,r_{n-1}\rangle$ with $n=0$ or $r_0<_{\Gamma(X)}\overline\varphi st$, and for $r=\overline\varphi s't'$ such that
\begin{itemize}
\item we have $s'<_{\Gamma(X)}s$ and $t'<_{\Gamma(X)}\overline\varphi st$,
\item or we have $s=s'$ and $t<_{\Gamma(X)}t'$,
\item or we have $\overline\varphi s't'\leq_{\Gamma(X)}t$.
\end{itemize}
\item We get $\langle s_0,\dots,s_{m-1}\rangle<_{\Gamma(X)}\langle t_0,\ldots,t_{n-1}\rangle$, not necessarily with $m,n>1$,~if
\begin{itemize}
\item we have $m<n$ and $s_i=t_i$ for all $i<m$,
\item or there is a $j<\min\{m,n\}$ with $s_j<_{\Gamma(X)}t_j$ and $s_i=t_i$ for all $i<j$.
\end{itemize}
\end{enumerate}
We will sometimes write $<$ rather than~$<_{\Gamma(X)}$ when no ambiguity arises.
\end{definition}

Note that clause~(iii') for $m=0$ yields $0<_{\Gamma(X)}t$ when $t\neq 0$. For $m=1$ we learn that $s<_{\Gamma(X)}\langle t_0,\ldots,t_{n-1}\rangle$ is equivalent to $s\leq_{\Gamma(X)}t_0$ when $s\in\hau$ and $n>1$. The reader may wish to reformulate the clause for $m>1$ and $n=1$ in a similar~way. Also note that $m,n=1$ makes~(iii') tautological, so that no new inequalities arise. Finally, observe that $\Gamma(X)$ is isomorphic to~$\omega(\hau)$, as defined in Section~\ref{sect:intro}.

To justify the simultaneous recursion in Definition~\ref{def:Gamma(X)}, let $\Gamma^+(X)\supseteq\Gamma(X)$ be generated by clauses~(i) to~(iii) but with all conditions that involve $<_{\Gamma(X)}$ ignored. Define $h:\Gamma^+(X)\to\Gamma^+(X)$ as above, and consider $L:\Gamma ^+(X)\to\mathbb N$ with
\begin{gather*}
L(0):=L(\Gamma_x):=0,\qquad L(\overline\varphi st):=L(s)+L(t)+1,\\
L(\langle t_0,\ldots,t_{n-1}\rangle):=L(t_0)+\ldots+L(t_{n-1})+1\quad\text{(for $n>1$)}.
\end{gather*}
Note that $L(h(t))\leq L(t)$ holds for all~$t\in\Gamma^+(X)$. One can now decide $r\in\Gamma(X)$ and $s<_{\Gamma(X)}t$ by simultaneous recursion on~$L(r)$ and~$L(s)+L(t)$, respectively. This decision procedure is implicit in part~(ii) of~\cite[Lemma~2.6]{rathjen-atr}. Part~(i) of the latter coincides with~(b) in the next result, up to the modified formulation of~(ii') above.

\begin{lemma}\label{lem:Gamma-linear}
The following holds for any linear order~$X$:
\begin{enumerate}[label=(\alph*)]
\item We have $s,t<_{\Gamma(X)}\overline\varphi st$ and $t_0<_{\Gamma(X)}\langle t_0,\ldots,t_{n-1}\rangle$ in case $n>1$.
\item The relation~$<_{\Gamma(X)}$ is a linear order on~$\Gamma(X)$.
\end{enumerate}
\end{lemma}
\begin{proof}
First show that $r<s$ and $s<t$ yield $r<t$, by induction on~$L(r)+L(s)+L(t)$ and a lengthy but straightforward case distinction. To establish~(a), we consider the subterm function $\sub:\Gamma(X)\to[\Gamma(X)]^{<\omega}$ with
\begin{gather*}
\sub(0):=\sub(\Gamma_x):=\emptyset,\qquad\sub(\overline\varphi st):=\{s,t\}\cup\sub(s)\cup\sub(t),\\
\sub(\langle t_0,\ldots,t_{n-1}\rangle):=\textstyle\bigcup_{i<n}\big(\{t_i\}\cup\sub(t_i)\big)\quad\text{(for $n>1$)}.
\end{gather*}
A straightforward induction on~$L(s)+L(t)$ shows that $s\in\sub(t)$ entails $s<_{\Gamma(X)}t$. Concerning the case of $t=\langle t_0,\ldots,t_{n-1}\rangle$, note that $s\leq t_i$ yields $s\leq t_0$, as tran\-si\-tivity has already been proved. To show $t\not< t$, one argues by induction on~$L(t)$. The only non-trivial task is to exclude $t=\overline\varphi t_0t_1\leq t_1$. The latter would imply $t_1<t_1$ by~(a) and transitivity, against the induction hypothesis. Finally, a straightforward induction on~$L(s)+L(t)$ shows that we always have $s<t$ or~$s=t$ or $s>t$.
\end{proof}

Concerning the following definition, it is immediate that the range of~$\Gamma(f)$ is contained in~$\Gamma^+(Y)\supseteq\Gamma(Y)$, as defined in the paragraph before Lemma~\ref{lem:Gamma-linear}. In the proof of Proposition~\ref{ref:Gamma-dilator} below, we show that it is indeed contained in~$\Gamma(Y)$.

\begin{definition}\label{def:Gamma-dilator}
For an embedding $f:X\to Y$, we define $\Gamma(f):\Gamma(X)\to\Gamma(Y)$ by
\begin{gather*}
\Gamma(f)(0):=0,\qquad\Gamma(f)(\Gamma_x):=\Gamma_{f(x)},\\
\Gamma(f)(\overline\varphi t_0t_1):=\overline\varphi t_0't_1'\text{ with }t_i':=\Gamma(f)(t_i),\\
\Gamma(f)(\langle t_0,\ldots,t_{n-1}\rangle):=\langle\Gamma(f)(t_0),\ldots,\Gamma(f)(t_{n-1})\rangle\quad\text{(for $n>1$)}.
\end{gather*}
We also define functions $\supp^\Gamma_X:\Gamma(X)\to[X]^{<\omega}$ by stipulating
\begin{gather*}
\supp^\Gamma_X(0):=\emptyset,\quad\supp^\Gamma_X(\Gamma_x):=\{x\},\quad\supp^\Gamma_X(\overline\varphi st):=\supp^\Gamma_X(s)\cup\supp^\Gamma_X(t),\\
\supp^\Gamma_X(\langle t_0,\ldots,t_{n-1}\rangle):=\textstyle\bigcup_{i<n}\supp^\Gamma_X(t_i)\quad\text{(for $n>1$).}
\end{gather*}
\end{definition}

In the following, a stronger metatheory is needed for matters of well foundedness. We rely on $\Pi^1_1$-comprehension, which will be available in our intended application (via Corollary~\ref{cor:1-fixed-points}). The proof shows that a somewhat weaker principle suffices.

\begin{proposition}\label{ref:Gamma-dilator}
The data from Definitions~\ref{def:Gamma(X)} and~\ref{def:Gamma-dilator} constitutes a predilator~$\Gamma$ (provably in $\rca_0$), which is in fact a dilator (in the presence of $\Pi^1_1$-comprehension).
\end{proposition}
\begin{proof}
Given an embedding $f:X\to Y$, let $\Gamma(f):\Gamma^+(X)\to\Gamma^+(Y)$ be defined by the clauses from Definition~\ref{def:Gamma-dilator}, applied to the larger sets $\Gamma^+(Z)\supseteq\Gamma(Z)$ from the paragraph before Lemma~\ref{lem:Gamma-linear}. For $r\in\Gamma^+(X)$ and $s,t\in\Gamma(X)$ one readily shows
\begin{align*}
r\in\Gamma(X)\quad&\Leftrightarrow\quad\Gamma(f)(r)\in\Gamma(Y),\\
s<_{\Gamma(X)}t\quad&\Leftrightarrow\quad\Gamma(f)(s)<_{\Gamma(Y)}\Gamma(f)(t)
\end{align*}
by simultaneous induction on $L(r)$ and $L(s)+L(t)$, respectively. Concerning the first equivalence, we note that $\Gamma(f)$ commutes with the functions $h:\Gamma(Z)\to\Gamma(Z)$ from Definition~\ref{def:Gamma(X)}. To establish the second equivalence, it suffices to show the implication from left to right, which yields the second implication in
\begin{equation*}
s\not<t\quad\Rightarrow\quad t\leq s\quad\Rightarrow\quad\Gamma(f)(t)\leq\Gamma(f)(s)\quad\Rightarrow\quad\Gamma(f)(s)\not<\Gamma(f)(t).
\end{equation*}
By a straightforward induction over terms, one checks that $\Gamma$ is functorial. A similar induction shows that supports are natural, in the sense that we have
\begin{equation*}
[f]^{<\omega}\circ\supp^\Gamma_X={\supp^\Gamma_Y}\circ\Gamma(f).
\end{equation*}
To conclude that $\Gamma$ is a predilator, it remains to prove
\begin{equation*}
\rng(\Gamma(f))=\{t\in\Gamma(Y)\,|\,\supp^\Gamma_Y(t)\subseteq\rng(f)\}.
\end{equation*}
The inclusion from left to right follows from naturality, as $t=\Gamma(f)(s)$ yields
\begin{equation*}
\supp^\Gamma_Y(t)={\supp^\Gamma_Y}\circ\Gamma(f)(s)=[f]^{<\omega}\circ\supp^\Gamma_X(s)\subseteq\rng(f).
\end{equation*}
In the converse direction, a straightforward induction on the term $t\in\Gamma(Y)$ shows that $\supp^\Gamma_Y(t)\subseteq\rng(f)$ entails $t=\Gamma(f)(s)$ for some~$s\in\Gamma^+(X)$. To get $s\in\Gamma(X)$, we invoke the first equivalence in this proof. If $\Pi^1_1$-comprehension is available, then any subset of~$\mathbb N$ is contained in a countable coded $\omega$-model of arithmetical transfinite recursion, by \cite[Theorems~VII.2.7 and~2.10]{simpson09}. This principle is equivalent to the statement that $\Gamma(X)$ is well founded for any well order~$X$, by~\cite[Theorem~1.4]{rathjen-atr}.
\end{proof}

From~\cite[Theorem~3.5]{freund-predicative-collapsing} we know that $\Gamma(X)$ is a minimal Bachmann-Howard fixed point of a dilator $D$ with $D(Y)=1+2\times Y^2+X$. By the first part of the present section, it should not be hard to characterize $\Gamma(X)$ as a $1$-fixed point. Together with Theorem~\ref{thm:Pi11-rec-to-wf}, this would yield another proof that $\Gamma$ is a dilator.

Recall that a function~$f$ from ordinals to ordinals is normal if it is strictly increasing and continuous, where the latter means that $f(\lambda)=\sup\{f(\alpha)\,|\,\alpha<\lambda\}$ holds when $\lambda$ is a limit. In~\cite{freund-rathjen_derivatives} we have combined previous work of P.~Aczel~\cite{aczel-normal-functors} and J.-Y.~Girard~\cite{girard-pi2}, to define a class of `normal dilators' that induce normal functions on the ordinals. Informally, normal dilators admit internal versions of themselves:

\begin{definition}\label{def:Gamma-normal}
For each linear order~$X$, define $\gamma_X:X\to\Gamma(X)$ by~$\gamma_X(x):=\Gamma_x$.
\end{definition}

The following means that $\Gamma$ is normal in the sense of~\cite{freund-rathjen_derivatives}.

\begin{lemma}\label{lem:Gamma-normal}
For all $s\in\Gamma(X)$ and $x\in X$ we have
\begin{equation*}
s<_{\Gamma(X)}\gamma_X(x)\quad\Leftrightarrow\quad\supp^\Gamma_X(s)\subseteq_X x.
\end{equation*}
Each function $\gamma_X$ is an embedding, we have $\supp^\Gamma_X(\gamma_X(x))=\{x\}$ for all~$x\in X$, and the naturality property $\Gamma(f)\circ\gamma_X=\gamma_Y\circ f$ holds for any embedding $f:X\to Y$.
\end{lemma}
\begin{proof}
The equivalence is readily established by induction on the term~$s$, while naturality holds by a straightforward computation.
\end{proof}

An initial segment of an order~$Y$ is a suborder~$Y_0\subseteq Y$ such that $y<_Y y'\in Y_0$ entails $y\in Y_0$. Let us record an important consequence of normality.

\begin{corollary}\label{cor:normal-initial}
If the range of $f:X\to Y$ is an initial segment of~$Y$, then the range of $\Gamma(f):\Gamma(X)\to\Gamma(Y)$ is an initial segment of~$\Gamma(Y)$.
\end{corollary}
\begin{proof}
Consider an inequality $s<t\in\rng(\Gamma(f))$. To get $s\in\rng(\Gamma(f))$ we need only show $\supp^\Gamma_Y(s)\subseteq\rng(f)$, due to the support condition from Definition~\ref{def:dilator}. Aiming at a contradiction, assume that we have an element $y\in\supp^\Gamma_Y(s)$ with $y\notin\rng(f)$. Given that $\rng(f)$ is an initial segment, we obtain $y'<y$ for all $y'\in\rng(f)$. In view of $t\in\rng(\Gamma(f))$ we can write $t=\Gamma(f)(t_0)$. The naturality of supports yields
\begin{equation*}
\supp^\Gamma_Y(t)={\supp^\Gamma_Y}\circ\Gamma(f)(t_0)=[f]^{<\omega}\circ\supp^\Gamma_X(t_0)\subseteq_Y y.
\end{equation*}
Also note that $y\not<y$ entails $\supp^\Gamma_Y(s)\not\subseteq_Y y$. Now the previous lemma allows us to infer $t<\gamma_Y(y)\leq s$, which contradicts the assumption $s<t$.
\end{proof}

We now represent the total Veblen function. In the following, the first two cases do not clash as $h(0)=0$, and the third case applies precisely when $\overline\varphi st$ is defined.

\begin{definition}\label{def:total-Veblen}
Let $\varphi:\Gamma(X)^2\to\Gamma(X)$ be given by
\begin{equation*}
\varphi_st:=\varphi st:=\varphi(s,t):=\begin{cases}
t & \text{if $s<_{\Gamma(X)}h(t)$},\\
s & \text{if $s\in\scr$ and $t=0$},\\
\overline\varphi st & \text{otherwise}.
\end{cases}
\end{equation*}
\end{definition}

Let us determine the range and fixed points of the Veblen function, its monotonicity properties, and comparisons with terms of the various forms.

\begin{proposition}\label{prop:total-Veblen}
We have $\hau=\{\varphi st\,|\,s,t\in\Gamma(X)\}$ and
\begin{equation*}
\scr=\{\Gamma_x\,|\,x\in X\}=\{s\in\Gamma(X)\,|\,\varphi s0=s\}.
\end{equation*}
Fixed points in the second argument are characterized by
\begin{equation*}
\varphi st=t\quad\Leftrightarrow\quad s<_{\Gamma(X)}t\in\scr\text{ or }t=\varphi t_0t_1\text{ for some }t_i\in\Gamma(X)\text{ with }s<_{\Gamma(X)}t_0.
\end{equation*}
For all $s,s',t,t'\in\Gamma(X)$ we have $s,t\leq_{\Gamma(X)}\varphi st$ and
\begin{equation*}
t'<_{\Gamma(X)}t\,\Rightarrow\,\varphi st'<_{\Gamma(X)}\varphi st\qquad\text{and}\qquad s'<_{\Gamma(X)}s\,\Rightarrow\,\varphi s't\leq_{\Gamma(X)}\varphi st.
\end{equation*}
Finally, we always have
\begin{align*}
\varphi st<_{\Gamma(X)}\Gamma_x\quad&\Leftrightarrow\quad s<_{\Gamma(X)}\Gamma_x\text{ and }t<_{\Gamma(X)}\Gamma_x,\\
\varphi s't'<_{\Gamma(X)}\varphi st\quad&\Leftrightarrow\quad\begin{cases}
s'<_{\Gamma(X)}s\text{ and }t'<_{\Gamma(X)}\varphi st,\\
\text{or }s'=s\text{ and }t'<_{\Gamma(X)}t,\\
\text{or }s<_{\Gamma(X)}s'\text{ and }\varphi s't'<_{\Gamma(X)}t.
\end{cases}
\end{align*}
\end{proposition}
\begin{proof}
To obtain the characterization of~$\hau$, it suffices to note that the first case in Definition~\ref{def:total-Veblen} can only apply when we have $h(t)\neq 0$ and hence $t\in\hau$. The characterization of~$\scr$ is immediate. In the first equivalence, the left side amounts to $s<h(t)$, from which the right side is readily inferred. For the other direction, we need only observe that we always have $t_0\leq h(\varphi t_0t_1)$. In view of Lemma~\ref{lem:Gamma-linear}(a), the claim that we have $s,t\leq\varphi st$ reduces to the following observation: Due to the same lemma, we always have $h(t)\leq t$, so that $s<h(t)$ entails $s<t=\varphi st$. Monotonicity in the second argument is established by a case distinction. In the most interesting case, we have $s\in\scr$ and $t'=0$, so that we get $\varphi st'=s\leq\varphi st$. Aiming at a contradiction, we assume $\varphi st=s$. This value cannot arise by the second or third case from Definition~\ref{def:total-Veblen}, as $t'<t$ entails $t\neq 0$ and since $s$ and $\overline\varphi st$ are different terms. In the remaining case, we would have $s<h(t)$ and  $\varphi st=t$. But this would yield $s=t$ and hence $s<h(s)$, against an observation above. A~similar case distinction yields weak monotonicity in the first argument (note that $s'<s$ and $\varphi st=t$ lead to $\varphi s't=\varphi s'(\varphi st)=\varphi st$ by the fixed point property). The equivalence that characterizes~$\varphi st<\Gamma_x$ is immediate except when we have $s<h(t)$. In this case, we observe that the left side of the equivalence entails
\begin{equation*}
s<h(t)\leq t=\varphi st<\Gamma_x.
\end{equation*}
In the final equivalence of the proposition, the implication from right to left follows from the fixed point and monotonicity properties, e.\,g., because we have
\begin{equation*}
s'<s\text{ and }t'<\varphi st\quad\Rightarrow\quad \varphi s't'<\varphi s'(\varphi st)=\varphi st.
\end{equation*}
Conversely, assume that the right side of the last equivalence in the lemma is false. If we have $s'<s$ and $\varphi st=t'$, then we get $\varphi s't'=\varphi s'(\varphi st)=\varphi st$, so that the left side is false as well. A similar argument applies when we have $s<s'$ and $t=\varphi s't'$.  If we have $s'=s$ and $t'=t$, then the claim is immediate. In all remaining cases, the right side will hold after we interchange $s$ with $s'$ as well as $t$ with $t'$. By the direction from right to left, we get $\varphi st<\varphi s't'$, so that $\varphi s't'<\varphi st$ fails again.
\end{proof}
To conclude this section, we discuss some ordinal arithmetic that will be used later. It may help to recall that $\Gamma(X)$ is isomorphic to the ordered set~$\omega(\hau)$ of finite non\-increasing sequences in~$\hau$, as observed in the paragraph after Definition~\ref{def:Gamma(X)}. Indeed, the following corresponds to the usual operation from ordinal arithmetic, if one thinks of $\langle t_0,\ldots,t_{n-1}\rangle$ as the Cantor normal form~$\omega^{t_0}+\ldots+\omega^{t_{n-1}}$.

\begin{definition}\label{def:addition}
Let $+:\Gamma(X)^2\to\Gamma(X)$ be given by
\begin{multline*}
\langle s_0,\ldots,s_{m-1}\rangle+\langle t_0,\ldots,t_{n-1}\rangle:=\langle s_0,\ldots,s_{i-1},t_0,\ldots,t_{n-1}\rangle\\
\text{with}\quad i:=\begin{cases}
m & \text{if $m=0$ or $n=0$ or $t_0\leq_{\Gamma(X)}s_{m-1}$},\\
\min\{i<m\,|\,s_i<_{\Gamma(X)}t_0\} & \text{otherwise}.
\end{cases}
\end{multline*}
\end{definition}

The following is readily verified and standard (see \cite[Chapter~V.14.3]{schuette77}).

\begin{lemma}\label{lem:addition}
For all $r,r',s,t\in\Gamma(X)$ the following holds:
\begin{enumerate}[label=(\alph*)]
\item We have $(r+s)+t=r+(s+t)$ and $t+0=t=0+t$.
\item Given $s<_{\Gamma(X)}t$, we get $r+s<_{\Gamma(X)}r+t$ and $s+r\leq_{\Gamma(X)}t+r$.
\item If we have $t\in\hau$, then $r<_{\Gamma(X)}r'+t$ and $s<_{\Gamma(X)}t$ entail $r+s<_{\Gamma(X)}r'+t$.
\item We have $r\leq_{\Gamma(X)}t$ if, and only if, there is an $s\in\Gamma(X)$ with $r+s=t$.
\end{enumerate}
\end{lemma}

As $\overline\varphi_00$ is the smallest element of~$\mathsf H\subseteq\Gamma(X)$, the map
\begin{equation*}
\mathbb N\ni n\mapsto n:=\overline n:=\underbrace{\langle\overline\varphi_00,\ldots,\overline\varphi_00\rangle}_{\text{$n$ entries}}\in\Gamma(X)
\end{equation*}
embeds $\mathbb N$ as an initial segment of~$\Gamma(X)$. Addition on $\mathbb N$ and $\Gamma(X)$ are related by
\begin{equation*}
\overline m+t=\begin{cases}
\overline{m+n} & \text{if $t=\overline n$},\\
t & \text{if $t\neq\overline n$ for all~$n\in\mathbb N$}.
\end{cases}
\end{equation*}
In particular, this makes it harmless to write $n$ at the place of~$\overline n$. Instead of a binary multiplication, we use $t\mapsto 1+t$ to define a unary operation $t\mapsto\omega\cdot t$ with
\begin{gather*}
\omega\cdot 0:=0,\quad \omega\cdot\Gamma_z:=\Gamma_z,\quad\omega\cdot\overline\varphi_0t:=\overline\varphi_0(1+t),\quad\omega\cdot\overline\varphi_st:=\overline\varphi_st\text{ for }s\neq 0,\\
\omega\cdot\langle t_0,\ldots,t_{n-1}\rangle:=\langle\omega\cdot t_0,\ldots,\omega\cdot t_{n-1}\rangle.
\end{gather*}
It is not hard to check that $t\mapsto\omega\cdot t$ is strictly increasing, that we have $t\leq\omega\cdot t$, and that $s<\omega\cdot t$ entails $s+n<\omega\cdot t$ for all~$n\in\mathbb N$. Finally, we record how the ordinal arithmetic interacts with supports. The following is immediate in view of Definitions~\ref{def:Gamma-dilator} and~\ref{def:total-Veblen}.

\begin{lemma}\label{lem:Veblen-addition-support}
For any $s,t\in\Gamma(X)$ we have
\begin{equation*}
\supp^\Gamma_X(\varphi st)\cup\supp^\Gamma_X(s+t)\cup\supp^\Gamma_X(\omega\cdot t)\subseteq\supp^\Gamma_X(s)\cup\supp^\Gamma_X(t).
\end{equation*}
\end{lemma}

\section{Hierarchies of admissible sets via search trees}\label{sect:search-trees}

Kurt Sch\"utte's method of search trees (also known as deduction chains) can be used to prove completeness and to construct models in various settings, including predicate and $\omega$-logic~\cite{schuette56,schuette77}, second order arithmetic~\cite{rathjen-afshari,jaeger-strahm-bi-reflection} and set theory~\cite{freund-reflection-induction}. In the present section, we use search trees to construct hierarchies of admissible sets. This extends the construction of a single admissible set in~\cite[Section~4]{freund-equivalence}.

We will search for admissible sets within the constructible hierarchy. Given a transitive set~$u$, set $\mathbb L^u_0:=u$, let $\mathbb L^u_{\alpha+1}$ consist of the $\Delta_0$-definable subsets of~$\mathbb L^u_\alpha$, and put $\mathbb L^u_\lambda:=\bigcup_{\alpha<\lambda}\mathbb L^u_\alpha$ when $\lambda$ is a limit. The restriction to $\Delta_0$-formulas is not essential but will have technical advantages. 

In many of our arguments, the actual hierarchy $\mathbb L^u$ will be represented by a functorial variant~$\mathbf L^u$. This ensures that we get a dilator, to which the well ordering principle from Definition~\ref{def:nu-collapse} can be applied. The functor $\mathbf L^u$ has been introduced in~\cite[Section~3]{freund-equivalence}, based on the first author's PhD thesis~\cite{freund-thesis}. Central facts are recalled in the following, but we refer to~\cite{freund-equivalence} for full details.

First, each linear order~$Y$ gives rise to a set $\mathbf L^u_Y$, which consists of `constant symbols' from~$u$ and terms of the form $L^u_s$ or $\{x\in L^u_s\,|\,\varphi(x,a_0,\ldots,a_{n-1})\}$, for an element~$s\in Y$, a $\Delta_0$-formula~$\varphi$ in the language of set theory, and previously constructed terms $a_i\in\mathbf L^u_Y$ that may only involve elements $r\in Y$ with $r<_Y s$ (so that we have $\supp^{\mathbf L}_Y(a_i)\subseteq_Ys$ in the notation below). To be more precise about the notion of formula, we declare that the signature is $\{\in,=\}$, that there are separate symbols for bounded quantifiers (which are thus distinguished from bounded occurrences of the usual quantifiers), and that formulas are in negation normal form. In view of the latter, negation and implication are defined operations that rely on de Morgan's rules and delete double negations. As usual, a formula is $\Delta_0$ or bounded if it only contains bounded quantifiers.

Prior to any functorial considerations, let us point out that we get an interpretation function $\llbracket\cdot\rrbracket:\mathbf L^u_\alpha\to\mathbb L^u_\alpha$ when $Y=\alpha$ is an ordinal. Here $\mathbf L^u_\alpha$ is the term system from above, while $\mathbb L^u_\alpha$ refers to the actual constructible hierarchy. On the functorial side, each order embedding~$f:Y\to Z$ induces a function $\mathbf L^u_f:\mathbf L^u_Y\to\mathbf L^u_Z$, which is defined by a straightforward recursion over terms. Another recursion yields support functions $\supp^{\mathbf L}_Y:\mathbf L^u_Y\to[Y]^{<\omega}$ with
\begin{gather*}
\supp^{\mathbf L}_Y(w)=\emptyset\quad\text{for each constant symbol }w\in u,\qquad\supp^{\mathbf L}_Y(L^u_s)=\{s\},\\\supp^{\mathbf L}_Y(\{x\in L^u_s\,|\,\varphi(x,a_0,\ldots,a_{n-1})\})=\{s\}\cup\textstyle\bigcup_{i<n}\supp^{\mathbf L}_Y(a_i).
\end{gather*}
In the last case, $s$ is the biggest element of the support, due to the aforementioned condition $\supp^{\mathbf L}_Y(a_i)\subseteq_Y s$. Assuming that $u=\{u_i\,|\,i\in\omega\}$ is countable with fixed enumeration, one can define coding and decoding maps
\begin{equation*}
\en^{\mathbf L}_Y:[Y]^{<\omega}\times\omega\to\mathbf L^u_Y\quad\text{and}\quad\operatorname{code}^{\mathbf L}_Y:[Y]^{<\omega}\times\mathbf L^u_Y\to\omega
\end{equation*}
that are natural in~$Y$ and satisfy $\en^{\mathbf L}_Y(y,\operatorname{code}^{\mathbf L}_Y(y,a))=a$ when $\supp^{\mathbf L}_Y(a)\subseteq y$ (for details see~\cite[Theorem~3.7]{freund-equivalence}). Using these codes, one can define orders $<^{\mathbf L}_Y$ on the sets~$\mathbf L^u_Y$, which are compatible with the functions~$\mathbf L^u_f$. This turns $\mathbf L^u$ into a dilator.

The previous constructions may not be too surprising, because there is little interaction between syntax and semantics. However, semantic aspects of the constructible hierarchy can also be recovered on the syntactic level, as we know from proof theoretic work of J\"ager~\cite{jaeger-KPN,jaeger-kripke-platek} (cf.~Sch\"utte's~\cite{schuette64} work on ramified analysis). The relevant considerations are also functorial, as shown in~\cite[Section~3]{freund-equivalence}: Consider the language that extends $\{\in,=\}$ by a constant symbol for each element of~$\mathbf L^u_Y$. By an $\mathbf L^u_Y$-formula we shall mean a formula in this language. The constant symbols that occur in an $\mathbf L^u_Y$-formula will also be called its parameters. Unless noted other\-wise, we assume that $\mathbf L^u_Y$-formulas are closed. Let us assume $\{0,1\}\subseteq u\subseteq\mathbf L^u_Y$, in order to have indices for binary connectives. Then~\cite[Definition~3.12]{freund-equivalence} associates each $\mathbf L^u_Y$-formula~$\varphi$ with a disjunction or conjunction
\begin{equation*}
\varphi\simeq\textstyle\bigvee_{a\in\iota(\varphi)}\varphi_a\quad\text{or}\quad\varphi\simeq\textstyle\bigwedge_{a\in\iota(\varphi)}\varphi_a.
\end{equation*}
Here $\iota(\varphi)=\iota_Y(\varphi)$ is a subset of $\mathbf L^u_Y$ (which may be empty or infinite) and $\varphi_a$ is an $\mathbf L^u_Y$-formula for each~$a\in\iota(\varphi)$. For full details we refer to the cited definition. As an example, we recall that $\varphi=(b\in\{x\in L^u_s\,|\,\theta(x,c)\})$ yields
\begin{equation*}
\varphi\simeq\textstyle\bigvee_{a\in\iota(\varphi)}\theta(a,c)\land a=b\quad\text{with}\quad\iota(\varphi)=\{a\in\mathbf L^u_Y\,|\,\supp^{\mathbf L}_Y(a)\subseteq_Y s\}.
\end{equation*}
If $Y=\alpha$ is an ordinal, then we get a well founded relation by declaring that $\varphi_a$ precedes~$\varphi$ for each $a\in\iota(\varphi)$. In this case, our disjunctions and conjunctions yield an inductive definition of truth for $\mathbf L^u_\alpha$-formulas. The latter coincides with satisfaction in the actual set $\mathbb L^u_\alpha$, under the aforementioned interpretation $\llbracket\cdot\rrbracket:\mathbf L^u_\alpha\to\mathbb L^u_\alpha$. Let us now state the crucial functorial property: For an embedding $f:Y\to Z$, let $\varphi[f]$ be the $\mathbf L^u_Z$-formula that results from a given $\mathbf L^u_Y$-formula~$\varphi$ when each parameter $a$ is replaced by $\mathbf L^u_f(a)$. Then $\varphi$ and $\varphi[f]$ are both disjunctive or both conjunctive, and \cite[Theorem~3.15]{freund-equivalence} yields
\begin{equation*}
\varphi_a[f]=\varphi[f]_{\mathbf L^u_f(a)}\quad\text{when }a\in\iota_Y(\varphi)\text{ or equivalently }\mathbf L^u_f(a)\in\iota_Z(\varphi[f]).
\end{equation*}
Using the constructions that we have just recalled, we will aim to build a hierarchy of $\nu$ admissible sets above a transitive~$u$. The following assumptions will be discharged in the proof of our main theorem. We write $\ord$ for the class of ordinals.

\begin{assumption}\label{ass:u}
Until the end of Section~\ref{sect:ordinal-analysis}, we fix a transitive set~$u$ and a limit ordinal~$\nu$, both countable with fixed enumerations $u=\{u_i\,|\,i\in\mathbb N\}$ and $\nu=\{\nu_i\,|\,i\in\mathbb N\}$ (no relation with the order). The height $o(u):=u\cap\ord$ is assumed to be a successor ordinal~$o(u)>1$. We also assume that $\Pi^1_1$-comprehension holds.
\end{assumption}

The assumption that $u$ and $\nu$ are countable is essential for our approach. On the other hand, the assumption about the height of $u$ has technical reasons and can later be discharged. It entails $\{0,1\}\subseteq u$, which provides the aforementioned indices for binary connectives. Furthermore, it ensures that $\alpha$ is a limit ordinal whenever the same holds for~$o(\mathbb L^u_\alpha)=o(u)+\alpha$ (otherwise we could have~$\alpha=0$). In this situation, the set $\mathbb L^u_\alpha\ni u$ is admissible if it satisfies the following axioms.

\begin{definition}\label{def:enum-ax}
Let $\langle\operatorname{Ax}_n\,|\,n\geq 1\rangle$ enumerate all instances of $\Delta_0$-collection, i.\,e., all sentences (in the signature $\{\in,=\}$ and without parameters) that have the form
\begin{equation*}
\forall z_1,\ldots, z_k\forall v(\forall x\in v\exists y\,\theta(x,y,z_1,\ldots,z_k)\to\exists w\forall x\in v\exists y\in w\,\theta(x,y,z_1,\ldots,z_k))
\end{equation*}
for a $\Delta_0$-formula~$\theta$. Furthermore, let $\operatorname{Ax}_0$ be the sentence $\forall x\exists y.\,x\in y$.
\end{definition}

Let us write $Z^{<\omega}$ for the tree of finite sequences with entries in~$Z$. In~\cite{freund-equivalence} we have built labelled trees $S_Y\subseteq(\mathbf L^u_Y)^{<\omega}$ for all linear orders~$Y$, which represent attempted proofs of contradiction from the axioms~$\ax_n$ and the rules associated with the infinite disjunctions $\varphi\simeq\bigvee_{a\in\iota_Y(\varphi)}\varphi_a$ and conjunctions $\varphi\simeq\bigwedge_{a\in\iota_Y(\varphi)}\varphi_a$ that were mentioned above. By a relativized ordinal analysis, we showed that $S_Y$ cannot be well founded for all well orders~$Y$, assuming a suitable well ordering principle. This allowed us to conclude that $S_Y$ has an infinite branch for some well order~$Y$. Analogous to other proofs of completeness, such a branch determined a model of the axioms~$\ax_n$, i.\,e., a single admissible set. The following construction of $\nu$ admissible sets is similar overall but different in one respect: we will obtain search trees $S_Y^R$ that depend not only on an order~$Y$ but also on an embedding~$R:\nu\to Y$. The latter determines the heights of the admissible sets in our hierarchy. On an intuitive level, one may think of~$R$ as enumerating regular cardinals (cf.~\cite[Definition~4.1]{buchholz-local-predicativity}).

To describe our search trees in detail, we fix some notation and terminology. Given a sequence $\sigma=\langle\sigma_0,\ldots,\sigma_{n-1}\rangle\in Z^{<\omega}$, write $\len(\sigma):=n$ for its length and put $\sigma\!\restriction\!k:=\langle\sigma_0,\ldots,\sigma_{k-1}\rangle$ for any $k\leq\len(\sigma)$. For $z\in Z$ and $\sigma\in Z^{<\omega}$ as before, set $\sigma^\frown z:=\langle\sigma_0,\ldots,\sigma_{n-1},z\rangle$. The support functions of $\mathbf L^u$ induce functions
\begin{gather*}
\supp^S_Y:(\mathbf L^u_Y)^{<\omega}\to[Y]^{<\omega},\\
\supp^S_Y(\langle\sigma_0,\ldots,\sigma_{n-1}\rangle):=\textstyle\bigcup_{i<n}\supp^{\mathbf L}_Y(\sigma_i).
\end{gather*}
Our search trees will be labelled by $\mathbf L^u_Y$-sequents, which are defined as finite sequences of $\mathbf L^u_Y$-formulas. Semantically, one should think of a sequent as the disjunction of its entries. As usual, we use the letters $\Gamma$ and $\Delta$ to denote sequents, and we write $\varphi_0,\ldots,\varphi_{n-1}$ and $\Gamma,\varphi$ at the place of $\langle\varphi_0,\ldots,\varphi_{n-1}\rangle$ and $\Gamma^\frown\varphi$. When the order and multiplicity of formulas do not matter, we treat sequents like finite sets and write, for example, $\varphi\in\Gamma$ to express that $\varphi$ is an entry of~$\Gamma$. The relativization of an $\mathbf L^u_Y$-formula~$\varphi$ to an element~$a\in\mathbf L^u_Y$ is the $\mathbf L^u_Y$-formula $\varphi^a$ that results from~$\varphi$ when we replace all occurrences $\forall x.\,\psi$ and $\exists x.\,\psi$ of unbounded quantifiers by bounded quantifiers $\forall x\in a.\,\psi$ and $\exists x\in a.\,\psi$, respectively. We do not relativize quantifiers that are already bounded, as this is superfluous when $a$ is transitive and contains the original bounds. Finally, we can describe our search trees in detail:

\begin{definition}\label{def:search}
Consider a linear order~$Y$ and a strictly increasing map~$R:\nu\to Y$. Based on the enumeration $\nu=\{\nu_i\,|\,i\in\mathbb N\}$ from Assumption~\ref{ass:u}, we put
\begin{equation*}
L(i):=L^u_{R(\nu_i)}\in\mathbf L^u_Y.
\end{equation*}
We define a tree $S^R_Y\subseteq(\mathbf L^u_Y)^{<\omega}$ and a labelling function $l_Y:S^R_Y\to\text{``$\mathbf L^u_Y$-sequents''}$ by recursion over sequences in $(\mathbf L^u_Y)^{<\omega}$. Concerning the base case, we declare that we have $\langle\rangle\in S^R_Y$ and $l_Y(\langle\rangle)=\langle\rangle$. In the recursion step, it suffices to consider the children of a previously constructed element $\sigma\in S^R_Y$, as we aim to build a tree. First assume $\len(\sigma)=2k$ is even. Assuming that~$k$ codes the pair~$\langle n,i\rangle$, we declare
\begin{equation*}
\sigma^\frown a\in S^R_Y\,:\Leftrightarrow\,a=L(i)\quad\text{and}\quad l_Y(\sigma^\frown L(i)):=l_Y(\sigma),\neg\ax_n^{L(i)}.
\end{equation*}
Here $a=L(i)$ asserts equality as terms, and the superscript refers to relativization. Now assume that $\len(\sigma)=2k+1$ is odd and that $k$ codes the triple $\langle l,m,n\rangle$. We assume that our coding ensures $l,m,n\leq k$. This entails $l<\len(l_Y(\sigma))$, as we append a formula at each even stage and do no delete any formulas in the following. Let $\varphi$ be the \mbox{$l$-th} formula in~$l_Y(\sigma)$. If $\varphi\simeq\bigwedge_{a\in\iota_Y(\varphi)}\varphi_a$ is conjunctive, we define
\begin{equation*}
\sigma^\frown a\in S^R_Y\,:\Leftrightarrow\,a\in\iota_Y(\varphi)\quad\text{and}\quad l_Y(\sigma^\frown a):=l_Y(\sigma),\varphi_a.
\end{equation*}
If $\varphi\simeq\bigvee_{a\in\iota_Y(\varphi)}\varphi_a$ is disjunctive, we put
\begin{equation*}
b:=\operatorname{en}^{\mathbf L}_Y(\supp^S_Y(\sigma\!\restriction\!m),n)\in\mathbf L^u_Y,
\end{equation*}
using the function $\operatorname{en}_Y^{\mathbf L}:[Y]^{<\omega}\times\omega\to\mathbf L^u_Y$ mentioned above. We then declare
\begin{equation*}
\sigma^\frown a\in S^R_Y\,:\Leftrightarrow\,a=0\quad\text{and}\quad l_Y(\sigma^\frown 0):=\begin{cases}
l_Y(\sigma),\varphi_b & \text{if }b\in\iota_Y(\varphi),\\
l_Y(\sigma) & \text{otherwise},
\end{cases}
\end{equation*}
for which we recall that $0\in u\subseteq\mathbf L^u_Y$ holds by Assumption~\ref{ass:u}.
\end{definition}

For $f:\mathbb N\to\mathbf L^u_Y$ we write $f\!\restriction\!k:=\langle f(0),\ldots,f(k-1)\rangle$ and put
\begin{equation*}
\supp^\infty_Y(f):=\textstyle\bigcup_{k\in\mathbb N}\supp^S_Y(f\!\restriction\!k)=\bigcup_{k\in\mathbb N}\supp^{\mathbf L}_Y(f(k))\subseteq Y.
\end{equation*}
Recall that $f$ is a branch of $S^R_Y$ if $f\!\restriction\!k\in S^R_Y$ holds for all $k\in\mathbb N$. Given $\alpha<\nu$, pick an $i\in\mathbb N$ with $\alpha=\nu_i$, and let $k$ code a pair $\langle n,i\rangle$ for some $n\in\mathbb N$. Assuming that $f$ is a branch, we must have $f(2k)=L^u_{R(\alpha)}$, by construction of the search tree. According to \cite[Definition~3.1]{freund-equivalence} we have $\supp^{\mathbf L}_Y(L^u_{R(\alpha)})=\{R(\alpha)\}$, so that we get
\begin{equation*}
R(\alpha)\in\supp^\infty_Y(f)\quad\text{for all $\alpha<\nu$}.
\end{equation*}
If $Y$ is well founded, then so is its suborder~$\supp^\infty_Y(f)$. In the base theory $\atrs$ from Theorem~\ref{thm:main}, we can use axiom beta to get a transitive collapse, i.\,e., an order preserving map from $\supp^\infty_Y(f)$ onto an ordinal. This yields the desired admissibles:

\begin{theorem}\label{thm:branch-to-model}
Assume that $f$ is a branch in~$S^R_Y$ for a well order~$Y$ and a strictly increasing map $R:\nu\to Y$. Let $c:\supp^\infty_Y(f)\to\ord$ be the transitive collapse. Then $\mathbb L^u_{c(R(\alpha))}\ni u$ is an admissible set for every~$\alpha<\nu$.
\end{theorem}

Before we give a proof, we show that our construction of search trees is functorial. This fact will facilitate the proof of our theorem, but its full significance will only become apparent in the next section.

\begin{definition}\label{def:search-trees-functor}
Consider an embedding $g:Y\to Z$ of linear orders. We define
\begin{gather*}
S_g:(\mathbf L^u_Y)^{<\omega}\to (\mathbf L^u_Z)^{<\omega},\\
S_g(\langle\sigma_0,\ldots,\sigma_{n-1}\rangle):=\langle\mathbf L^u_g(\sigma_0),\ldots,\mathbf L^u_g(\sigma_{n-1})\rangle.
\end{gather*}
Under the assumptions of the following proposition, we also write $S_g:S^P_Y\to S^R_Z$ for the restriction with the indicated (co)domain. Furthermore, let us define~$<^S_Y$ as the Kleene-Brouwer order on $(\mathbf L^u_Y)^{<\omega}$ (also called Lusin-Sierpi\'nski order), which is generated by the clauses $\sigma^\frown a<^S_Y\sigma$ and $\sigma^\frown a<^S_Y\sigma^\frown b$ for $a<^{\mathbf L}_Y b$. We also write $<^S_Y$ for the restriction of this relation to a search tree~$S^P_Y$.
\end{definition}

Due to the corresponding properties of~$\mathbf L^u$, it is immediate that the definition turns $Y\mapsto (\mathbf L^u_Y)^{<\omega}$ into a predilator. In particular, we have the support property
\begin{equation*}
\{S_g(\sigma)\,|\,\sigma\in (\mathbf L^u_Y)^{<\omega}\}=\{\tau\in (\mathbf L^u_Z)^{<\omega}\,|\,\supp^S_Z(\tau)\subseteq\rng(g)\}.
\end{equation*}
Under the assumptions of the following proposition, this equation remains valid when we replace $(\mathbf L^u_Y)^{<\omega}$ and $(\mathbf L^u_Z)^{<\omega}$ by $S^P_Y$ and $S^R_Z$, respectively.

\begin{proposition}\label{prop:search-trees-supports}
Consider linear orders $Y$ and $Z$ with embeddings $P:\nu\to Y$ and $R:\nu\to Z$. If the embedding $g:Y\to Z$ satisfies $g\circ P=R$, then
\begin{equation*}
\sigma\in S^P_Y\quad\Leftrightarrow\quad S_g(\sigma)\in S^R_Z
\end{equation*}
holds for all $\sigma\in(\mathbf L^u_Y)^{<\omega}$.
\end{proposition}
\begin{proof}
Recall that we have a map $\varphi\mapsto\varphi[g]$ from $\mathbf L^u_Y$-formulas to $\mathbf L^u_Z$-formulas. We extend this map to sequents, by setting
\begin{equation*}
\Gamma[g]:=\varphi_0[g],\ldots,\varphi_{n-1}[g]\quad\text{for}\quad\Gamma=\varphi_0,\ldots,\varphi_{n-1}.
\end{equation*}
By induction over the sequence~$\sigma$, we prove the equivalence from the proposition and simultaneously
\begin{equation*}
l_Y(\sigma)[g]=l_Z(S_g(\sigma))\quad\text{when}\quad\sigma\in S^P_Y.
\end{equation*}
The base case with $\sigma=\langle\rangle$ is immediate. In the induction step, we may assume that we have $\sigma\in S^P_Y$ or equivalently $S_g(\sigma)\in S^R_Z$, as we are concerned with trees. First assume that $\len(\sigma)=\len(S_g(\sigma))=2k$ is even, where $k$ codes~$\langle n,i\rangle$. Refining the notation from Definition~\ref{def:search}, we write $L[Y](i):=L^u_{P(\nu_i)}$ and $L[Z](i):=L^u_{R(\nu_i)}$. As \cite[Definition~3.5]{freund-equivalence} yields $\mathbf L^u_g(L^u_a)=L^u_{g(a)}$, we get
\begin{equation*}
\mathbf L^u_g(L[Y](i))=L^u_{g\circ P(\nu_i)}=L^u_{R(\nu_i)}=L[Z](i).
\end{equation*}
Since $\mathbf L^u_g$ is injective on terms (recall that it respects~$<^{\mathbf L}$), we can conclude
\begin{align*}
\sigma^\frown a\in S^P_Y\quad\Leftrightarrow\quad a=L[Y](i)\quad&\Leftrightarrow\quad\mathbf L^u_g(a)=L[Z](i)\\
{}&\Leftrightarrow\quad S_g(\sigma^\frown a)=S_g(\sigma)^\frown\mathbf L^u_g(a)\in S^R_Z.
\end{align*}
In order to see that the desired relation between the sequent labels is preserved, it suffices to observe that we get $\ax_n^{L[Y](i)}[g]=\ax_n^{L[Z](i)}$ from the above (since the operation $\varphi\mapsto\varphi[g]$ replaces any parameter $a$ by $\mathbf L^u_g(a)$). For the case in which the sequences $\sigma$ and $S_g(\sigma)$ have odd length $2k+1$, we refer to the detailed argument in the proof of~\cite[Proposition~4.8]{freund-equivalence} (where the tuple $\langle l,m,n\rangle$ with code~$k$ is written as $\langle\pi_0(n),\pi_1(n),\pi_2(n)\rangle$ with code~$n$).
\end{proof}

Let us now establish the theorem that was stated above.

\begin{proof}[Proof of Theorem~\ref{thm:branch-to-model}]
As preparation, we provide a reduction to the case where the inclusion $\supp^\infty_Y(f)\subseteq Y$ is an equality. Let $g:\kappa\to Y$ be the increasing enumeration of $\supp^\infty_Y(f)$, so that we have $c(g(\gamma))=\gamma$ for $\gamma<\kappa$. Define $P:\nu\to\kappa$ by stipulating $g\circ P=R$, which yields $P(\alpha)=c(R(\alpha))$. For each~$k\in\mathbb N$ we have
\begin{equation*}
\supp^{\mathbf L}_Y(f(k))\subseteq\supp^\infty_Y(f)=\rng(g).
\end{equation*}
By the support property for~$\mathbf L^u$ (see~\cite[Proposition~3.6]{freund-equivalence}), it follows that $f(k)$ lies in the range of~$\mathbf L^u_g:\mathbf L^u_\kappa\to\mathbf L^u_Y$. We thus get an $h:\mathbb N\to\mathbf L^u_\kappa$ with $\mathbf L^u_g\circ h=f$. Since
\begin{equation*}
S_g(h\!\restriction\!k)=\langle\mathbf L^u_g\circ h(0),\ldots,\mathbf L^u_g\circ h(k-1)\rangle=f\!\restriction\!k\in S^R_Y
\end{equation*}
holds for all~$k\in\mathbb N$, we can use Proposition~\ref{prop:search-trees-supports} to conclude that $h$ is a branch of~$S^P_\kappa$. By the naturality of supports for~$\mathbf L^u$ (see again~\cite[Proposition~3.6]{freund-equivalence}), we get
\begin{multline*}
\{g(\gamma)\,|\,\gamma\in\supp^\infty_\kappa(h)\}=\textstyle\bigcup_{k\in\mathbb N}[g]^{<\omega}(\supp^{\mathbf L}_\kappa(h(k))\\
{}=\textstyle\bigcup_{k\in\mathbb N}\supp^{\mathbf L}_Y(\mathbf L^u_g\circ h(k))=\supp^\infty_Y(f).
\end{multline*}
This shows $\supp^\infty_\kappa(h)=\kappa$, which was the purpose of our preparatory construction. To formulate the central claim of this proof, we say that an $\mathbf L^u_\kappa$-formula $\varphi$ occurs on $h$ if we have $\varphi\in l_\kappa(h\!\restriction\!k)$ for some~$k\in\mathbb N$. Let us also recall that we can evaluate $\mathbf L^u_\kappa$-formulas in $\mathbb L^u_\kappa$, via the aforementioned interpretation $\llbracket\cdot\rrbracket:\mathbf L^u_\kappa\to\mathbb L^u_\kappa$. Crucially, we shall show that $\mathbb L^u_\kappa$ satisfies $\neg\varphi$ whenever $\varphi$ occurs on~$h$. According to~\cite[Theorem~3.14]{freund-equivalence}, this reduces to the following claims:
\begin{enumerate}[label=(\roman*)]
\item if $\varphi\simeq\bigwedge_{a\in\iota_\kappa(\varphi)}\varphi_a$ occurs on~$h$, then so does $\varphi_a$ for some $a\in\iota_\kappa(\varphi)$,
\item if $\varphi\simeq\bigvee_{a\in\iota_\kappa(\varphi)}\varphi_a$ occurs on~$h$, then so does $\varphi_a$ for all $a\in\iota_\kappa(\varphi)$.
\end{enumerate}
Indeed, we get a well founded relation on $\mathbf L^u_\kappa$-formulas by declaring that each $\varphi_a$ precedes~$\varphi$, as mentioned above. Given~(i) and~(ii), transfinite induction over this relation shows that each $\varphi$ on~$h$ must fail in~$\mathbb L^u_\kappa$. The proof of~\cite[Theorem~3.14]{freund-equivalence} shows that this inductive argument goes through in our base theory. Before we establish~(i) and~(ii), let us explain how to derive the theorem: Given any $\alpha<\nu$ and $n\in\mathbb N$, let $k$ be the code of a pair~$\langle n,i\rangle$ with~$\nu_i=\alpha$. By construction of our search trees, the formula $\neg\ax_n^{L(i)}$ occurs in $l_\kappa(h\!\restriction\!(2k+1))$ and hence on~$h$. In view of~\cite[Definition~3.2]{freund-equivalence} we have
\begin{equation*}
\llbracket L(i)\rrbracket=\llbracket L^u_{P(\alpha)}\rrbracket=\mathbb L^u_{P(\alpha)}.
\end{equation*}
Hence our central claim entails that $\mathbb L^u_\kappa$ satisfies the relativization of $\ax_n$ to $\mathbb L^u_{P(\alpha)}$. But this simply means that $\mathbb L^u_{P(\alpha)}$ satisfies $\ax_n$. It follows that $\mathbb L^u_{P(\alpha)}=\mathbb L^u_{c(R(\alpha))}$ is admissible (cf.~the paragraph before Definition~\ref{def:enum-ax}), as required by our theorem. Claims~(i) and~(ii) are established as in the proof of~\cite[Theorem~4.6]{freund-equivalence}. However, the fact that we have $\supp^\infty_\kappa(h)=\kappa$ does simplify matters. We provide details for the more difficult claim~(ii): Assume that the disjunctive formula~$\varphi$ occurs on~$h$, say as the \mbox{$l$-th} formula in~$l_\kappa(h\!\restriction\!m_0)$. Given an arbitrary $a\in\iota_\kappa(\varphi)$, we observe
\begin{equation*}
\supp^{\mathbf L}_\kappa(a)\subseteq\kappa=\supp^\infty_\kappa(h)=\textstyle\bigcup_{k\in\mathbb N}\supp^S_\kappa(h\!\restriction\!k).
\end{equation*}
Since the last union is increasing, we may pick a number $m\geq m_0$ such that the finite set $\supp^{\mathbf L}_\kappa(a)$ is contained in $\supp^S_\kappa(h\!\restriction\!m)$. We then have
\begin{equation*}
a=\en^{\mathbf L}_\kappa(\supp^S_\kappa(h\!\restriction\!m),n)\quad\text{for}\quad n:=\operatorname{code}^{\mathbf L}_\kappa(\supp^S_\kappa(h\!\restriction\!m),a),
\end{equation*}
by~\cite[Theorem~3.7]{freund-equivalence} or the discussion above. Let us now define $k$ as the code of the triple~$\langle l,m,n\rangle$. As in Definition~\ref{def:search}, we may assume that our coding of tuples ensures $m\leq k$ and hence $m_0<2k+1$. When we build our search trees, we extend sequents at the end, but we never delete or permute formulas. Thus $\varphi$ is still the \mbox{$l$-th} formula in~$l_\kappa(h\!\restriction\!(2k+1))$. By construction we get
\begin{equation*}
l_\kappa(h\!\restriction\!(2k+2))=l_\kappa(h\!\restriction\!(2k+1)),\varphi_a.
\end{equation*}
Hence $\varphi_a$ occurs on~$h$, as desired.
\end{proof}

Using methods from ordinal analysis, we will show that the well ordering principle from Definition~\ref{def:nu-collapse} entails the following: it cannot be the case that $S^R_Y$ is well founded whenever $Y$ is a well order. Once this is known, Theorem~\ref{thm:branch-to-model} will yield a hierarchy of $\nu$ admissible sets, as needed for the crucial direction of Theorem~\ref{thm:main}. To conclude, we record a fact that will be needed later (cf.~\cite[Corollary~4.10]{freund-equivalence}):

\begin{corollary}\label{cor:search-parameters}
Consider a linear order~$Z$ and an embedding $R:\nu\to Z$. We have
\begin{equation*}
\supp^{\mathbf L}_Z(b)\subseteq\supp^S_Z(\sigma)\cup\{R(\alpha)\,|\,\alpha<\nu\}
\end{equation*}
for any node $\sigma\in S^R_Z$ and any parameter $b$ that occurs in some formula of~$l_Z(\sigma)$.
\end{corollary}
\begin{proof}
Let $Y$ be the set on the right of the desired inclusion, considered as a suborder of~$Z$. Write $\iota:Y\hookrightarrow Z$ for the inclusion, and define $P:\nu\to Y$ by $\iota\circ P=R$. In view of $\supp^S_Z(\sigma)\subseteq\rng(\iota)$ we obtain $\sigma=S_\iota(\rho)$ for some node $\rho\in S^P_Y$, due to Proposition~\ref{prop:search-trees-supports}. By the proof of the latter, we have $l_Y(\rho)[\iota]=l_Z(\sigma)$. We can thus write $b=\mathbf L^u_\iota(a)$ with $a\in\mathbf L^u_Y$, so that
\begin{equation*}
\supp^{\mathbf L}_Z(b)={\supp^{\mathbf L}_Z}\circ\mathbf L^u_\iota(a)=[\iota]^{<\omega}\circ\supp^{\mathbf L}_Y(a)\subseteq\rng(\iota)=Y
\end{equation*}
follows by the naturality of supports.
\end{proof}

\section{From search tree to collapsing functions}\label{sect:search-to-notations}

In this section, we apply the well ordering principle from Definition~\ref{def:nu-collapse} to the search trees $S^R_Y$ that were constructed in Definition~\ref{def:search}. The result is an order~$\ot$, which is quite close to the relativized ordinal notation system in~\cite[Definition~6.4]{rathjen-axiomatic-thinking} (cf.~also \cite{buchholz-new-system} and~\cite[Section~12.2]{thomson-rathjen-Pi-1-1}). We will later use $\ot$ as a basis for the ordinal analysis that proves the implication from~(iii) to~(iv) in Theorem~\ref{thm:main}.

Recall the dilator $\Gamma$ and the functions $\gamma_X:X\to\Gamma(X)$ from Section~\ref{sect:booting-up}. The desired order~$\ot$ will be constructed as part of a system of orders and embeddings, which can be depicted as follows (where a hooked arrow indicates that the range is an initial segment of the codomain, while $\pto$ refers to a partial surjective function):
\begin{equation*}
\begin{tikzcd}
\nu\times\ot\ar[r,twoheadrightarrow,"\psi^{\mathbf X}","p"' {xshift=1.6em}] & \mathbf X \ar[r,hook,"I"]\ar[d,"\gamma_{\mathbf X}"] & \mathbf X+S^{\mathbf R}_{\Gamma(\mathbf X)}=:\mathbf K\ar[r,"\gamma_{\mathbf K}"] & \Gamma(\mathbf K)=:\mathbf O.\\
& \Gamma(\mathbf X)\ar[rru,hook,swap,"\Gamma(I)"] & &
\end{tikzcd}
\end{equation*}
Before we give a formal construction of these objects, let us explain their intuitive meaning. In view of Section~\ref{sect:booting-up}, the order $\mathbf O=\Gamma(\mathbf K)$ is closed under the binary Veblen function and includes the first $\mathbf K$ strongly critical ordinals, which are represented by the elements $\gamma_{\mathbf K}(z)\in\mathbf O$ with~$z\in\mathbf K$ (we choose $\mathbf K$ for `kritisch'). By composing all vertical arrows, we obtain $\nu$-many partial but order preserving `collapsing functions' from~$\mathbf O$ to itself. The values of these functions are represented by the elements of a set~$\mathbf X$. We have a map $I$ that realizes this set as an initial segment of~$\mathbf K$. Since $\Gamma$ is a functor and normal, we also obtain an identification~$\Gamma(I)$ of the set $\Gamma(\mathbf X)$ with an initial segment of~$\mathbf O$ (see Corollary~\ref{cor:normal-initial}). This means, first, that the collapsing values form an initial segment of the strongly critical ordinals. Moreover, it means that the ordinals generated from the collapsing values form an initial segment of the full system~$\mathbf O$. Both properties are typical for ordinal notation systems (see again the examples in~\cite{buchholz-new-system,thomson-rathjen-Pi-1-1}). It is also typical that there are strongly critical ordinals that lie above all collapsing values. In our case, these `large' ordinals correspond to the nodes of a certain search tree~$S^{\mathbf R}_{\Gamma(\mathbf X)}$ (cf.~the elements $\mathfrak E_\sigma$ in~\cite[Definition~5.2]{freund-equivalence}). For our ordinal analysis, it will be crucial that this search tree is built over the lower part~$\Gamma(\mathbf X)$ of the order~$\mathbf O$, with respect to a map $\mathbf R:\nu\to\Gamma(\mathbf X)$ that has a meaningful connection to the collapsing functions. Concerning the latter, we will obtain $\mathbf R(\alpha)=\gamma_{\mathbf X}\circ\psi^{\mathbf X}(\alpha+1,0)$ for $0\in\Gamma(\mathbf K)=\ot$, which evokes $\psi_{\alpha+1}0=\Omega_{\alpha+1}\in R$ from~\cite[Lemma~1.7]{buchholz-new-system} and~\cite[Definition~4.1]{buchholz-local-predicativity}.

We would like to define $\psi^{\mathbf X}:\nu\times\mathbf O\pto\mathbf X$ as the partial inverse of a function~$\pi$ as in Definition~\ref{def:nu-collapse}. Before we can apply the latter, however, we must overcome a significant obstacle. The issue is that Definition~\ref{def:nu-collapse} requires a dilator as input, while the construction of search trees in Definition~\ref{def:search} does not provide one, at least not directly: the tree $S^R_Y$ depends not only on the order~$Y$ but also on a given embedding~$R:\nu\to Y$. This issue will occupy us for most of the present section, and its resolution may at times appear technical. At the same time, we believe that the issue itself is not technical but has real mathematical substance. In particular, it distinguishes the construction of a single admissible set in~\cite{freund-equivalence} -- where no similar issue arose -- from the construction of an infinite hierarchy of admissible sets.

In order to resolve the issue that was mentioned in the previous paragraph, we will precompose the construction of search trees with the order transformation
\begin{equation*}
X\mapsto J(X):=\nu\times\Gamma(X).
\end{equation*}
Recall that products were discussed in the paragraph before Definition~\ref{def:nu-collapse}, which does also explain $J(f):=\nu\times\Gamma(f)$ for an order embedding~$f$. It is straightforward to check that we get a dilator if we provide supports by
\begin{equation*}
\supp^J_X:J(X)\to[X]^{<\omega}\quad\text{with}\quad\supp^J_X(\alpha,\sigma):=\supp^{\Gamma}_X(\sigma).
\end{equation*}
As $\nu$ is a limit by Assumption~\ref{ass:u}, we may consider the embeddings
\begin{equation*}
j[X]:\nu\to J(X)\quad\text{with}\quad j[X](\alpha):=(\alpha+1,0).
\end{equation*}
These are natural in the sense that $J(f)\circ j[X]=j[Y]$ holds for any embedding~$f$, as we have $\Gamma(f)(0)=0$ by Definition~\ref{def:Gamma-dilator}. We can now describe the preprocessed search trees that were mentioned above:

\begin{definition}\label{def:S^E}
Consider the order transformation
\begin{equation*}
X\mapsto\mathbf S_0(X):=S^{j[X]}_{J(X)},
\end{equation*}
where the definiens refers to Definitions~\ref{def:search} and~\ref{def:search-trees-functor}. Invoking the latter in conjunction with Proposition~\ref{prop:search-trees-supports}, we map each embedding $f:X\to Y$ to the embedding
\begin{equation*}
\mathbf S_0(f):=S_{J(f)}:\mathbf S_0(X)\to \mathbf S_0(Y).
\end{equation*}
Note that the cited proposition can be applied because we have $J(f)\circ j[X]=j[Y]$, as seen above. Finally, we define functions $\supp^0_X:\mathbf S_0(X)\to[X]^{<\omega}$ by setting
\begin{equation*}
\supp^0_X(\sigma):=\bigcup\{\supp^{J}_X(\rho)\,|\,\rho\in\supp^S_{J(X)}(\sigma)\}.
\end{equation*}
This relies on the definition of~$\supp^S$ in the paragraph before Definition~\ref{def:search}.
\end{definition}

As we had hoped, our preprocessed search trees form a dilator, at least when statement~(iv) from Theorem~\ref{thm:main} is violated.

\begin{proposition}\label{prop:S^E-dilator}
The constructions from Definition~\ref{def:S^E} yield a predilator~$\mathbf S_0$. The latter is a dilator if there is no sequence of admissible sets $\mathsf{Ad}_\alpha$ with $u\in\mathsf{Ad}_\alpha\in\mathsf{Ad}_\beta$ for $\alpha<\beta<\nu$ (with $u$ and $\nu$ as fixed in Assumption~\ref{ass:u}).
\end{proposition}
\begin{proof}
Let us observe that the first map in
\begin{equation*}
X\mapsto(\mathbf L^u_{J(X)})^{<\omega}\supseteq S^{j[X]}_{J(X)}=\mathbf S_0(X)
\end{equation*}
is the composition of predilators and hence a predilator itself, by the paragraph before Proposition~\ref{prop:search-trees-supports}. Using the latter, we can conclude that $\mathbf S_0$ is also a predilator. To provide details for the crucial step, we show that the support property
\begin{equation*}
\supp^0_Y(\sigma)\subseteq\rng(f)\quad\To\quad\sigma\in\rng(\mathbf S_0(f))
\end{equation*}
holds for any embedding~$f:X\to Y$ and any $\sigma\in\mathbf S_0(Y)$. Given the antecedent of our implication, the definition of~$\supp^0_Y$ and the support property for~$J$ yield
\begin{equation*}
\supp^S_{J(Y)}(\sigma)\subseteq\rng(J(f)).
\end{equation*}
This allows us to write
\begin{equation*}
\sigma=S_{J(f)}(\sigma_0)\quad\text{for some}\quad\sigma_0\in(\mathbf L^u_{J(X)})^{<\omega},
\end{equation*}
by the paragraph before Proposition~\ref{prop:search-trees-supports}. Now the latter ensures that $\sigma\in\mathbf S_0(Y)$ entails $\sigma_0\in\mathbf S_0(X)$ and hence $\sigma=\mathbf S_0(f)(\sigma_0)\in\rng(\mathbf S_0(f))$, as desired. Under the assumption from the proposition, we now show that $\mathbf S_0$ is a dilator. Given a well order~$X$, we must establish that $\mathbf S_0(X)$ is well founded. As $\Pi^1_1$-comprehension is available by Assumption~\ref{ass:u}, we can infer that $\Gamma(X)$ and~$J(X)$ are well orders, by Proposition~\ref{ref:Gamma-dilator} or directly by \cite[Theorem~1.4]{rathjen-atr}. According to \cite[Lemma~3.10]{freund-equivalence}, it follows that $\mathbf L^u_{J(X)}$ is well founded (see the beginning of Section~\ref{sect:search-trees} and compare with the usual constructible hierarchy). Hence $\mathbf S_0(X)$ is well founded (with respect to the Kleene-Brouwer order from Definition~\ref{def:search-trees-functor}) unless it has a branch. In the latter case, Theorem~\ref{thm:branch-to-model} would yield a hierarchy of~$\nu$ admissible sets above~$u$, against the assumption of the present proposition.
\end{proof}

Following the informal explanation at the beginning of this section, we now add space for collapsing values below the elements of our search tree. Sums of linear orders and embeddings are defined as in the proof of Proposition~\ref{prop:omega-1fp}. Recall that elements of $Z_0+Z_1$ are written as $z_0$ and $Z_0+z_1$ with $z_i\in Z_i$.

\begin{definition}\label{def:X+S^E}
For each linear order~$X$ and each embedding~$f:X\to Y$, we put $\mathbf S(X):=X+\mathbf S_0(X)$ and define $\mathbf S(f):\mathbf S(X)\to\mathbf S(Y)$ by $\mathbf S(f):=f+\mathbf S_0(f)$. By
\begin{equation*}
\supp^{\mathbf S}_X(x):=\{x\}\quad\text{and}\quad\supp^{\mathbf S}_X(X+\sigma):=\supp^0_X(\sigma)
\end{equation*}
we define a family of functions $\supp^{\mathbf S}_X:\mathbf S(X)\to[X]^{<\omega}$.
\end{definition}

To prove the crucial direction of Theorem~\ref{thm:main} by contradiction, we will assume that statement~(iv) fails. In view of Proposition~\ref{prop:S^E-dilator}, this will have the effect that $\mathbf S_0$ is a dilator. It is easy to conclude that $\mathbf S$ and $\Gamma\circ\mathbf S$ are dilators as well (recall how composition is defined in the paragraph before Proposition~\ref{thm:1fp-to-Bachmann-Howard}). We bring in statement~(ii) of Theorem~\ref{thm:main} in the form of the following assumption.

\begin{assumption}\label{ass:collapse-search}
Until the end of Section~\ref{sect:ordinal-analysis}, we assume that $\Gamma\circ\mathbf S$ is a dilator. Furthermore, we assume that we have a fixed well order~$\mathbf Y$ and $\nu$-collapse
\begin{equation*}
\pi_{\mathbf Y}:\mathbf Y\to\nu\times(\Gamma\circ\mathbf S)(\mathbf Y)
\end{equation*}
in the sense of Definition~\ref{def:nu-collapse} (with $\nu$ and the suppressed $u$ as in Assumption~\ref{ass:u}).
\end{assumption}

The inverse of $\pi_{\mathbf Y}$ is a partial embedding
\begin{equation*}
\nu\times(\Gamma\circ\mathbf S)(\mathbf Y)=\nu\times\Gamma\left(\mathbf Y+S^{j[\mathbf Y]}_{J(\mathbf Y)}\right)\pto\mathbf Y.
\end{equation*}
This looks a lot like the function
\begin{equation*}
\psi^{\mathbf X}:\nu\times\Gamma\left(\mathbf X+S^{\mathbf R}_{\Gamma(\mathbf X)}\right)\pto\mathbf X
\end{equation*}
that was promised at the beginning of this section. However, one important point remains to be improved: the collapse $\psi^{\mathbf X}$ and the embedding $\mathbf R:\nu\to\Gamma(X)$ were supposed to be connected in a meaningful way, while the function $j[\mathbf Y]:\nu\to J(\mathbf Y)$ and the order~$J(\mathbf Y)$ appear rather ad hoc and unrelated to~$\pi_{\mathbf Y}$. Perhaps surprisingly, we can use $\pi_{\mathbf Y}$ to `infuse meaning' ex post. The following is a preparation.

\begin{lemma}\label{lem:Omega-alpha-defined}
We have $(\alpha,0)\in\rng(\pi_{\mathbf Y})$ for all $\alpha<\nu$.
\end{lemma}
\begin{proof}
By Definition~\ref{def:Gamma-dilator} we have $\supp^{\Gamma}_{\mathbf S(\mathbf Y)}(0)=\emptyset$, which entails
\begin{equation*}
\supp^{\Gamma\circ\mathbf S}_{\mathbf Y}(0)=\bigcup\{\supp^{\mathbf S}_{\mathbf Y}(\rho)\,|\,\rho\in\supp^{\Gamma}_{\mathbf S(\mathbf Y)}(0)\}=\emptyset.
\end{equation*}
In the notation from Definition~\ref{def:nu-collapse}, we get
\begin{equation*}
G_\alpha(0)=\bigcup\{G^{\Gamma\circ\mathbf S}_\alpha(s)\,|\,s\in\supp^{\Gamma\circ\mathbf S}_{\mathbf Y}(0)\}=\emptyset\subseteq_{\Gamma\circ\mathbf S(\mathbf Y)}0.
\end{equation*}
The claim follows by Definition~\ref{def:nu-collapse}.
\end{proof}

Recall that the normal dilator~$\Gamma$ comes with an embedding~$\gamma_{\mathbf Y}:\mathbf Y\to\Gamma(\mathbf Y)$, which is given by Definition~\ref{def:Gamma-normal}.

\begin{definition}
In view of the previous lemma, let the embedding $\mathbf P_0:\nu\to\mathbf Y$ be determined by $\pi_{\mathbf Y}\circ\mathbf P_0(\alpha)=(\alpha+1,0)$. We also put $\mathbf P:=\gamma_{\mathbf Y}\circ\mathbf P_0:\nu\to\Gamma(\mathbf Y)$.
\end{definition}

Given $s\in\Gamma(\mathbf Y)$, let $y\in\mathbf Y$ be the maximal element of $\supp^\Gamma_{\mathbf Y}(s)\cup\{\mathbf P_0(0)\}$, which is finite and non-empty. Write $\pi_{\mathbf Y}(y)=(\alpha,\sigma)$. Since $\pi_{\mathbf Y}$ is an embedding, we get
\begin{equation*}
\supp^\Gamma_{\mathbf Y}(s)\subseteq_{\mathbf Y}\mathbf P_0(\alpha)\quad\text{and thus}\quad s<_{\Gamma(\mathbf Y)}\gamma_{\mathbf Y}\circ\mathbf P_0(\alpha)=\mathbf P(\alpha),
\end{equation*}
using Lemma~\ref{lem:Gamma-normal}. This observation ensures that the following is well defined.

\begin{definition}\label{def:Y-Gamma-J}
We define $Y:\Gamma(\mathbf Y)\to J(\mathbf Y)$ by $Y(\mathbf P(\alpha)):=(\alpha+1,0)$ and
\begin{equation*}
Y(s):=(\alpha,s)\quad\text{with}\quad\alpha=\min\{\gamma<\nu\,|\,s<_{\Gamma(\mathbf Y)}\mathbf P(\gamma)\}
\end{equation*}
for any $s\in\Gamma(\mathbf Y)$ that does not lie in the range of~$\mathbf P$.
\end{definition}

It is not hard to see that $Y$ is an order embedding, and we have $Y\circ \mathbf P=j[\mathbf Y]$ by construction. We can thus invoke Proposition~\ref{prop:search-trees-supports} to obtain embeddings
\begin{alignat*}{3}
S_Y:{}&&S^{\mathbf P}_{\Gamma(\mathbf Y)}&\to S^{j[\mathbf Y]}_{J(\mathbf Y)}=\mathbf S_0(\mathbf Y),\\
\mathbf Y+S_Y:{}&&\mathbf Y+S^{\mathbf P}_{\Gamma(\mathbf Y)}&\to\mathbf Y+\mathbf S_0(\mathbf Y)=\mathbf S(\mathbf Y).
\end{alignat*}
In contrast to $j[\mathbf Y]:\nu\to J(\mathbf Y)$, the map $\mathbf P:\nu\to\Gamma(\mathbf Y)$ has a `natural' codomain and a meaningful connection to~$\pi_{\mathbf Y}$. With respect to the informal discussion at the beginning of this section, it may thus be tempting to define $\mathbf X$ as $\mathbf Y$. The partial function $\psi^{\mathbf X}$ from this discussion should then be inverse to the dashed arrow in
\begin{equation*}
\begin{tikzcd}[column sep=large]
\mathbf Y\arrow[r,dashed,"?"]\arrow[rd,swap,"\pi_{\mathbf Y}"] & \nu\times\Gamma(\mathbf Y+S^{\mathbf P}_{\Gamma(\mathbf Y)})\arrow[d,"\nu\times\Gamma(\mathbf Y+S_Y)"]\\
& \nu\times(\Gamma\circ\mathbf S)(\mathbf Y).
\end{tikzcd}
\end{equation*}
However, it seems that the range of~$\pi_{\mathbf Y}$ need not be contained in the range of the vertical arrow, so that the dashed arrow may not exist. To resolve this issue, we define a suborder that guarantees the desired inclusion in a hereditary way.

\begin{definition}\label{def:XsubY}
Let us write $\tl$ for the well founded relation on $\mathbf Y$ that is provided by Definition~\ref{def:nu-collapse}, which means that we have
\begin{equation*}
x\tl y\quad\Leftrightarrow\quad x\in\supp^{\Gamma\circ\mathbf S}_{\mathbf Y}(s)\text{ for }\pi_{\mathbf Y}(y)=(\alpha,s).
\end{equation*}
By recursion over this relation, we define a suborder $\mathbf X\subseteq\mathbf Y$ with
\begin{equation*}
y\in\mathbf X\quad:\Leftrightarrow\quad\pi_{\mathbf Y}(y)\in\rng(\nu\times\Gamma(\mathbf Y+S_Y))\text{ and }x\in\mathbf X\text{ for all }x\tl y.
\end{equation*}
We will write $\iota:\mathbf X\rightarrow\mathbf Y$ for the inclusion.
\end{definition}

Let us complement Lemma~\ref{lem:Omega-alpha-defined} as follows.

\begin{lemma}\label{lem:Omega-alpha-X}
If $\pi_{\mathbf Y}(y)=(\alpha,0)$ holds for some $\alpha<\nu$, then we have $y\in\mathbf X=\rng(\iota)$.
\end{lemma}
\begin{proof}
It suffices to recall that $0=\Gamma(f)(0)\in\rng(\Gamma(f))$ holds for any embedding~$f$, and that $\supp^{\Gamma\circ\mathbf S}_{\mathbf Y}(0)=\emptyset$ was shown in the proof of Lemma~\ref{lem:Omega-alpha-defined}.
\end{proof}

To define the other objects that were promised at the beginning of this section, we repeat some of the previous constructions, but now with $\mathbf X$ at the place of~$\mathbf Y$.

\begin{definition}\label{def:Beth-O}
Determine $\mathbf R_0:\nu\to\mathbf X$ and $\mathbf R:\nu\to\Gamma(\mathbf X)$ by
\begin{equation*}
\pi_{\mathbf Y}\circ\iota\circ\mathbf R_0(\alpha)=(\alpha+1,0)\quad\text{and}\quad\mathbf R:=\gamma_{\mathbf X}\circ\mathbf R_0.
\end{equation*}
For the order $S^{\mathbf R}_{\Gamma(\mathbf X)}$ given by Definitions~\ref{def:search} and~\ref{def:search-trees-functor}, we now put
\begin{equation*}
\mathbf K:=\mathbf X+S^{\mathbf R}_{\Gamma(\mathbf X)}\quad\text{and}\quad\mathbf O:=\Gamma(\mathbf K).
\end{equation*}
\end{definition}

Note that we have $\iota\circ\mathbf R_0=\mathbf P_0$, as $\pi_{\mathbf Y}$ is order preserving and hence injective. From Lemma~\ref{lem:Gamma-normal} we know that $\gamma$ is natural with respect to~$\iota:\mathbf X\rightarrow\mathbf Y$. We get
\begin{equation*}
\Gamma(\iota)\circ\mathbf R=\Gamma(\iota)\circ\gamma_{\mathbf X}\circ\mathbf R_0=\gamma_{\mathbf Y}\circ\iota\circ\mathbf R_0=\gamma_{\mathbf Y}\circ\mathbf P_0=\mathbf P.
\end{equation*}
Thus Proposition~\ref{prop:search-trees-supports} yields an embedding $S_{\Gamma(\iota)}:S^{\mathbf R}_{\Gamma(\mathbf X)}\to S^{\mathbf P}_{\Gamma(\mathbf Y)}$. By composing with another map from above, we obtain embeddings
\begin{align*}
(\mathbf Y+S_Y)\circ(\iota+S_{\Gamma(\iota)})=\iota+S_{Y\circ\Gamma(\iota)}&:\mathbf K\to\mathbf S(\mathbf Y),\\
\Gamma(\iota+S_{Y\circ\Gamma(\iota)})&:\mathbf O=\Gamma(\mathbf K)\to\Gamma\circ\mathbf S(\mathbf Y).
\end{align*}
In particular, we can conclude that $\mathbf O$ is a well order, as $\Gamma\circ\mathbf S(\mathbf Y)$ is well founded by Assumption~\ref{ass:collapse-search}. The following resolves an issue that was mentioned above. It may help to read the lemma in conjunction with the definition that follows it.

\begin{lemma}\label{lem:def-pi_X}
The range of~$\pi_{\mathbf Y}\circ\iota$ is contained in the range of $\nu\times\Gamma(\iota+S_{Y\circ\Gamma(\iota)})$.
\end{lemma}
\begin{proof}
The crucial step is to show that any $s\in\mathbf Y+S^{\mathbf P}_{\Gamma(\mathbf Y)}$ validates
\begin{equation*}
s\in\rng(\iota+S_{\Gamma(\iota)})\quad\Leftrightarrow\quad{\supp^{\mathbf S}_{\mathbf Y}}\circ(\mathbf Y+S_Y)(s)\subseteq\mathbf X=\rng(\iota).
\end{equation*}
Even though we will not use this fact, we note that the equivalence means that
\begin{equation*}
\begin{tikzcd}[column sep=large]
\mathbf X+S^{\mathbf R}_{\Gamma(\mathbf X)}\ar[r,"\mathbf X+S_X"]\ar[d,swap,"\iota+S_{\Gamma(\iota)}"] & \mathbf S(\mathbf X)\ar[d,"\mathbf S(\iota)"]\\
\mathbf Y+S^{\mathbf P}_{\Gamma(\mathbf Y)}\ar[r,"\mathbf Y+S_Y"] & \mathbf S(\mathbf Y)
\end{tikzcd}
\end{equation*}
is a pullback, where $X:\Gamma(\mathbf X)\to J(\mathbf X)$ is constructed analogous to Definition~\ref{def:Y-Gamma-J}. For $s=y\in\mathbf Y\subseteq\mathbf Y+S^{\mathbf P}_{\Gamma(\mathbf Y)}$ we can invoke Definition~\ref{def:X+S^E} to get
\begin{equation*}
{\supp^{\mathbf S}_{\mathbf Y}}\circ(\mathbf Y+S_Y)(s)=\supp^{\mathbf S}_{\mathbf Y}(y)=\{y\}.
\end{equation*}
So both sides of our equivalence amount to $y\in\rng(\iota)$. For $s=\mathbf Y+\sigma$ we have
\begin{equation*}
s\in\rng(\iota+S_{\Gamma(\iota)})\quad\Leftrightarrow\quad\sigma\in\rng(S_{\Gamma(\iota)})\quad\Leftrightarrow\quad \supp^S_{\Gamma(\mathbf Y)}(\sigma)\subseteq\rng(\Gamma(\iota)),
\end{equation*}
where the second equivalence holds by Proposition~\ref{prop:search-trees-supports} and the paragraph before it. On the other hand, Definitions~\ref{def:S^E} and~\ref{def:X+S^E} yield
\begin{multline*}
{\supp^{\mathbf S}_{\mathbf Y}}\circ(\mathbf Y+S_Y)(s)=\supp^{\mathbf S}_{\mathbf Y}(\mathbf Y+S_Y(\sigma))=\supp^0_{\mathbf Y}(S_Y(\sigma))\\
\quad{}=\bigcup\{\supp^J_{\mathbf Y}(\rho)\,|\,\rho\in\supp^S_{J(\mathbf Y)}(S_Y(\sigma))\}=\bigcup\{\supp^J_{\mathbf Y}(Y(\tau))\,|\,\tau\in\supp^S_{\Gamma(\mathbf Y)}(\sigma)\}.
\end{multline*}
Here the last equality relies on the fact that $\supp^S$ is a natural transformation. By the previous lines of equivalences and equations, the desired equivalence reduces to
\begin{equation*}
\tau\in\rng(\Gamma(\iota))\quad\Leftrightarrow\quad\supp^J_{\mathbf Y}(Y(\tau))\subseteq\mathbf X=\rng(\iota).
\end{equation*}
Considering the definition of~$Y$, we distinguish two cases: For $\tau=\mathbf P(\alpha)$, the paragraph after Definition~\ref{def:Beth-O} yields $\tau=\Gamma(\iota)\circ\mathbf R(\alpha)\in\rng(\Gamma(\iota))$. We also have
\begin{equation*}
\supp^J_{\mathbf Y}(Y(\tau))=\supp^J_{\mathbf Y}(\alpha+1,0)=\supp^\Gamma_{\mathbf Y}(0)=\emptyset\subseteq\mathbf X.
\end{equation*}
If $\tau$ does not lie in the range of~$\mathbf P$, then we have $Y(\tau)=(\alpha,\tau)$ for some~$\alpha<\nu$. In this case we get $\supp^J_{\mathbf Y}(Y(\tau))=\supp^\Gamma_{\mathbf Y}(\tau)$, so that the open equivalence coincides with the support property of the dilator~$\Gamma$. Thus the equivalence from the beginning of the proof is established. For $s\in\Gamma(\mathbf L)$ with $\mathbf L:=\mathbf Y+S^{\mathbf P}_{\Gamma(\mathbf Y)}$ we now observe
\begin{equation*}
s\in\rng(\Gamma(\iota+S_{\Gamma(\iota)}))\quad\Leftrightarrow\quad\supp^\Gamma_{\mathbf L}(s)\subseteq\rng(\iota+S_{\Gamma(\iota)}),
\end{equation*}
also by the support condition for~$\Gamma$. Furthermore, we compute
\begin{align*}
{\supp^{\Gamma\circ\mathbf S}_{\mathbf Y}}\circ\Gamma(\mathbf Y+S_Y)(s)&=\bigcup\{\supp^{\mathbf S}_{\mathbf Y}(\rho)\,|\,\rho\in{\supp^\Gamma_{\mathbf S(\mathbf Y)}}\circ\Gamma(\mathbf Y+S_Y)(s)\}\\
&=\bigcup\{{\supp^{\mathbf S}_{\mathbf Y}}\circ(\mathbf Y+S_Y)(\tau)\,|\,\tau\in\supp^\Gamma_{\mathbf L}(s)\}.
\end{align*}
Using the equivalence from the beginning of the proof, one can now derive
\begin{equation*}
s\in\rng(\Gamma(\iota+S_{\Gamma(\iota)}))\quad\Leftrightarrow\quad{\supp^{\Gamma\circ\mathbf S}_{\mathbf Y}}\circ\Gamma(\mathbf Y+S_Y)(s)\subseteq\mathbf X=\rng(\iota).
\end{equation*}
Even though we will not use this, we note that this step corresponds to the fact that $\Gamma$ preserves pullbacks. It is straightforward to derive the lemma: Given~$y\in\mathbf X$, we write $\pi_{\mathbf Y}\circ\iota(y)=(\alpha,t)$. The definition of $\mathbf X\subseteq\mathbf Y$ yields $\supp^{\Gamma\circ\mathbf S}_{\mathbf Y}(t)\subseteq\mathbf X$ as well as $t=\Gamma(\mathbf Y+S_Y)(s)$ for some $s\in\Gamma(\mathbf L)$. By the equivalence above, we can conclude that $s=\Gamma(\iota+S_{\Gamma(\iota)})(r)$ holds for some $r\in\Gamma(\mathbf K)=\mathbf O$. We thus get
\begin{equation*}
t=\Gamma(\mathbf Y+S_Y)\circ\Gamma(\iota+S_{\Gamma(\iota)})(r)=\Gamma\left(\iota+S_{Y\circ\Gamma(\iota)}\right)(r).
\end{equation*}
So $\pi_{\mathbf Y}\circ\iota(y)=(\alpha,t)$ is the image of $(\alpha,r)$ under $\nu\times\Gamma(\iota+S_{Y\circ\Gamma(\iota)})$.
\end{proof}

The following completes the constructions that were sketched at the beginning of the present section. We point out that $\pi_{\mathbf X}$ is analogous to the dashed arrow from the diagramm before Lemma~\ref{def:XsubY}.

\begin{definition}\label{def:pi_X}
Invoking Lemma~\ref{lem:def-pi_X}, let $\pi_{\mathbf X}$ be the unique embedding such that
\begin{equation*}
\begin{tikzcd}[column sep=large]
\mathbf X\arrow[r,"\pi_{\mathbf X}"]\arrow[d,swap,"\iota"] & \nu\times\mathbf O\arrow[d,"\nu\times\Gamma(\iota+S_{Y\circ\Gamma(\iota)})"]\\
\mathbf Y\arrow[r,"\pi_{\mathbf Y}"] & \nu\times(\Gamma\circ\mathbf S)(\mathbf Y)
\end{tikzcd}
\end{equation*}
is a commutative diagram. To define a partial function $\psi^{\mathbf X}:\nu\times\mathbf O\pto\mathbf X$ that is surjective and order preserving, we put
\begin{equation*}
\psi^{\mathbf X}_\alpha s:=\psi^{\mathbf X}(\alpha,s):=\begin{cases}
x & \text{if $\pi_{\mathbf X}(x)=(\alpha,s)$},\\
\text{undefined} & \text{if $(\alpha,s)\notin\rng(\pi_{\mathbf X})$}.
\end{cases}
\end{equation*}
We will write $\dom(\psi^{\mathbf X}):=\rng(\pi_{\mathbf X})$ for the domain of this partial function. Also, let $I:\mathbf X\to\mathbf X+S^{\mathbf R}_{\Gamma(\mathbf X)}=\mathbf K$ with $I(x):=x$ be the map onto the first summand.
\end{definition}

Crucially, the search tree $S^{\mathbf R}_{\Gamma(\mathbf X)}$ depends on an embedding $\mathbf R:\nu\to\Gamma(\mathbf X)$ that has a meaningful connection to the collapsing function~$\psi^{\mathbf X}$.

\begin{lemma}\label{lem:psi-X_R}
We have $(\alpha,0)\in\dom(\psi^{\mathbf X})$ and $\gamma_{\mathbf X}(\psi^{\mathbf X}_{\alpha+1}0)=\mathbf R(\alpha)$ for all $\alpha<\nu$.
\end{lemma}
\begin{proof}
Note that we have distinct elements $0\in\Gamma(\mathbf K)=\mathbf O$ and $0\in\mathbf \Gamma(\mathbf S(\mathbf Y))$. In view of Definitions~\ref{def:Gamma-dilator} and~\ref{def:Beth-O}, we get
\begin{align*}
\left(\nu\times\Gamma(\iota+S_{Y\circ\Gamma(\iota)})\right)(\alpha+1,0)&=(\alpha+1,0)=\pi_{\mathbf Y}\circ\iota\circ\mathbf R_0(\alpha)\\
{}&=\left(\nu\times\Gamma(\iota+S_{Y\circ\Gamma(\iota)})\right)\circ\pi_X\circ\mathbf R_0(\alpha).
\end{align*}
This entails $(\alpha+1,0)=\pi_{\mathbf X}\circ\mathbf R_0(\alpha)\in\rng(\pi_{\mathbf X})=\dom(\psi^{\mathbf X})$ and $\psi^{\mathbf X}_{\alpha+1}0=\mathbf R_0(\alpha)$, so that we get $\gamma_{\mathbf X}(\psi^{\mathbf X}_{\alpha+1}0)=\gamma_{\mathbf X}\circ\mathbf R_0(\alpha)=\mathbf R(\alpha)$. To show $(\alpha,0)\in\dom(\psi^{\mathbf X})$ with~$\alpha$ at the place of~$\alpha+1$, use Lemmas~\ref{lem:Omega-alpha-defined} and~\ref{lem:Omega-alpha-X} to write $(\alpha,0)=\pi_{\mathbf Y}\circ\iota(x)$ with~$x\in\mathbf X$. Then argue as before, with $\alpha$ and $y$ at the place of $\alpha+1$ and~$\mathbf R_0(\alpha)$.
\end{proof}

In the rest of this section we characterize the range of~$\pi_{\mathbf X}$ or, in other words, the domain of the partial function~$\psi^{\mathbf X}$. As a first step, we assign supports to the elements of~$\mathbf K$ and~$\mathbf O$. To avoid misunderstanding, we point out that the following support functions do not belong to a dilator. Let us also recall that $\supp^S$ was defined in the paragraph before Definition~\ref{def:search}.

\begin{definition}\label{def:supp-O}
Let $\supp^{\mathbf K}:\mathbf K=\mathbf X+S^{\mathbf R}_{\Gamma(\mathbf X)}\to[\mathbf X]^{<\omega}$ be given by
\begin{equation*}
\supp^{\mathbf K}(x)=\{x\},\quad\supp^{\mathbf K}(\mathbf X+\sigma)=\bigcup\{\supp^\Gamma_{\mathbf X}(\rho)\,|\,\rho\in\supp^S_{\Gamma(\mathbf X)}(\sigma)\backslash\rng(\mathbf R)\}.
\end{equation*}
Furthermore, define $\supp^{\mathbf O}:\mathbf O=\Gamma(\mathbf K)\to[\mathbf X]^{<\omega}$ by setting
\begin{equation*}
\supp^{\mathbf O}(s)=\bigcup\{\supp^{\mathbf K}(\rho)\,|\,\rho\in\supp^\Gamma_{\mathbf K}(s)\}.
\end{equation*}
\end{definition}

The given definition -- and in particular the exclusion of $\rng(\mathbf R)$ -- is justified by the following connection with the support functions of our dilators~$\mathbf S$ and $\Gamma\circ\mathbf S$.

\begin{lemma}\label{lem:relate-supports}
Each of the diagrams
\begin{equation*}
\begin{tikzcd}[column sep=large]
\mathbf K\ar[r,"\iota+S_{Y\circ\Gamma(\iota)}"]\ar[d,swap,"\supp^{\mathbf K}"] & \mathbf S(\mathbf Y)\ar[d,"\supp^{\mathbf S}_{\mathbf Y}"]\\
{[\mathbf X]^{<\omega}}\ar[r,"{[\iota]^{<\omega}}"] & {[\mathbf Y]^{<\omega}}
\end{tikzcd}
\qquad\text{and}\qquad
\begin{tikzcd}
\mathbf O\ar[rr,"\Gamma(\iota+S_{Y\circ\Gamma(\iota)})"]\ar[d,swap,"\supp^{\mathbf O}"] && \Gamma\circ\mathbf S(\mathbf Y)\ar[d,"\supp^{\Gamma\circ\mathbf S}_{\mathbf Y}"]\\
{[\mathbf X]^{<\omega}}\ar[rr,"{[\iota]^{<\omega}}"] && {[\mathbf Y]^{<\omega}}
\end{tikzcd}
\end{equation*}
commutes.
\end{lemma}
\begin{proof}
Let us abbreviate $f:=\iota+S_{Y\circ\Gamma(\iota)}$. Using Definitions~\ref{def:S^E} and~\ref{def:X+S^E} as well as the naturality of supports, we get
\begin{align*}
\supp^{\mathbf S}_{\mathbf Y}\circ f(\mathbf X+\sigma)&=\bigcup\{\supp^J_{\mathbf Y}(Y\circ\Gamma(\iota)(\rho))\,|\,\rho\in\supp^S_{\Gamma(\mathbf X)}(\sigma)\},\\
[\iota]^{<\omega}\circ\supp^{\mathbf K}(\mathbf X+\sigma)&=\bigcup\{\supp^\Gamma_{\mathbf Y}(\Gamma(\iota)(\rho))\,|\,\rho\in\supp^S_{\Gamma(\mathbf X)}(\sigma)\backslash\rng(\mathbf R)\}.
\end{align*}
To see why the range of~$\mathbf R$ is excluded, note that $\rho=\mathbf R(\alpha)$ entails $\Gamma(\iota)(\rho)=\mathbf P(\alpha)$, so that Definition~\ref{def:Y-Gamma-J} yields $Y\circ\Gamma(\iota)(\rho)=(\alpha+1,0)$ and thus
\begin{equation*}
\supp^J_{\mathbf Y}(Y\circ\Gamma(\iota)(\rho))=\supp^{\Gamma}_{\mathbf Y}(0)=\emptyset.
\end{equation*}
As a straightforward consequence, the left diagram commutes if we have
\begin{equation*}
\supp^J_{\mathbf Y}(Y\circ\Gamma(\iota)(\rho))=\supp^\Gamma_{\mathbf Y}(\Gamma(\iota)(\rho))\quad\text{for}\quad\rho\in\Gamma(\mathbf X)\backslash\rng(\mathbf R).
\end{equation*}
Even though we do not need this fact, it is instructive to observe that the equation fails for $\rho=\mathbf R(\alpha)=\gamma_{\mathbf X}(\psi^{\mathbf X}_{\alpha+1}0)$, where Lemma~\ref{lem:Gamma-normal} yields
\begin{equation*}
\supp^\Gamma_{\mathbf Y}(\Gamma(\iota)(\rho))=[\iota]^{<\omega}\left(\supp^\Gamma_{\mathbf X}(\gamma_{\mathbf X}(\psi^{\mathbf X}_{\alpha+1}0))\right)=[\iota]^{<\omega}\left(\{\psi^{\mathbf X}_{\alpha+1}0\}\right)\neq\emptyset.
\end{equation*}
On the other hand, $\rho\notin\rng(\mathbf R)$ entails $\Gamma(\iota)(\rho)\notin\rng(\mathbf P)$, as we have $\Gamma(\iota)\circ\mathbf R=\mathbf P$ and $\Gamma(\iota)$ is injective. We then get $Y\circ\Gamma(\iota)(\rho)=(\alpha,\Gamma(\iota)(\rho))$ for some~$\alpha<\nu$. In this case, the desired equality is immediate by the definition of the support for~$J$. The right diagram is readily reduced to the left one.
\end{proof}

Our well founded `subterm' relation on~$\mathbf Y$ can now be transferred to~$\mathbf X$.

\begin{lemma}\label{lem:wf-X}
For any $x\in X$ and $(\alpha,s)\in\dom(\psi^{\mathbf X})$ we have
\begin{equation*}
\iota(x)\tl\iota(\psi^{\mathbf X}_\alpha s)\quad\Leftrightarrow\quad x\in\supp^{\mathbf O}(s),
\end{equation*}
where $\tl$ is the well founded relation on~$\mathbf Y$ that was specified in Definition~\ref{def:XsubY}.
\end{lemma}
\begin{proof}
When $\psi^{\mathbf X}_\alpha s$ is defined, we have $\pi_{\mathbf X}(\psi^{\mathbf X}_\alpha s)=(\alpha,s)$ and hence
\begin{equation*}
\pi_{\mathbf Y}\circ\iota(\psi^{\mathbf X}_\alpha s)=\left(\nu\times\Gamma(\iota+S_{Y\circ\Gamma(\iota)})\right)\circ\pi_{\mathbf X}(\psi^{\mathbf X}_\alpha s)=\left(\alpha,\Gamma(\iota+S_{Y\circ\Gamma(\iota)})(s)\right).
\end{equation*}
Together with the previous lemma, it follows that $\iota(x)\tl\iota(\psi^{\mathbf X}_\alpha s)$ amounts to
\begin{align*}
\iota(x)\in&{\supp^{\Gamma\circ\mathbf S}_{\mathbf Y}}\circ\Gamma(\iota+S_{Y\circ\Gamma(\iota)})(s)={}\\
&\bigcup\left\{\supp^{\mathbf S}_{\mathbf Y}(\tau)\,\left|\,\tau\in{\supp^\Gamma_{\mathbf S(\mathbf Y)}}\circ\Gamma(\iota+S_{Y\circ\Gamma(\iota)})(s)\right.\right\}=\\
&\bigcup\left\{{\supp^{\mathbf S}_{\mathbf Y}}\circ(\iota+S_{Y\circ\Gamma(\iota)})(\rho)\,\left|\,\rho\in\supp^\Gamma_{\mathbf K}(s)\right.\right\}=\\
&\bigcup\left\{[\iota]^{<\omega}\left(\supp^{\mathbf K}(\rho)\right)\,\left|\,\rho\in\supp^\Gamma_{\mathbf K}(s)\right.\right\}=[\iota]^{<\omega}\left(\supp^{\mathbf O}(s)\right),
\end{align*}
which is clearly equivalent to $x\in\supp^{\mathbf O}(s)$.
\end{proof}

Given that $\iota:X\to Y$ is an inclusion map, we will also refer to $\tl$ as a well founded relation on~$\mathbf X$. The following definition uses recursion along this relation. It also exploits that any element of~$\mathbf X$ can be uniquely written as $\psi^{\mathbf X}_\alpha s$, since the partial function $\psi^{\mathbf X}:\nu\times\mathbf O\pto\mathbf X$ is surjective and order preserving. When we refer to $\psi^{\mathbf X}_\alpha s$ as a given element of~$\mathbf X$, we always assume~$(\alpha,s)\in\dom(\psi^{\mathbf X})$.

\begin{definition}\label{def:K_beta}
For $\gamma<\nu$ we define $K^-_\gamma:\mathbf X\to[\mathbf O]^{<\omega}$ and $K_\gamma:\mathbf O\to[\mathbf O]^{<\omega}$ by
\begin{align*}
K^-_\gamma(\psi^{\mathbf X}_\alpha s)&:=\begin{cases}
\{s\}\cup K_\gamma(s) & \text{if $\gamma\leq\alpha$},\\
\emptyset & \text{otherwise},
\end{cases}\\
K_\gamma(s)&:=\bigcup\{K^-_\gamma(x)\,|\,x\in\supp^{\mathbf O}(s)\}.
\end{align*}
\end{definition}

As promised, we can now characterize the domain of our collapsing function.

\begin{proposition}\label{prop:dom_psi^X}
For any $\gamma<\nu$ and $s\in\mathbf O$ we have
\begin{equation*}
(\gamma,s)\in\dom(\psi^{\mathbf X})=\rng(\pi_{\mathbf X})\quad\Leftrightarrow\quad K_\gamma(s)\subseteq_{\mathbf O} s.
\end{equation*}
\end{proposition}
\begin{proof}
Let $G_\gamma$ and $G^{\Gamma\circ\mathbf S}_\gamma$ be the maps that arise from Definition~\ref{def:nu-collapse} in conjunction with Assumption~\ref{ass:collapse-search}. We abbreviate $f:=\iota+S_{Y\circ\Gamma(\iota)}:\mathbf K\to\mathbf S(\mathbf Y)$ and show that
\begin{equation*}
\begin{tikzcd}[column sep=large]
\mathbf X\ar[r,"K^-_\gamma"]\ar[d,swap,"\iota"] & {[\mathbf O]^{<\omega}}\ar[d,"{[\Gamma(f)]^{<\omega}}"] & \mathbf O\ar[l,swap,"K_\gamma"]\ar[d,"\Gamma(f)"]\\
\mathbf Y\ar[r,"G^{\Gamma\circ\mathbf S}_\gamma"] & {[\Gamma\circ\mathbf S(\mathbf Y)]^{<\omega}} & \Gamma\circ\mathbf S(\mathbf Y)\ar[l,swap,"G_\gamma"]
\end{tikzcd}
\end{equation*}
is commutative. To prove that the left quare commutes, we employ induction over the well founded relation from Lemma~\ref{lem:wf-X}. For the induction step, recall that the proof of Lemma~\ref{lem:wf-X} yields $\pi_{\mathbf Y}\circ\iota(\psi^{\mathbf X}_\alpha s)=(\alpha,\Gamma(f)(s))$. By Definition~\ref{def:nu-collapse} we get
\begin{equation*}
G^{\Gamma\circ\mathbf S}_\gamma\circ\iota(\psi^{\mathbf X}_\alpha s)=\begin{cases}
\{\Gamma(f)(s)\}\cup G_\gamma\circ\Gamma(f)(s) & \text{if $\alpha\geq\gamma$},\\
\emptyset & \text{otherwise}.
\end{cases}
\end{equation*}
To complete the induction step, use the hypothesis and Lemma~\ref{lem:relate-supports} to compute
\begin{align*}
[\Gamma(f)]^{<\omega}\circ K_\gamma(s)&=\bigcup\{[\Gamma(f)]^{<\omega}\circ K^-_\gamma(x)\,|\,x\in\supp^{\mathbf O}(s)\}\\
&=\bigcup\{G^{\Gamma\circ\mathbf S}_\gamma(x)\,|\,x\in[\iota]^{<\omega}\circ\supp^{\mathbf O}(s)\}\\
&=\bigcup\{G^{\Gamma\circ\mathbf S}_\gamma(x)\,|\,\tau\in{\supp^{\Gamma\circ\mathbf S}_{\mathbf Y}}\circ\Gamma(f)(s)\}=G_\gamma\circ\Gamma(f)(s).
\end{align*}
Note that this proves that the right square commutes. Definition~\ref{def:nu-collapse} does now yield
\begin{align*}
K_\gamma(s)\subseteq_{\mathbf O} s\quad&\Leftrightarrow\quad G_\gamma\circ\Gamma(f)(s)=[\Gamma(f)]^{<\omega}\circ K_\gamma(s)\subseteq_{\Gamma\circ\mathbf S(\mathbf Y)}\Gamma(f)(s)\\
{}&\Leftrightarrow\quad(\gamma,\Gamma(f)(s))\in\rng(\pi_{\mathbf Y}).
\end{align*}
To complete the proof, we show that $(\gamma,\Gamma(f)(s))\in\rng(\pi_{\mathbf Y})$ and $(\gamma,s)\in\rng(\pi_{\mathbf X})$ are equivalent, which means that the diagram from Definition~\ref{def:pi_X} is a pullback. Concerning the easier direction, we note that $(\gamma,s)=\pi_{\mathbf X}(x)$ entails
\begin{equation*}
(\gamma,\Gamma(f)(s))=(\nu\times\Gamma(f))\circ\pi_X(x)=\pi_{\mathbf Y}\circ\iota(x)\in\rng(\pi_{\mathbf Y}).
\end{equation*}
To prove the converse, we assume $(\gamma,\Gamma(f)(s))=\pi_{\mathbf Y}(y)$ and derive $y\in\mathbf X=\rng(\iota)$. In view of $f=(\mathbf Y+S_Y)\circ(\iota+S_{\Gamma(\iota)})$ we set $t:=\Gamma(\iota+S_{\Gamma(\iota)})(s)$ to obtain
\begin{equation*}
\pi_{\mathbf Y}(y)=(\gamma,\Gamma(\mathbf Y+S_Y)(t))\in\rng(\nu\times\Gamma(\mathbf Y+S_Y)).
\end{equation*}
The proof of Lemma~\ref{lem:def-pi_X} shows that $t\in\rng(\Gamma(\iota+S_{\Gamma(\iota)}))$ entails
\begin{equation*}
\supp^{\Gamma\circ\mathbf S}_{\mathbf Y}\left(\Gamma(\mathbf Y+S_Y)(t)\right)\subseteq\mathbf X.
\end{equation*}
By Definition~\ref{def:XsubY} we now get $y\in\mathbf X$, as desired.
\end{proof}

Let us also record a basic observation that will be needed later:

\begin{lemma}\label{lem:supp-O}
We have ${\supp^{\mathbf O}}\circ\Gamma(I)=\supp^\Gamma_\mathbf X$.
\end{lemma}
\begin{proof}
First recall that ${\supp^\Gamma_{\mathbf K}}\circ\Gamma(I)=[I]^{<\omega}\circ\supp^\Gamma_{\mathbf X}$ holds by naturality. In view of Definition~\ref{def:supp-O} we have ${\supp^{\mathbf K}}\circ I(x)=\{x\}$ and thus
\begin{align*}
{\supp^{\mathbf O}}\circ\Gamma(I)(\rho)&=\bigcup\{\supp^{\mathbf K}(\tau)\,|\,\tau\in{\supp^\Gamma_{\mathbf K}}\circ\Gamma(I)(\rho)\}\\
{}&=\bigcup\{{\supp^{\mathbf K}}\circ I(x)\,|\,x\in\supp^\Gamma_{\mathbf X}(\rho)\}=\supp^\Gamma_{\mathbf X}(\rho),
\end{align*}
as desired.
\end{proof}

We conclude this section with an observation about the order on~$\mathbf O$.

\begin{lemma}\label{lem:psi-less-Gamma}
We have $\Gamma(I)(s)<_{\mathbf O}\Gamma_{\mathbf X+\sigma}$ for all $s\in\Gamma(\mathbf X)$ and $\sigma\in S^{\mathbf R}_{\Gamma(\mathbf X)}$.
\end{lemma}
\begin{proof}
As $I$ maps into the first summand of $\mathbf X+ S^{\mathbf R}_{\Gamma(\mathbf X)}$, we see that $\Gamma_{\mathbf X+\sigma}$ lies outside the range of~$\Gamma(I)$. But the latter is an initial segment of~$\mathbf O$, by Corollary~\ref{cor:normal-initial}.
\end{proof}

\section{Operator control and infinite proofs}\label{sect:operators}

From the previous section we have a function
\begin{equation*}
\psi^{\mathbf X}:\nu\times\mathbf O=\nu\times\Gamma\left(\mathbf X+S^{\mathbf R}_{\Gamma(\mathbf X)}\right)\pto\mathbf X
\end{equation*}
that is surjective and order preserving but partial, i.\,e., not always defined. In the present section, we transform $\psi^{\mathbf X}$ into a function $\psi:\nu\times\mathbf O\to\mathbf O$ that is total but not always order preserving. We then define an abstract variant of the operator controlled proofs that have been introduced by Buchholz~\cite{buchholz-local-predicativity}. Finally, we construct an operator controlled proof that embeds the search tree $S^{\mathbf R}_{\Gamma(\mathbf X)}$ from Section~\ref{sect:search-trees}.

As a first step, we transform $\psi^{\mathbf X}$ into a function~$\psi^+$ that remains partial but has co\-domain~$\mathbf O$. Note that the following definition composes arrows from the diagram at the beginning of Section~\ref{sect:search-to-notations}. This diagram commutes by Lemma~\ref{lem:Gamma-normal}, which means that $\gamma_{\mathbf K}\circ I$ equals $\Gamma(I)\circ\gamma_{\mathbf X}$. The maps $\gamma_Z:Z\to\Gamma(Z)$ and $I:\mathbf X\to\mathbf K$ are given by Definitions~\ref{def:Gamma-normal} and~\ref{def:pi_X}, while $\supp^{\mathbf O}:\mathbf O\to[\mathbf X]^{<\omega}$ comes from Definition~\ref{def:supp-O}.

\begin{definition}\label{def:psi+}
The partial function $\psi^+:\nu\times\mathbf O\to_p\mathbf O$ is given by
\begin{equation*}
\psi^+_\alpha s:=\psi^+(\alpha,s):=\begin{cases}
\gamma_{\mathbf K}\circ I(\psi^{\mathbf X}_\alpha s) & \text{if $(\alpha,s)\in\dom(\psi^{\mathbf X})=:\dom(\psi^+)$},\\
\text{undefined} & \text{otherwise}.
\end{cases}
\end{equation*}
To define $\supp^+:\mathbf O\to[\mathbf O]^{<\omega}$, we set $\supp^+:=[\gamma_{\mathbf K}\circ I]^{<\omega}\circ\supp^{\mathbf O}$.
\end{definition}

For an arbitrary dilator~$D$, no family of embeddings $Z\to D(Z)$ needs to exist. This explains why Definition~\ref{def:nu-collapse} involves two families of functions $G^D_\gamma$ and $G_\gamma$ with domain $X$ and $D(X)$, respectively. In Definition~\ref{def:K_beta} we have constructed corresponding functions $K^-_\gamma:\mathbf X\to[\mathbf O]^{<\omega}$ and $K_\gamma:\mathbf O\to[\mathbf O]^{<\omega}$. In the present case, however, we do have an embedding $\gamma_{\mathbf X}\circ I:\mathbf X\to\mathbf O$ (amongst others because of the maps $\gamma_Z:Z\to\Gamma(Z)$ that make $\Gamma$ normal). As the following shows, this allows us to eliminate $K^-_\gamma$ in favour of~$K_\gamma$. Similarly, the functions $G^D_\gamma$ and $G_\gamma$ are unified in traditional ordinal notation systems, as we have seen in Example~\ref{ex:Buchholz-psi}.

\begin{proposition}\label{prop:elim-K_beta^-}
For any $\gamma<\nu$ and $(\alpha,s)\in\dom(\psi^+)$ we have
\begin{equation*}
K_\gamma(\psi^+_\alpha s)=\begin{cases}
\{s\}\cup K_\gamma(s) &\text{if $\gamma\leq\alpha$},\\
\emptyset & \text{otherwise}.
\end{cases}
\end{equation*}
Furthermore, we have $K_\gamma(t)=\bigcup\{K_\gamma(r)\,|\,r\in\supp^+(t)\}$ for any $t\in\mathbf O$.
\end{proposition}
\begin{proof}
The naturality of supports and Lemma~\ref{lem:Gamma-normal} yield
\begin{equation*}
{\supp^\Gamma_{\mathbf K}}\circ\Gamma(I)\circ\gamma_{\mathbf X}(\psi^{\mathbf X}_\alpha s)=[I]^{<\omega}\circ{\supp^\Gamma_{\mathbf X}}\circ\gamma_{\mathbf X}(\psi^{\mathbf X}_\alpha s)=\{I(\psi^{\mathbf X}_\alpha s)\}.
\end{equation*}
In view of Definition~\ref{def:supp-O}, we can derive
\begin{multline*}
\supp^{\mathbf O}(\psi^+_\alpha s)={\supp^{\mathbf O}}\circ\gamma_{\mathbf K}\circ I(\psi^{\mathbf X}_\alpha s)={\supp^{\mathbf O}}\circ\Gamma(I)\circ\gamma_{\mathbf X}(\psi^{\mathbf X}_\alpha s)\\
{}=\bigcup\{\supp^{\mathbf K}(\rho)\,|\,\rho\in{\supp^\Gamma_{\mathbf K}}\circ\Gamma(I)\circ\gamma_{\mathbf X}(\psi^{\mathbf X}_\alpha s)\}={\supp^{\mathbf K}}\circ I(\psi^{\mathbf X}_\alpha s)=\{\psi^{\mathbf X}_\alpha s\}.
\end{multline*}
For later reference, we record that this entails
\begin{equation*}
\supp^+(\psi^+_\alpha s)=\{\gamma_{\mathbf K}\circ I(\psi^{\mathbf X}_\alpha s)\}=\{\psi^+_\alpha s\}.
\end{equation*}
Due to Definition~\ref{def:K_beta} we obtain
\begin{equation*}
K_\gamma(\psi^+_\alpha s)=\bigcup\{K^-_\gamma(x)\,|\,x\in\supp^{\mathbf O}(\psi^+_\alpha s)\}=K_\gamma^-(\psi^{\mathbf X}_\alpha s).
\end{equation*}
The first claim of the proposition is now immediate by Definition~\ref{def:K_beta}. In the paragraph before this definition, we have observed that any element $x\in\mathbf X$ can be written as $x=\psi^{\mathbf X}_\alpha s$ for some $(\alpha,s)\in\dom(\psi^{\mathbf X})$. We get $\gamma_{\mathbf K}\circ I(x)=\psi^+_\alpha s$, which means that the previous observation can be reformulated as
\begin{equation*}
K_\gamma\circ\gamma_{\mathbf K}\circ I=K^-_\gamma.
\end{equation*}
In view of Definition~\ref{def:K_beta}, we can deduce
\begin{equation*}
K_\gamma(t)=\bigcup\{K_\gamma^-(x)\,|\,x\in\supp^{\mathbf O}(t)\}=\bigcup\{K_\gamma(r)\,|\,r\in[\gamma_{\mathbf K}\circ I]^{<\omega}\circ\supp^{\mathbf O}(t)\}.
\end{equation*}
Considering the definition of~$\supp^+$, this coincides with the remaining claim.
\end{proof}

The following result will be used to extend~$\psi^+$ into a total function.

\begin{proposition}\label{prop:psi-X_total}
Given any $\alpha<\nu$ and $s\in\mathbf O$, we get $(\alpha,t)\in\dom(\psi^+)$ for some element $t\in\{s\}\cup K_\alpha(s)$ with $s\leq_{\mathbf O}t$.
\end{proposition}
\begin{proof}
The main task will be to show that $r\in K_\alpha(s)$ entails $r\notin K_\alpha(r)\subseteq K_\alpha(s)$. Once this is achieved, we can conclude by induction on the cardinality of the finite set~$K_\alpha(s)$. Indeed, for $K_\alpha(s)\subseteq_{\mathbf O} s$ we get $(\alpha,s)\in\dom(\psi^+)$ by Proposition~\ref{prop:dom_psi^X}, so we can take $t=s$. If $K_\alpha(s)\subseteq_{\mathbf O} s$ fails, we can pick an $r\in K_\alpha(s)$ with $s\leq r$. By the initial claim, $K_\alpha(r)$ has fewer elements than~$K_\alpha(s)$. Inductively, we thus get $(\alpha,t)\in\dom(\psi^+)$ for some $t\in\{r\}\cup K_\alpha(r)\subseteq K_\alpha(s)$ with $s\leq r\leq t$. To prove the initial claim, recall that Lemma~\ref{lem:wf-X} provides a well founded relation~$\tl$ on~$\mathbf X\subseteq\mathbf Y$. It will be convenient to consider the associated height function $h:\mathbf X\to\mathbb N$ with
\begin{equation*}
h\left(\psi^{\mathbf X}_\gamma t\right)=\max\left(\{0\}\cup\{h(z)+1\,|\,z\in\supp^{\mathbf O}(t)\}\right).
\end{equation*}
Aiming at $r\notin K_\alpha(r)$, we fix an arbitrary element $x\in\supp^{\mathbf O}(r)$. We use induction on $h(y)\leq h(x)$ to prove $r\notin K^-_\alpha(y)$. Writing $y=\psi^{\mathbf X}_\gamma t$, we note that $h(x)\geq h(y)$ forces $x\notin\supp^{\mathbf O}(t)$ and hence $r\neq t$. With the induction hypothesis, this yields
\begin{equation*}
r\notin\{t\}\cup\bigcup\{K^-_\alpha(z)\,|\,z\in\supp^{\mathbf O}(t)\}=\{t\}\cup K_\alpha(t)\supseteq K^-_\alpha(\psi^{\mathbf X}_\gamma t)=K^-_\alpha(y).
\end{equation*}
Since $x\in\supp^{\mathbf O}(r)$ was arbitrary, we get
\begin{equation*}
r\notin\bigcup\{K^- _\alpha(x)\,|\,x\in\supp^{\mathbf O}(r)\}=K_\alpha(r).
\end{equation*}
Another induction on~$h(x)$ shows that $r\in K^-_\alpha(x)$ entails $K_\alpha(r)\subseteq K^-_\alpha(x)$. It is straightfoward to conclude that $r\in K_\alpha(s)$ entails $K_\alpha(r)\subseteq K_\alpha(s)$.
\end{proof}

We can now define the total extension of~$\psi^{\mathbf X}$ that was promised above.

\begin{definition}\label{def:psi-C}
To obtain a total function $\psi:\nu\times\mathbf O\to\mathbf O$, we put
\begin{multline*}
\psi_\alpha s:=\psi(\alpha,s):=\psi^+_\alpha t\quad\text{for the $<_{\mathbf O}$-minimal $t\in\{s\}\cup K_\alpha(s)$ with}\\
\text{$s\leq_{\mathbf O}t$ and $(\alpha,t)\in\dom(\psi^+)$}.
\end{multline*}
Let us also define $C_\alpha(t):=\{s\in\mathbf O\,|\,K_\alpha(s)\subseteq_{\mathbf O}t\}$ for all $\alpha<\nu$ and $t\in\mathbf O$.
\end{definition}

Note that we immediately get $\psi_\alpha s=\psi^+_\alpha s$ for $(\alpha,s)\in\dom(\psi^+)$. The sets~$C_\alpha(t)$ and the following proposition evoke traditional constructions of ordinal notation systems in terms of set theory (see e.\,g.~\cite[Definition~4.2]{buchholz-local-predicativity}). In contrast to these constructions, our functions $\psi_\alpha$ do not seem to be weakly increasing. Indeed, if we have $t<r<t'$ with $(\alpha,r)\in\dom(\psi^+)$ but $\psi_\alpha t=\psi^+_\alpha t'$ due to $r\notin K_\alpha(t)$, then we get $\psi_\alpha r=\psi^+_\alpha r<\psi^+_\alpha t'=\psi_\alpha t$. At the same time, Corollary~\ref{cor:psi-C-monotone} will ensure that the order is preserved in relevant cases.

\begin{proposition}\label{prop:C-psi}
The following holds for all $\alpha<\nu$ and $s,t\in\mathbf O$:
\begin{enumerate}[label=(\alph*)]
\item Given $s\in C_\alpha(t)$ with $s<t$, we get $\psi_\gamma s\in C_\alpha(t)$ for any~$\gamma<\nu$.
\item If we have $s<\psi_{\alpha+1}0$, then $s\in C_\alpha(t)$ implies~$s<\psi_\alpha t$.
\item If we have $t\in C_\alpha(t)$, then $s<\psi_\alpha t$ implies $s\in C_\alpha(t)$.
\end{enumerate}
\end{proposition}
\begin{proof}
(a) For $\gamma<\alpha$ we have $K_\alpha(\psi_\gamma s)=\emptyset$, so that $\psi_\gamma s\in C_\alpha(t)$ is immediate. Let us now assume $\gamma\geq\alpha$. With $h:\mathbf X\to\mathbb N$ as in the proof of Proposition~\ref{prop:psi-X_total}, an easy induction on~$h(x)$ yields $K^-_\gamma(x)\subseteq K^-_\alpha(x)$ and simultaneously $K_\gamma(s)\subseteq K_\alpha(s)$. We note that this entails $C_\alpha(t)\subseteq C_\gamma(t)$. Given that we have $s\in C_\alpha(t)$ and~$s<t$, we learn that $\psi_\gamma s=\psi^+_\gamma t'$ holds for some
\begin{equation*}
t'\in\{s\}\cup K_\gamma(s)\subseteq\{s\}\cup K_\alpha(s)\subseteq_{\mathbf O} t.
\end{equation*}
As in the proof of Proposition~\ref{prop:psi-X_total}, we get $K_\alpha(t')\subseteq K_\alpha(s)$ and hence
\begin{equation*}
K_\alpha(\psi_\gamma s)=K_\alpha(\psi^+_\gamma t')=\{t'\}\cup K_\alpha(t')\subseteq\{t'\}\cup K_\alpha(s)\subseteq_{\mathbf O} t.
\end{equation*}
This amounts to $\psi_\gamma s\in C_\alpha(t)$, as desired.

(b) We use induction on the build-up of $s\in\Gamma(\mathbf K)$ according to Definition~\ref{def:Gamma(X)}. In the crucial case, we have $s=\Gamma_z$ for some $z\in\mathbf K=\mathbf X+S^{\mathbf R}_{\Gamma(\mathbf X)}$. As Lemma~\ref{lem:psi-X_R} ensures $(\alpha+1,0)\in\dom(\psi^+)$, the assumption $s<\psi_{\alpha+1}0$ yields
\begin{equation*}
\gamma_{\mathbf K}(z)=s<_{\mathbf O}\psi_{\alpha+1}0=\psi^+_{\alpha+1}0=\gamma_{\mathbf K}\circ I(\psi^{\mathbf X}_{\alpha+1}0).
\end{equation*}
The range of~$I:\mathbf X\to\mathbf K$ is an initial segment, so $z=I(x)$ holds for some $x<\psi^{\mathbf X}_{\alpha+1}0$. Like any other element of~$\mathbf X$, the latter can be written in the form $x=\psi^{\mathbf X}_\gamma r$, which yields $s=\psi^+_\gamma r$. We must have $\gamma\leq\alpha$, as $\psi^{\mathbf X}$ is order preserving. If we have $\gamma<\alpha$, then $s<\psi_\alpha t$ is immediate. Let us now assume~$\gamma=\alpha$. We then have $r\in K_\alpha(s)$, so that $s\in C_\alpha(t)$ yields $r<t$. For the appropriate $t'\geq t$ we get
\begin{equation*}
s=\psi^+_\alpha r<\psi^+_\alpha t'=\psi_\alpha t.
\end{equation*}
In the case of a term $s=\overline\varphi s_0s_1$, we recall that Definition~\ref{def:Gamma-dilator} yields
\begin{equation*}
\supp^\Gamma_{\mathbf K}(s)=\supp^\Gamma_{\mathbf K}(s_0)\cup\supp^\Gamma_{\mathbf K}(s_1).
\end{equation*}
The equality remains valid when we replace $\supp^\Gamma_{\mathbf K}$ by $\supp^{\mathbf O}$ or $\supp^+$ or $K_\alpha$, due to Definition~\ref{def:supp-O} and Proposition~\ref{prop:elim-K_beta^-}. So $s\in C_\alpha(t)$ is equivalent to $s_0,s_1\in C_\alpha(t)$. Also note that $s<\psi_{\alpha+1}0$ and $s_0,s_1<\psi_{\alpha+1}0$ are equivalent by Definition~\ref{def:Gamma(X)}, as
\begin{equation*}
\psi_{\alpha+1}0\in\rng(\gamma_{\mathbf K})=\{\Gamma_z\,|\,z\in\mathbf K\}
\end{equation*}
is strongly critical. We can thus invoke the induction hypothesis to get $s_0,s_1<\psi_\alpha t$. The latter entails $s<\psi_\alpha t$, because $\psi_\alpha t$ is strongly critical as well. For a term of the form $s=\langle s_0,\ldots,s_{n-1}\rangle$, the argument is similar.

(c) As in the proof of~(b), we argue by induction on the build-up of $s\in\Gamma(\mathbf K)$. Let us first assume that we have $s=\Gamma_z$ for some~$z\in\mathbf K$. Given $s<\psi_\alpha t<\psi_{\alpha+1}0$, we can once again write $s=\psi^+_\gamma r$ with $\gamma\leq\alpha$. If the last inequality is strict, we obtain $K_\alpha(s)=\emptyset$, so that $s\in C_\alpha(t)$ is immediate. Now assume $\gamma=\alpha$ and recall that $(\alpha,r)\in\dom(\psi^+)$ entails $K_\alpha(r)\subseteq_{\mathbf O}r$. Given $t\in C_\alpha(t)$, we have $\psi_\alpha t=\psi^+_\alpha t$, so that $s<\psi_\alpha t$ entails $r<t$. Together we get
\begin{equation*}
K_\alpha(s)=\{r\}\cup K_\alpha(r)\subseteq_{\mathbf O}t
\end{equation*}
and hence $s\in C_\alpha(t)$, as desired. Let us also consider a term $s=\overline\varphi s_0s_1<\psi_\alpha t$. For each $i\leq 1$ we get $s_i<\psi_\alpha t$, so that the induction hypothesis yields $s_i\in C_\alpha(t)$. We can conclude $s\in C_\alpha(t)$, as noted in the proof of~(b). An analogous argument applies in the case of a term $s=\langle s_0,\ldots,s_{n-1}\rangle$ with $n>1$. For $s=0$, it suffices to observe that $K_\alpha(0)$ is empty, since the same holds for~$\supp^\Gamma_{\mathbf K}(0)$.
\end{proof}

As observed in part~(b) of the previous proof, all values $\psi_\alpha t$ are strongly critical. The next result provides inequalities between different values of~$\psi$.

\begin{corollary}\label{cor:psi-C-monotone}
The following holds for all $s,t\in\mathbf O$:
\begin{enumerate}[label=(\alph*)]
\item For $t\neq 0$ we have $\psi_\alpha 0<\psi_\alpha t<\psi_{\alpha+1} 0=\Gamma(I)\circ\mathbf R(\alpha)$.
\item If we have $s\in C_\alpha(t)$, then $s<t$ implies $\psi_\alpha s<\psi_\alpha t$.
\end{enumerate}
\end{corollary}
\begin{proof}
Concerning part~(a), let us first observe that Lemmas~\ref{lem:Gamma-normal} and~\ref{lem:psi-X_R} yield
\begin{equation*}
\Gamma(I)\circ\mathbf R(\alpha)=\Gamma(I)\circ\gamma_{\mathbf X}(\psi^{\mathbf X}_{\alpha+1})=\gamma_{\mathbf K}\circ I(\psi^{\mathbf X}_{\alpha+1}0)=\psi^+_{\alpha+1}0=\psi_{\alpha+1}0.
\end{equation*}
The second inequality in part~(a) is immediate, while the first one reduces to~(b), as $\supp^\Gamma_{\mathbf K}(0)=\emptyset$ entails $K_\alpha(0)=\emptyset$ and hence $0\in C_\alpha(t)$. Let us now establish part~(b). Given $s\in C_\alpha(t)$ and $s<t$, we get $\psi_\alpha s\in C_\alpha(t)$ by part~(a) of the previous proposition. Part~(b) of the latter yields $\psi_\alpha s<\psi_\alpha t$, as we have $\psi_\alpha s<\psi_{\alpha+1}0$.
\end{proof}

With the sets $C_\alpha(t)$ at hand, we can recover the operators $\mathcal H_s$ of Buchholz~\cite{buchholz-local-predicativity}.

\begin{definition}\label{def:opertors-H}
For $s\in\mathbf O$ and $a\in[\mathbf O]^{<\omega}$ we set
\begin{equation*}
\mathcal H_s(a):=\bigcap\{C_\alpha(t)\,|\,\alpha<\nu\text{ and }t\in\mathbf O\text{ with }s<_{\mathbf O}t\text{ and }a\subseteq C_\alpha(t)\}\subseteq\mathbf O.
\end{equation*}
\end{definition}

Note that the intersection is taken over a non-empty family, because $a\subseteq C_\alpha(t)$ amounts to $b\subseteq_{\mathbf O}t$ for the finite set $b=\bigcup_{r\in a}K_\alpha(r)$. The following is immediate.

\begin{lemma}\label{lem:H-closure}
The following holds for all $s,t\in\mathbf O$ and $a,b\in[\mathbf O]^{<\omega}$:
\begin{enumerate}[label=(\alph*)]
\item We have $a\subseteq\mathcal H_s(a)$.
\item Given $a\subseteq\mathcal H_s(b)$, we get $\mathcal H_s(a)\subseteq\mathcal H_s(b)$.
\item For $s<t$ we have $\mathcal H_s(a)\subseteq\mathcal H_t(a)$.
\end{enumerate}
\end{lemma}

Parts~(a) and~(b) express that $\mathcal H_s$ is a closure operator. Together, they ensure that $a\subseteq b$ implies $\mathcal H_s(a)\subseteq\mathcal H_s(b)$. As we will see, the following is an abstract way to say that $\mathcal H_s$ is nice in the sense of~\cite[Definition~3.5]{buchholz-local-predicativity}.

\begin{proposition}\label{prop:supp-H}
For all $s,t\in\mathbf O$ and $a\in[\mathbf O]^{<\omega}$ we have
\begin{equation*}
s\in\mathcal H_t(a)\quad\Leftrightarrow\quad\supp^+(s)\subseteq\mathcal H_t(a).
\end{equation*}
\end{proposition}
\begin{proof}
For each $\alpha<\nu$, Proposition~\ref{prop:elim-K_beta^-} yields
\begin{equation*}
K_\alpha(s)\subseteq_{\mathbf O}t\quad\Leftrightarrow\quad K_\alpha(r)\subseteq_{\mathbf O}t\text{ for all }r\in\supp^+(s).
\end{equation*}
Hence the equivalence from the proposition holds with $C_\alpha(t)$ at the place of~$\mathcal H_t(a)$. This pointwise version is stronger than the claim itself.
\end{proof}

The corollary below encapsulates various closure properties, such as
\begin{equation*}
\overline\varphi r_0r_1\in\mathcal H_r(a)\quad\Leftrightarrow\quad\{r_0,r_1\}\subseteq\mathcal H_r(a).
\end{equation*}
In view of Definition~\ref{def:Gamma-dilator}, the direction from right to left follows from the corollary for $s_i=r_i$ and $t_0=\overline\varphi r_0r_1$ (with $m=2$ and $n=1$). The converse direction follows when we take $s_0=\overline\varphi r_0r_1$ and $t_i=r_i$. We get an analogous equivalence for terms of the form~$\langle r_0,\ldots,r_{k-1}\rangle$. Due to Lemma~\ref{lem:Veblen-addition-support}, we also learn that $\{s_0,s_1\}\subseteq\mathcal H_r(a)$ entails $s_0+s_1\in\mathcal H_r(a)$ and $\varphi s_0s_1\in\mathcal H_r(a)$, where $\varphi$ is our total extension of~$\overline\varphi$. One can also take $m=0$, to obtain $0,1\in\mathcal H_r(a)$ from $\supp^\Gamma_{\mathbf K}(0)=\emptyset$.

\begin{corollary}\label{cor:H-nice}
Consider any $s_0,\ldots,s_{m-1}$ and $t_0,\ldots,t_{n-1}$ in~$\mathbf O$. If we have
\begin{equation*}
\textstyle\bigcup_{i<m}\supp^\Gamma_{\mathbf K}(s_i)\supseteq\textstyle\bigcup_{j<n}\supp^\Gamma_{\mathbf K}(t_j),
\end{equation*}
then $\{s_0,\ldots,s_{m-1}\}\subseteq\mathcal H_r(a)$ implies $\{t_0,\ldots,t_{n-1}\}\subseteq\mathcal H_r(a)$.
\end{corollary}
\begin{proof}
As in the proof of Proposition~\ref{prop:C-psi}, the given inclusion remains valid when we replace $\supp^\Gamma_{\mathbf K}$ by $\supp^+$. We can conclude by the previous proposition.
\end{proof}

The following result on collapsing functions (cf.~\cite[Lemma~4.6]{buchholz-local-predicativity}) completes our list of closure properties. In particular, it yields $\psi_\alpha 0\in\mathcal H_t(a)$ for all~$\alpha<\nu$.

\begin{corollary}\label{cor:operator-collapsing}
Given $s\in\mathcal H_t(a)$ with $s\leq_{\mathbf O}t$, we get $\psi_\alpha s\in\mathcal H_t(a)$ for all $\alpha<\nu$.
\end{corollary}
\begin{proof}
To obtain $\psi_\alpha s\in\mathcal H_t(a)$, we need to establish $\psi_\alpha s\in C_\beta(t')$ for arbitrary $\beta<\nu$ and $t'>t$ with $a\subseteq C_\beta(t')$. The assumption $s\in\mathcal H_t(\alpha)$ ensures $s\in C_\beta(t')$. Given that we have $s\leq t<t'$, Proposition~\ref{prop:C-psi} yields $\psi_\alpha s\in C_\beta(t')$, as required.
\end{proof}

The rest of this section concerns a notion of infinite proof that is heavily inspired by work of Buchholz~\cite{buchholz-local-predicativity}. As preparation, we introduce notation that relates to the parameters and the rank of formulas. In Section~\ref{sect:search-trees} and Definition~\ref{def:psi+}, we have explained $\supp^{\mathbf L}_{\Gamma(\mathbf X)}(a)\in[\Gamma(\mathbf X)]^{<\omega}$ and $\supp^+(s)\in[\mathbf O]^{<\omega}$ for $a\in\mathbf L^u_{\Gamma(\mathbf X)}$ and~$s\in\mathbf O$, respectively. The following definition overloads this notation by admitting arguments of different types. To interpret the notation correctly, one will need to infer the type of the argument from the context.

\begin{definition}\label{def:supp-extend}
For an $\mathbf L^u_{\Gamma(\mathbf X)}$-formula~$\varphi$ and an $\mathbf L^u_{\Gamma(\mathbf X)}$-sequent $\Gamma$, we put
\begin{align*}
\supp^{\mathbf L}_{\Gamma(\mathbf X)}(\varphi)&:=\bigcup\{\supp^{\mathbf L}_{\Gamma(\mathbf X)}(a)\,|\,a\in\mathbf L^u_{\Gamma(\mathbf X)}\text{ is a parameter of }\varphi\},\\
\supp^{\mathbf L}_{\Gamma(\mathbf X)}(\Gamma)&:=\textstyle\bigcup_{i<n}\supp^{\mathbf L}_{\Gamma(\mathbf X)}(\varphi_i)\quad\text{for}\quad\Gamma=\varphi_0,\ldots,\varphi_{n-1}.
\end{align*}
When $\sigma$ is an element of~$\mathbf L^u_{\Gamma(\mathbf X)}$, an $\mathbf L^u_{\Gamma(\mathbf X)}$-formula or an $\mathbf L^u_{\Gamma(\mathbf X)}$-sequent, we define
\begin{equation*}
\supp^+(\sigma):=[\Gamma(I)]^{<\omega}\circ\supp^{\mathbf L}_{\Gamma(\mathbf X)}(\sigma)\in[\mathbf O]^{<\omega}.
\end{equation*}
For $\alpha<\nu$, an $\mathbf L^u_{\Gamma(\mathbf X)}$-formula~$\varphi$ is called a $\Sigma(\alpha)$-formula if all universal quantifiers in~$\varphi$ are bounded and we have
\begin{equation*}
\supp^{\mathbf L}_{\Gamma(\mathbf X)}(\varphi)\subseteq_{\Gamma(\mathbf X)}\mathbf R(\alpha).
\end{equation*}
Let us also agree to abbreviate $L[\alpha]:=L^u_{\mathbf R(\alpha)}\in\mathbf L^u_{\Gamma(\mathbf X)}$ for~$\alpha<\nu$. 
\end{definition}

To motivate the new notation, we recall that Definition~\ref{def:search} involves relativized axioms $\ax_n^{L(i)}$ with $L(i)=L[\alpha]$ for $\alpha=\nu_i$. We are particularly interested in the case of $\Delta_0$-collection, where~$\ax_n^{L[\alpha]}$ has instances of the form
\begin{equation*}
\forall x\in a_0\exists y\in L[\alpha]\,\theta(x,y,a_1,\ldots,a_n)\to\exists w\in L[\alpha]\forall x\in a_0\exists y\in w\,\theta(x,y,a_1,\ldots,a_n).
\end{equation*}
In the relevant cases, we will have $a_i\in\mathbf L^u_{\Gamma(\mathbf X)}$ and $\supp^{\mathbf L}_{\Gamma(\mathbf X)}(a_i)\subseteq_{\Gamma(\mathbf X)}\mathbf R(\alpha)$. On an intuitive level, this means that the parameters come from the $\mathbf R(\alpha)$-th stage of the constructible hierarchy, i.\,e., from $L[\alpha]$. The given condition ensures that
\begin{equation*}
\varphi:=\forall x\in a_0\exists y\,\theta(x,y,a_1,\ldots,a_n)
\end{equation*}
is a $\Sigma(\alpha)$-formula. Our instance of $\Delta_0$-collection can now be written as
\begin{equation*}
\varphi^{L[\alpha]}\to\exists w\in L[\alpha].\,\varphi^w.
\end{equation*}
For an arbitrary $\Sigma(\alpha)$-formula, this implication can be deduced from $\Delta_0$-collection in~$L[\alpha]$, at least for the actual constructible hierarchy (see~\cite[Theorem~I.4.3]{barwise-admissible}). This fact will not be used in the following, but it does explain the role of~$\Sigma(\alpha)$-formulas.

As a final ingredient for our infinite proofs, we assign formula ranks that will be used to control cut inferences. In order to explain the following definition, we recall that $\mathbf L^u_{\Gamma(\mathbf X)}$ is built over a set~$u\ni 0$ of urelements (fixed in Assumption~\ref{ass:u}). According to Section~\ref{sect:search-trees}, our $\mathbf L^u_{\Gamma(\mathbf X)}$-formulas are closed (unless noted otherwise) and in negation normal form. The required ordinal arithmetic on~$\mathbf O=\Gamma(\mathbf K)$ was discussed at the end of Section~\ref{sect:booting-up}. Let us point out that $\Gamma(I):\Gamma(\mathbf X)\to\mathbf O$ commutes with basic ordinal arithmetic. It follows that all ranks lie in the range of~$\Gamma(I)$. For notational reasons, it will still be convenient to have ranks in~$\mathbf O$ rather than~$\Gamma(\mathbf X)$.

\begin{definition}\label{def:ranks}
The function $\rk:\mathbf L^u_{\Gamma(\mathbf X)}\to\mathbf O$ is given by
\begin{gather*}
\rk(w):=0\text{ for }w\in u,\qquad \rk(L^u_s):=\omega\cdot(1+\Gamma(I)(s)),\\
\rk(\{x\in L^u_s\,|\,\varphi(x,a_1,\ldots,a_n)\}):=\rk(L^u_s)+1.
\end{gather*}
To each bounded $\mathbf L^u_{\Gamma(\mathbf X)}$-formula~$\varphi$, we assign a rank $\rk(\varphi)\in\mathbf O$ by setting
\begin{align*}
\rk(a\in b)&:=\rk(\neg\, a\in b):=\max\{\rk(a)+6,\rk(b)+1\},\\
\rk(a=b)&:=\rk(\neg\, a=b):=\max\{\rk(a),\rk(b),5\}+4,\\
\rk(\varphi_0\lor\varphi_1)&:=\rk(\varphi_0\land\varphi_1):=\max\{\rk(\varphi_0),\rk(\varphi_1)\}+1,\\
\rk(\exists x\in a.\,\varphi(x))&:=\rk(\forall x\in a.\,\varphi(x)):=\max\{\rk(a),\rk(\varphi(0))+2\}.
\end{align*}
\end{definition}

Note that we get $\rk(\varphi)=\rk(\neg\varphi)$ for any bounded $\mathbf L^u_{\Gamma(\mathbf X)}$-formula~$\varphi$, because of our treatment of negation as a defined operation. Let us record a basic property:

\begin{lemma}\label{lem:ranks-operators}
For all $b\in\mathbf L^u_{\Gamma(\mathbf X)}$ and $t\in\mathbf O$ we have $\rk(b)\in\mathcal H_0(\supp^+(b))$ and
\begin{equation*}
\supp^+(b)\subseteq_{\mathbf O}t\quad\Leftrightarrow\quad\rk(b)<_{\mathbf O}\omega\cdot(1+t).
\end{equation*}
Both properties remain valid when we replace~$b$ by a bounded $\mathbf L^u_{\Gamma(\mathbf X)}$-formula~$\varphi$.
\end{lemma}
\begin{proof}
For $b\in u$ it suffices to observe $\rk(b)=0$ and $\supp^+(b)=\emptyset$. In the remaining cases, the equivalence holds since we have $\rk(b)=\omega\cdot(1+\Gamma(I)(s))+i$ for some~$i\leq 1$, where $s$ is the largest element of $\supp^{\mathbf L}_{\Gamma(\mathbf X)}(b)$. We also get
\begin{equation*}
\Gamma(I)(s)\in\supp^+(b)\subseteq\mathcal H_0(\supp^+(b)).
\end{equation*}
In view of $1=\varphi_00$, Lemma~\ref{lem:Veblen-addition-support} yields
\begin{equation*}
\supp^\Gamma_{\mathbf K}(\rk(b))=\supp^\Gamma_{\mathbf K}(\omega\cdot(1+\Gamma(I)(s))+i)\subseteq\supp^\Gamma_{\mathbf K}(\Gamma(I)(s)).
\end{equation*}
Thus $\rk(b)\in\mathcal H_0(\supp^+(b))$ follows by Corollary~\ref{cor:H-nice}. A straightforward induction over formulas shows that we can write $\rk(\varphi)=\rk(b)+n$ with $n\in\mathbb N$, where~$b$ is a parameter of~$\varphi$ or equal to $0\in u\subseteq\mathbf L^u_{\Gamma(\mathbf X)}$. In both cases we get
\begin{equation*}
\rk(\varphi)=\rk(b)+n\in\mathcal H_0(\supp^+(b))\subseteq\mathcal H_0(\supp^+(\varphi))
\end{equation*}
due to Corollary~\ref{cor:H-nice} and Lemma~\ref{lem:H-closure}. By another induction over formulas, we see that $\rk(b)\leq\rk(\varphi)$ holds for any parameter~$b$ of the formula $\varphi$. Given that $r<\omega\cdot s$ entails $r+n<\omega\cdot s$, this ensures that the equivalence remains valid.
\end{proof}

To justify the focus on bounded formulas, we recall that any $\mathbf L^u_{\Gamma(\mathbf X)}$-formula~$\varphi$ is associated with a disjunction $\bigvee_{a\in\iota(\varphi)}\varphi_a$ or conjunction $\bigwedge_{a\in\iota(\varphi)}\varphi_a$, as explained in Section~\ref{sect:search-trees}. If $\varphi$ is bounded, so is $\varphi_a$ for every~$a\in\iota(\varphi)$, due to~\cite[Definition~3.12]{freund-equivalence}. Thus all formulas in Definition~\ref{def:search} are bounded, and the same will hold for the formulas in our infinite proofs. We say that an $\mathbf L^u_{\Gamma(\mathbf X)}$-sequent is bounded if it consists of bounded formulas only. The assignment of ranks is designed to validate the following, which is shown in the proof of~\cite[Theorem~3.14]{freund-equivalence} (see also~\cite[Lemma~3]{buchholz-local-predicativity}).

\begin{lemma}\label{lem:rank-disj-conj}
Given any bounded $\mathbf L^u_{\Gamma(\mathbf X)}$-formula $\varphi$, we have
\begin{equation*}
\rk(\varphi_a)<_{\mathbf O}\rk(\varphi)\quad\text{for all $a\in\iota(\varphi)=\iota_{\Gamma(\mathbf X)}(\varphi)$}.
\end{equation*}
\end{lemma}

In the paragraph before Lemma~\ref{lem:def-pi_X}, we have observed that $\mathbf O$ is well founded, which justifies the following recursion. Intuitively, we have $(r,a)\vdash_s^t\Gamma$ if the sequent $\Gamma$ has an infinite proof with height at most~$t$, where $\mathcal H_r(a)$ and~$s$ control relevant parameters and cuts. The given definition is inspired by~\cite[Theorem~3.8]{buchholz-local-predicativity}.

\begin{definition}\label{def:H-controlled-deriv}
By recursion on~$t$, we declare that the relation
\begin{equation*}
(r,a)\vdash^t_s\Gamma
\end{equation*}
between elements $r,s,t\in\mathbf O$, $a\in[\mathbf O]^{<\omega}$ and a bounded $\mathbf L^u_{\Gamma(\mathbf X)}$-sequent~$\Gamma$ holds precisely if we have
\begin{equation*}
\{t\}\cup\supp^+(\Gamma)\subseteq\mathcal H_r(a)
\end{equation*}
and one of the following clauses applies:
\begin{enumerate}[label=(\roman*),topsep=5pt,itemsep=5pt]
\item for some conjunctive $\varphi\simeq\bigwedge_{b\in\iota(\varphi)}\varphi_b\in\Gamma$ and every $b\in\iota(\varphi)\subseteq\mathbf L^u_{\Gamma(\mathbf X)}$, there is a $t(b)<t$ such that we have $(r,a\cup\supp^+(b))\vdash^{t(b)}_s\Gamma,\varphi_b$,
\item for some disjunctive $\varphi\simeq\bigvee_{b\in\iota(\varphi)}\varphi_b\in\Gamma$ and some $b\in\iota(\varphi)\subseteq\mathbf L^u_{\Gamma(\mathbf X)}$ such that we have $\supp^+(b)\subseteq_{\mathbf O} t$, there is a $t(0)<t$ with $(r,a)\vdash^{t(0)}_s\Gamma,\varphi_b$,
\item for some bounded $\mathbf L^u_{\Gamma(\mathbf X)}$-formula~$\psi$ with $\rk(\psi)<s$, there is a $t(0)<t$ such that we have $(r,a)\vdash^{t(0)}_s\Gamma,\psi$ and $(r,a)\vdash^{t(0)}_s\Gamma,\neg\psi$,
\item for some $\alpha<\nu$ and some $\Sigma(\alpha)$-formula~$\varphi$ with $\exists z\in L[\alpha].\,\varphi^z\in\Gamma$, there is an element $t(0)<t$ with $(r,a)\vdash^{t(0)}_s\Gamma,\varphi^{L[\alpha]}$.
\end{enumerate}
\end{definition}

Sometimes one wants to apply the given clauses in a modified form, e.\,g., to derive $(r,a)\vdash_s^t\Gamma_0,\Gamma_1$ from $(r,a)\vdash_s^{t(0)}\Gamma_0,\psi$ and $(r,a)\vdash_s^{t(1)}\Gamma_1,\neg\psi$ with $t(0)\neq t(1)<t$. This is possible due to the following standard result (cf.~\cite[Lemma~3.9(a)]{buchholz-local-predicativity}).

\begin{lemma}[Weakening]\label{lem:weakening}
Given $r\leq r',s\leq s',t\leq t'$ and $a\subseteq\mathcal H_{r'}(a')$, we have
\begin{equation*}
(r,a)\vdash^t_s\Gamma\quad\text{and}\quad\{t'\}\cup\supp^+(\Delta)\subseteq\mathcal H_{r'}(a')\qquad\Rightarrow\qquad(r',a')\vdash^{t'}_{s'}\Delta,\Gamma.
\end{equation*}
\end{lemma}
\begin{proof}
One argues by induction on~$t\in\mathbf O$ and distinguishes cases that correspond to the clauses from Definition~\ref{def:H-controlled-deriv}. In each case, one uses the induction hypothesis and reapplies the same clause. This is possible because Lemma~\ref{lem:H-closure} yields
\begin{equation*}
a\cup c\subseteq\mathcal H_r(a\cup c)\subseteq\mathcal H_{r'}(a'\cup c),
\end{equation*}
where one takes $c=\supp^+(b)$ for clause~(i) and $c=\emptyset$ in the other cases.
\end{proof}

We always refer to the lemma as `weakening', even when $a'$ is a proper subset~of~$a$, where we get an apparent strengthening. In the following result, the bound $\omega\cdot\rk(\varphi)$ could be improved to~$2\cdot\rk(\varphi)$. We keep the suboptimal bound because only $t\mapsto\omega\cdot t$ has been defined in the present paper.

\begin{lemma}\label{lem:a=a}
For any bounded $\mathbf L^u_{\Gamma(\mathbf X)}$-formula~$\varphi$ and any $a\in\mathbf L^u_{\Gamma(\mathbf X)}$ we have
\begin{equation*}
(0,\supp^+(\varphi))\vdash^{\omega\cdot\rk(\varphi)}_0\varphi,\neg\varphi\qquad\text{and}\qquad(0,\supp^+(a))\vdash^{\omega\cdot\rk(a)+2}_0 a=a.
\end{equation*}
\end{lemma}
\begin{proof}
To establish the first claim, we argue by induction on~$\rk(\varphi)$. First observe that $\mathcal H_0(\supp^+(\varphi))$ contains $\rk(\varphi)$ and hence also $\omega\cdot\rk(\varphi)$, due to Lemma~\ref{lem:ranks-operators} and its proof. As disjunction and conjunction are dual (see \cite[Definition~3.12]{freund-equivalence}), we may assume $\varphi\simeq\bigvee_{b\in\iota(\varphi)}\varphi_b$ to get $\neg\varphi\simeq\bigwedge_{b\in\iota(\varphi)}\neg\varphi_b$, or in other words $\iota(\neg\varphi)=\iota(\varphi)$ and $\neg(\varphi_b)=(\neg\varphi)_b$. In view of Lemma~\ref{lem:rank-disj-conj}, we use the induction hypothesis to get
\begin{equation*}
(0,\supp^+(\varphi_b))\vdash^{\omega\cdot\rk(\varphi_b)}_0\varphi_b,\neg\varphi_b\qquad\text{for each }b\in\iota(\varphi).
\end{equation*}
To prepare an application of weakening, we observe that \cite[Definition~3.12]{freund-equivalence} yields
\begin{equation*}
\supp^{\mathbf L}_{\Gamma(\mathbf X)}(\varphi_b)\subseteq\supp^{\mathbf L}_{\Gamma(\mathbf X)}(\varphi)\cup\supp^{\mathbf L}_{\Gamma(\mathbf X)}(b).
\end{equation*}
This inclusion remains valid when we apply $[\Gamma(I)]^{<\omega}$ to both sides, i.\,e., when we replace $\supp^{\mathbf L}_{\Gamma(\mathbf X)}$ by $\supp^+$. For each $b\in\iota(\varphi)$, we can use Lemma~\ref{lem:ranks-operators} to derive
\begin{equation*}
t(b):=\max\{\rk(\varphi_b),\rk(b)\}\in\mathcal H_0(\supp^+(\varphi)\cup\supp^+(b)).
\end{equation*}
As announced, we now apply weakening to get
\begin{equation*}
(0,\supp^+(\varphi)\cup\supp^+(b))\vdash_0^{\omega\cdot t(b)}\varphi,\neg\varphi,\varphi_b,\neg\varphi_b.
\end{equation*}
The choice of~$t(b)$ and Lemma~\ref{lem:ranks-operators} ensure $\supp^+(b)\subseteq_{\mathbf O}\omega\cdot t(b)+1$, as required in clause~(ii) of Definition~\ref{def:H-controlled-deriv}. By the latter, we thus obtain
\begin{equation*}
(0,\supp^+(\varphi)\cup\supp^+(b))\vdash_0^{\omega\cdot t(b)+1}\varphi,\neg\varphi,\neg\varphi_b\qquad\text{for each }b\in\iota(\varphi)=\iota(\neg\varphi).
\end{equation*}
Based on~\cite[Definition~3.12]{freund-equivalence} and Lemma~\ref{lem:ranks-operators}, it is not hard to check that $b\in\iota(\varphi)$ entails $\rk(b)<\rk(\varphi)$, so that we get $t(b)<\rk(\varphi)$ by Lemma~\ref{lem:rank-disj-conj}. We can thus apply clause~(i) of Definition~\ref{def:H-controlled-deriv}, in order to complete the proof of the first claim from the lemma. To derive the second claim, we show
\begin{equation*}
(0,\supp^+(a))\vdash_0^{\omega\cdot\rk(a)+1}\forall x\in a.\,x\in a
\end{equation*}
by induction on~$\rk(a)$. Let us consider a term of the form $a=\{x\in L^u_s\,|\,\theta(x,\mathbf d)\}$. For $a\in u$ and $a=L^u_s$ the argument is easier (but note that $a\in u$ leads to the bound $\omega\cdot\rk(a)+1=1$ rather than $\omega\cdot\rk(a)=0$). By~\cite[Definition~3.12]{freund-equivalence} we have
\begin{gather*}
\forall x\in a.\,x\in a\simeq\textstyle\bigwedge_{b\in\iota}\neg\theta(b,\mathbf d)\lor b\in a\quad\text{and}\quad b\in a\simeq\textstyle\bigvee_{c\in\iota}\theta(c,\mathbf d)\land c=b\\
\text{with }\iota=\{b\in\mathbf L^u_{\Gamma(\mathbf X)}\,|\,\supp^{\mathbf L}_{\Gamma(\mathbf X)}(b)\subseteq_{\Gamma(\mathbf X)}s\}.
\end{gather*}
In the clause for $b\in a$, we will take $c$ to be the same term as~$b$. To derive $b=b$, we recall the general clause
\begin{equation*}
(b_0=b_1)\simeq\textstyle\bigwedge_{i\in\{0,1\}}\forall x\in b_i.\,x\in b_{1-i}.
\end{equation*}
When $b_0$ and $b_1$ are the same term~$b$, then the two conjuncts coincide, but we still need a step to introduce the conjunction. So the induction hypothesis and clause~(i) of Definition~\ref{def:H-controlled-deriv} yield
\begin{equation*}
(0,\supp^+(b))\vdash_0^{\omega\cdot\rk(b)+2}b=b.
\end{equation*}
This shows the second claim of the lemma, once the present induction is completed. We have $\supp^+(\theta(b,\mathbf d))\subseteq\supp^+(a)\cup\supp^+(b)$, and Lemma~\ref{lem:ranks-operators} provides
\begin{equation*}
s(b)<\rk(a)\quad\text{for}\quad s(b):=\max\{\rk(b)+1,\rk(\theta(b,\mathbf d))\}.
\end{equation*}
Using the first part of the present lemma, we can thus derive
\begin{equation*}
(0,\supp^+(a)\cup\supp^+(b))\vdash^{\omega\cdot s(b)+1}_0\neg\theta(b,\mathbf d),\theta(b,\mathbf d)\land b=b.
\end{equation*}
We now use clause~(ii) of Definition~\ref{def:H-controlled-deriv} three times, once to get $b\in a$ and twice to combine the disjuncts, so that we obtain
\begin{equation*}
(0,\supp^+(a)\cup\supp^+(b))\vdash^{\omega\cdot s(b)+4}_0\neg\theta(b,\mathbf d)\lor b\in a.
\end{equation*}
To complete the induction step, one applies clause~(i) of the same definition.
\end{proof}

In the rest of this section, we show how the search tree $S^{\mathbf R}_{\Gamma(\mathbf X)}$ from Definition~\ref{def:search} can be transformed into an infinite proof. We begin with the crucial axioms.

\begin{proposition}\label{prop:ax1+n}
For each of the $\Delta_0$-collection axioms $\ax_{1+n}$ from Definition~\ref{def:enum-ax} and any $\alpha<\nu$, we have
\begin{equation*}
(0,\emptyset)\vdash_0^t\ax_{1+n}^{L[\alpha]}\quad\text{with}\quad t:=\psi_{\alpha+1}0+\omega\cdot 3.
\end{equation*}
\end{proposition}
\begin{proof}
Corollaries~\ref{cor:H-nice} and~\ref{cor:operator-collapsing} provide $\psi_{\alpha+1}0+\omega\cdot m+n\in\mathcal H_0(\emptyset)$ for $m,n\in\mathbb N$. We recall $L[\alpha]=L^u_{\mathbf R(\alpha)}$ and $\psi_{\alpha+1}0=\Gamma(I)\circ\mathbf R(\alpha)$ as well as $\supp^{\mathbf L}_{\Gamma(\mathbf X)}(L^u_s)=\{s\}$. The initial condition from Definition~\ref{def:H-controlled-deriv} can now be derived as
\begin{equation*}
\supp^+(\ax_{1+n}^{L[\alpha]})=[\Gamma(I)]^{<\omega}\circ\supp^{\mathbf L}_{\Gamma(\mathbf X)}(L^u_{\mathbf R(\alpha)})=\{\psi_{\alpha+1}0\}\subseteq\mathcal H_0(\emptyset).
\end{equation*}
As in the  paragraph that follows Definition~\ref{def:supp-extend}, we write collection in the form
\begin{equation*}
\ax_{1+n}=\forall\mathbf z\forall v\,(\psi\to\exists w\,\psi^w)\quad\text{with}\quad\psi(v,\mathbf z):=\forall x\in v\exists y\,\theta(x,y,\mathbf z),
\end{equation*}
for a $\Delta_0$-formula~$\theta$ and variables $\mathbf z=z_1,\ldots,z_k$. Note that we get
\begin{equation*}
\ax_{1+n}^{L[\alpha]}=\forall z_1\in L[\alpha]\ldots\forall z_k\in L[\alpha]\forall v\in L[\alpha]\,(\psi^{L[\alpha]}\to\exists w\in L[\alpha].\,\psi^w).
\end{equation*}
Let us now recall that~\cite[Definition~3.12]{freund-equivalence} yields
\begin{equation*}
\forall y\in L[\alpha].\,\varphi(y)\,\simeq\,\textstyle\bigwedge_{a\in\iota}\varphi(a)\quad\text{for}\quad\iota:=\{a\in\mathbf L^u_{\Gamma(\mathbf X)}\,|\,\supp^{\mathbf L}_{\Gamma(\mathbf X)}(a)\subseteq_{\Gamma(\mathbf X)}\mathbf R(\alpha)\}.
\end{equation*}
To conclude by $k+1$ applications of clause~(i) from Definition~\ref{def:H-controlled-deriv}, we shall thus show the following: For $s:=\psi_{\alpha+1}0+\omega\cdot 2+3$ and arbitrary $a_0,\ldots,a_k\in\iota$, we have
\begin{equation*}
(0,\supp^+(\varphi))\vdash^s_0 \varphi^{L[\alpha]}\to\exists w\in L[\alpha].\,\varphi^w\quad\text{with}\quad\varphi:=\psi(a_0,\ldots,a_k).
\end{equation*}
In the proof of Proposition~\ref{prop:C-psi} we have observed that $\psi_{\alpha+1}0$ is strongly critical. This justifies the last step in the computation
\begin{equation*}
\rk(L[\alpha])=\omega\cdot(1+\Gamma(I)\circ\mathbf R(\alpha))=\omega\cdot(1+\psi_{\alpha+1}0)=\psi_{\alpha+1}0.
\end{equation*}
By Definition~\ref{def:ranks} in conjunction with Lemma~\ref{lem:ranks-operators}, we get $\rk(\varphi^{L[\alpha]})=\psi_{\alpha+1}0+2$. We can thus use Lemma~\ref{lem:a=a} to obtain
\begin{equation*}
(0,\supp^+(\varphi^{L[\alpha]}))\vdash_0^r\neg\varphi^{L[\alpha]},\varphi^{L[\alpha]}\quad\text{with}\quad r:=\omega\cdot(\psi_{\alpha+1}0+2)=\psi_{\alpha+1}0+\omega\cdot 2.
\end{equation*}
Weakening allows us to replace $\supp^+(\varphi^{L[\alpha]})$ by $\supp^+(\varphi)$, as we have
\begin{equation*}
\supp^+(\varphi^{L[\alpha]})\subseteq\supp^+(\varphi)\cup\supp^+(L[\alpha])\subseteq\mathcal H_0(\supp^+(\varphi)).
\end{equation*}
Now $\varphi$ is a $\Sigma(\alpha)$-formula, due to~$a_i\in\iota$. Thus clause~(iv) of Definition~\ref{def:H-controlled-deriv} yields
\begin{equation*}
(0,\supp^+(\varphi))\vdash_0^{r+1}\neg\varphi^{L[\alpha]},\exists w\in L[\alpha].\,\varphi^w.
\end{equation*}
From Section~\ref{sect:search-trees} we recall that $\varphi^{L[\alpha]}\to\exists w\in L[\alpha].\,\varphi^w$ and $\neg\varphi^{L[\alpha]}\lor\exists w\in L[\alpha].\,\varphi^w$ denote the same formula in negation normal form. We can thus conclude by two applications of clause~(ii) from Definition~\ref{def:H-controlled-deriv}.
\end{proof}

On an intuitive level, the following holds because the stage $\mathbf R(\alpha)$ of $L[\alpha]$ is a limit (in fact $\Gamma(I)\circ\mathbf R(\alpha)=\psi_{\alpha+1}$ is strongly critical).

\begin{proposition}\label{prop:ax0}
Consider the axiom $\ax_0=\forall x\exists y.\, x\in y$ from Definition~\ref{def:enum-ax}. For any $\alpha<\nu$ we have $(0,\emptyset)\vdash^t_0\ax_0^{L[\alpha]}$ with $t:=\psi_{\alpha+1}0$.
\end{proposition}
\begin{proof}
First note that we have
\begin{equation*}
\{t\}\cup\supp^+(\ax_0^{L[\alpha]})=\{\psi_{\alpha+1}0\}\subseteq\mathcal H_0(\emptyset),
\end{equation*}
as in the previous proof. To conclude by clauses~(i) and~(ii) of Definition~\ref{def:H-controlled-deriv}, we write $\iota=\{a\in\mathbf L^u_{\Gamma(\mathbf X)}\,|\,\supp^{\mathbf L}_{\Gamma(\mathbf X)}(a)\subseteq_{\Gamma(\mathbf X)}\mathbf R(\alpha)\}$ and observe
\begin{align*}
\ax_0^{L[\alpha]}=\forall x\in L[\alpha]\exists y\in L[\alpha].\,x\in y\,&\simeq\,\textstyle\bigwedge_{a\in\iota}\exists y\in L[\alpha].\,a\in y,\\
\exists y\in L[\alpha].\,a\in y\,&\simeq\,\textstyle\bigvee_{b\in\iota} a\in b.
\end{align*}
Given an arbitrary $a\in\iota$, we must thus derive $a\in b$ for a suitable~$b\in\iota$. Let us set
\begin{equation*}
b:=L^u_r\quad\text{with}\quad r:=\begin{cases}
0 & \text{if }a\in u,\\
s+1 & \text{if }a=L^u_s\text{ or }a=\{x\in L^u_s\,|\,\theta(x,\mathbf c)\}.
\end{cases}
\end{equation*}
In the more interesting second case, we note that $a\in\iota$ and $s\in\supp^{\mathbf L}_{\Gamma(\mathbf X)}(a)$ entail
\begin{equation*}
s<_{\Gamma(\mathbf X)}\mathbf R(\alpha)\in\rng(\gamma_{\mathbf X})=\{\Gamma_x\,|\,x\in\mathbf X\}.
\end{equation*}
We can infer $s+1<\mathbf R(\alpha)$ by Lemma~\ref{lem:addition} (recall $1=\overline\varphi_00$). Let us rewrite this as
\begin{equation*}
\supp^{\mathbf L}_{\Gamma(\mathbf X)}(b)\subseteq_{\Gamma(\mathbf X)}\mathbf R(\alpha),
\end{equation*}
which also holds when we have $a\in u$ and hence $r=0$. As in the previous proof, we use Lemma~\ref{lem:ranks-operators} to conclude that $\omega\cdot\rk(b)+n<\psi_{\alpha+1}0$ holds for all~$n\in\mathbb N$. Given that Definition~\ref{def:Gamma-dilator} yields $\Gamma(I)(0)=0$ and $\Gamma(I)(s+1)=\Gamma(I)(s)+1$, we can employ Corollary~\ref{cor:H-nice} to get $\supp^+(b)=\{\Gamma(I)(r)\}\subseteq\mathcal H_0(\supp^+(a))$ and hence
\begin{equation*}
\rk(b)\in\mathcal H_0(\supp^+(b))\subseteq\mathcal H_0(\supp^+(a)).
\end{equation*}
Let us now recall that \cite[Definition~3.12]{freund-equivalence} yields
\begin{equation*}
a\in b\,\simeq\,\textstyle\bigvee_{c\in\kappa}c=a\quad\text{with}\quad\kappa=\{c\in\mathbf L^u_{\Gamma(\mathbf X)}\,|\,\supp^{\mathbf L}_{\Gamma(\mathbf X)}(c)\subseteq_{\Gamma(\mathbf X)} r\}.
\end{equation*}
As Lemma~\ref{lem:a=a} provides a derivation of~$a=a$, we take~$c$ to be the term~$a$. Note that the choice of~$r$ ensures $a\in\kappa$ and $\supp^+(a)\subseteq_{\mathbf O}\Gamma(I)(r)\leq\omega\cdot\rk(b)$. We may thus apply clause~(ii) of Definition~\ref{def:H-controlled-deriv}, to get
\begin{equation*}
(0,\supp^+(a))\vdash^{\omega\cdot\rk(b)}_0 a\in b.
\end{equation*}
In view of $b\in\iota$ and $\supp^+(b)\subseteq_{\mathbf O}\omega\cdot\rk(b)+1$, the same clause now yields
\begin{equation*}
(0,\supp^+(a))\vdash^{\omega\cdot\rk(b)+1}_0 \exists y\in L[\alpha].\,a\in y.
\end{equation*}
Since $a\in\iota$ was arbitrary and we always have $\omega\cdot\rk(b)+1<\omega\cdot\psi_{\alpha+1}0=\psi_{\alpha+1}0$, we can conclude by clause~(i) of Definition~\ref{def:H-controlled-deriv}.
\end{proof}

To conclude this section, we show that the search tree $S^{\mathbf R}_{\Gamma(\mathbf X)}$ from Definition~\ref{def:search} can be converted into an infinite proof. We are particularly interested in the root node~$\langle\rangle\in S^{\mathbf R}_{\Gamma(\mathbf X)}$, which gives rise to elements
\begin{equation*}
\mathbf X+\langle\rangle\in\mathbf X+S^{\mathbf R}_{\Gamma(\mathbf X)}=\mathbf K\qquad\text{and}\qquad\Gamma_{\mathbf X+\langle\rangle}\in\Gamma(\mathbf K)=\mathbf O.
\end{equation*}
The label $l_{\Gamma(\mathbf X)}(\langle\rangle)$ at the root is the empty sequent, which we denote by~$\langle\rangle$ as well.

\begin{theorem}[Embedding]\label{thm:embedding}
We have $(0,\emptyset)\vdash^t_t\langle\rangle$ for $t=\Gamma_{\mathbf X+\langle\rangle}\in\mathbf O$.
\end{theorem}
\begin{proof}
For $\sigma=\langle\sigma_0,\ldots,\sigma_{n-1}\rangle\in S^{\mathbf R}_{\Gamma(\mathbf X)}\subseteq(\mathbf L^u_{\Gamma(\mathbf X)})^{<\omega}$ we extend Definition~\ref{def:supp-extend} by
\begin{equation*}
\supp^+(\sigma):=[\Gamma(I)]^{<\omega}\circ\supp^S_{\Gamma(\mathbf X)}(\sigma)=\textstyle\bigcup_{i<n}\supp^+(\sigma_i)\in[\mathbf O]^{<\omega},
\end{equation*}
where the second equality uses~$\supp^S_{\Gamma(\mathbf X)}(\sigma)=\bigcup_{i<n}\supp^{\mathbf L}_{\Gamma(\mathbf X)}(\sigma_i)$ from Section~\ref{sect:search-trees}. Let us write $l(\sigma)=l_{\Gamma(\mathbf X)}(\sigma)$ for the sequent label from Definition~\ref{def:search}. We will show
\begin{equation*}
(0,\supp^+(\sigma))\vdash^s_s l(\sigma)\quad\text{with}\quad s:=\Gamma_{\mathbf X+\sigma}
\end{equation*}
by induction on~$\sigma\in S^{\mathbf R}_{\Gamma(\mathbf X)}$ in the Kleene-Brouwer order, which is well founded due to the embedding $\sigma\mapsto\Gamma_{\mathbf X+\sigma}$ into the well order~$\mathbf O$. Note that the theorem is the case of the root~$\sigma=\langle\rangle$. Considering Definition~\ref{def:H-controlled-deriv}, we first show
\begin{equation*}
\{\Gamma_{\mathbf X+\sigma}\}\cup\supp^+(l(\sigma))\subseteq\mathcal H_0(\supp^+(\sigma)).
\end{equation*}
In view of $\Gamma(I)\circ\mathbf R(\alpha)=\psi_{\alpha+1}0\in\mathcal H_0(\emptyset)$, the claim about $\supp^+(l(\sigma))$ reduces to Corollary~\ref{cor:search-parameters}. To conclude via Proposition~\ref{prop:supp-H}, we assume $r\in\supp^+(\Gamma_{\mathbf X+\sigma})$ and derive $r\in\mathcal H_0(\supp^+(\sigma))$. Definitions~\ref{def:supp-O} and~\ref{def:psi+} yield $r=\gamma_{\mathbf K}\circ I(x)$ for some
\begin{align*}
x\in\supp^{\mathbf O}(\Gamma_{\mathbf X+\sigma})&=\bigcup\{\supp^{\mathbf K}(\rho)\,|\,\rho\in\supp^{\Gamma}_{\mathbf K}(\Gamma_{\mathbf X+\sigma})\}\\
{}&=\supp^{\mathbf K}(\mathbf X+\sigma)=\bigcup\{\supp^\Gamma_{\mathbf X}(\rho)\,|\,\rho\in\supp^S_{\Gamma(\mathbf X)}(\sigma)\backslash\rng(\mathbf R)\}.
\end{align*}
We thus get $x\in\supp^\Gamma_{\mathbf X}(\rho)$ with $\rho\in\supp^S_{\Gamma(\mathbf X)}(\sigma)$ and hence $\Gamma(I)(\rho)\in\supp^+(\sigma)$. By Lemma~\ref{lem:supp-O} and the other direction of Proposition~\ref{prop:supp-H}, we obtain
\begin{multline*}
r\in[\gamma_{\mathbf K}\circ I]^{<\omega}\circ\supp^\Gamma_{\mathbf X}(\rho)=[\gamma_{\mathbf K}\circ I]^{<\omega}\circ\supp^{\mathbf O}(\Gamma(I)(\rho))\\
{}=\supp^+(\Gamma(I)(\rho))\subseteq\mathcal H_0(\supp^+(\sigma)).
\end{multline*}
In our induction along the Kleene-Brouwer order, we distinguish cases according to Definition~\ref{def:search}. Let us first assume that $\sigma$ has even length~$2k$, where $k$ codes a pair~$\langle n,i\rangle$. For $\alpha=\nu_i$, the cited definition provides $\sigma^\frown L[\alpha]\in S^{\mathbf R}_{\Gamma(\mathbf X)}$, and the induction hypothesis yields
\begin{equation*}
(0,\supp^+(\sigma)\cup\supp^+(L[\alpha]))\vdash^r_r l(\sigma),\neg\ax_n^{L[\alpha]}\quad\text{with}\quad r=\Gamma_{\mathbf X+\sigma^\frown L[\alpha]}<\Gamma_{\mathbf X+\sigma}.
\end{equation*}
Here we can omit $\supp^+(L[\alpha])\subseteq\mathcal H_0(\emptyset)$ by `weakening'. From Lemma~\ref{lem:psi-less-Gamma} we get
\begin{equation*}
\psi_{\alpha+1}0=\Gamma(I)\circ\mathbf R(\alpha)<_{\mathbf O}\Gamma_{\mathbf X+\sigma^\frown L[\alpha]}=r.
\end{equation*}
Due to Propositions~\ref{prop:ax1+n} (for $n>0$) and~\ref{prop:ax0} (for $n=0$), we thus have
\begin{equation*}
(0,\supp^+(\sigma))\vdash^r_r l(\sigma),\ax_n^{L[\alpha]}.
\end{equation*}
As $\psi_{\alpha+1}0<s=\Gamma_{\mathbf X+\sigma}$ entails $\rk(\ax_n^{L[\alpha]})<s$, we can complete the induction step by clause~(iii) of Definition~\ref{def:H-controlled-deriv} (`cut rule'). The other cases from Definition~\ref{def:search} correspond directly to clauses~(i) and~(ii). Concerning the disjunctive case, we note that $\supp^+(b)\subseteq\rng(\Gamma(I))$ entails $\supp^+(b)\subseteq_{\mathbf O}\Gamma_{\mathbf X+\sigma}$, again by Lemma~\ref{lem:psi-less-Gamma}.
\end{proof}

\section{An abstract ordinal analysis}\label{sect:ordinal-analysis}

In this section, we show cut elimination and collapsing results that entail the consistency of our infinite proof system. On the one hand, these results resemble the known ordinal analysis of iterated admissibility~\cite{buchholz-local-predicativity,jaeger-iterated-admissible,pohlers-iterated-admissibles,Rathjen_PhD_2013}. On the other hand, our setting here is more abstract, since we work relative to the given dilator $\Gamma\circ\mathbf S$ from Assumption~\ref{ass:collapse-search} (recall that $\mathbf S$ arises from the search trees of Definition~\ref{def:search}). Once consistency is available, it will be straightforward to deduce the main result of our paper, as we shall see in the next section. We begin with a standard ingredient for cut elimination (cf.~\cite[Lemma~3.13]{buchholz-local-predicativity}):

\begin{lemma}[Inversion]\label{lem:inversion}
If $\varphi\simeq\bigwedge_{b\in\iota(\varphi)}\varphi_b$ is conjunctive, then we have
\begin{equation*}
(r,a)\vdash_s^t\Gamma,\varphi\qquad\Rightarrow\qquad(r,a\cup\supp^+(b))\vdash_s^t\Gamma,\varphi_b\quad\text{for any}\quad b\in\iota(\varphi).
\end{equation*}
\end{lemma}
\begin{proof}
Due to the initial condition from Definition~\ref{def:H-controlled-deriv}, the premise of the desired implication entails $\supp^+(\varphi)\subseteq\mathcal H_r(a)$. As in the proof of Lemma~\ref{lem:a=a} we get
\begin{equation*}
\supp^+(\varphi_b)\subseteq\supp^+(\varphi)\cup\supp^+(b)\subseteq\mathcal H_r(a\cup\supp^+(b)),
\end{equation*}
which ensures that the same initial condition holds for the conclusion. We now argue by induction on~$t\in\mathbf O$. In the crucial case, clause~(i) of Definition~\ref{def:H-controlled-deriv} was applied to the distinguished formula~$\varphi$, so that we have
\begin{equation*}
(r,a\cup\supp^+(b))\vdash_s^{t(b)}\Gamma,\varphi,\varphi_b
\end{equation*}
for some~$t(b)<t$. Here we can omit $\varphi$ due to the induction hypothesis. Weakening (Lemma~\ref{lem:weakening}) allows us to increase~$t(b)$ to~$t$, which yields the desired conclusion. In all other cases, one uses the induction hypothesis and reapplies the same clause. The latter is possible because clauses~(ii) and~(iv) concern formulas that are disjunctive and hence different from~$\varphi$.
\end{proof}

The following result (cf.~\cite[Lemma~3.14]{buchholz-local-predicativity}) shows how certain applications of the cut rule can be avoided. Let us point out that we cannot conclude by clause~(iii) of Definition~\ref{def:H-controlled-deriv}, since the latter would require $\rk(\psi)=\rk(\neg\psi)<\rk(\psi)$.

\begin{lemma}[Reduction]\label{lem:reduction}
For disjunctive~$\psi$ with $\rk(\psi)\notin\{\psi_{\alpha+1}0\,|\,\alpha<\nu\}$ we have
\begin{equation*}
(r,a)\vdash_{\rk(\psi)}^{t(0)}\Gamma,\neg\psi\quad\text{and}\quad(r,a)\vdash_{\rk(\psi)}^{t(1)}\Gamma,\psi\qquad\Rightarrow\qquad(r,a)\vdash_{\rk(\psi)}^{t(0)+t(1)}\Gamma.
\end{equation*}
\end{lemma}
\begin{proof}
The premise of the desired implication entails $t(i)\in\mathcal H_r(a)$ for~$i\in\{0,1\}$, as in the previous proof. Thus $t(0)+t(1)\in\mathcal H_r(a)$ holds by Corollary~\ref{cor:H-nice} in conjunction with Lemma~\ref{lem:Veblen-addition-support}. We now argue by induction on~$t(1)$ and distinguish cases according to the clause of Definition~\ref{def:H-controlled-deriv} that was used to derive $\Gamma,\psi$. In the crucial case, the formula $\psi\simeq\bigvee_{b\in\iota(\psi)}\psi_b$ itself was derived by clause~(ii), which means that we have
\begin{equation*}
(r,a)\vdash_{\rk(\psi)}^{s}\Gamma,\psi,\psi_b\quad\text{for some }b\in\iota(\psi)\text{ and }s<_{\mathbf O}t(1).
\end{equation*}
In particular, this means that we have $\supp^+(\psi_b)\subseteq\mathcal H_r(a)$, by the initial condition from Definition~\ref{def:H-controlled-deriv}. We may also assume $\supp^+(b)\subseteq\supp^+(\psi_b)$. Indeed, this is immediate if~$b$ occurs in~$\psi_b$. If it does not, then we have $\psi_b=\psi_i$ for some index $i\in\{0,1\}\subseteq\iota(\psi)$, as a glance at~\cite[Definition~3.12]{freund-equivalence} reveals. In this case we may thus redefine $b:=i\in u\subseteq\mathbf L^u_{\Gamma(\mathbf X)}$ to get $\supp^+(b)=\emptyset$. Let us now apply weakening to the given derivation of $\Gamma,\neg\psi$, so that we obtain
\begin{equation*}
(r,a)\vdash_{\rk(\psi)}^{t(0)}\Gamma,\neg\psi,\psi_b.
\end{equation*}
By the induction hypothesis, we can then infer
\begin{equation*}
(r,a)\vdash_{\rk(\psi)}^{t(0)+s}\Gamma,\psi_b.
\end{equation*}
From~\cite[Definition~3.12]{freund-equivalence} we know that $\neg\psi$ is conjunctive with $(\neg\psi)_b=\neg(\psi_b)$ for all $b\in\iota(\neg\psi)=\iota(\psi)$. We may thus apply inversion (Lemma~\ref{lem:inversion}) to the given derivation of $\Gamma,\neg\psi$, in order to get
\begin{equation*}
(r,a\cup\supp^+(b))\vdash^{t(0)}_{\rk(\psi)}\Gamma,\neg\psi_b.
\end{equation*}
For $b$ as above, we may omit $\supp^+(b)\subseteq\mathcal H_r(a)$ by weakening. As Lemma~\ref{lem:rank-disj-conj} ensures $\rk(\psi_b)<\rk(\psi)$, we can conclude by clause~(iii) of Definition~\ref{def:H-controlled-deriv}. In all other cases, one uses the induction hypothesis and reapplies the same clause. Here it is crucial to observe that clause~(iv) cannot be applied with $\psi=(\exists z\in L[\alpha].\varphi^z)$. Indeed, given that~$\varphi$ is a $\Sigma(\alpha)$-formula, we have
\begin{equation*}
\supp^+(\varphi)\subseteq_{\mathbf O}\Gamma(I)\circ\mathbf R(\alpha)=\psi_{\alpha+1}0.
\end{equation*}
We may replace~$\varphi$ by the `trivial' relativization~$\varphi^0$, since we have~$\supp^+(0)=\emptyset$. As $\psi_{\alpha+1}0$ is strongly critical (cf.~the proof of Proposition~\ref{prop:C-psi}), Lemma~\ref{lem:ranks-operators} yields
\begin{equation*}
\rk(\varphi^0)+2<_{\mathbf O}\omega\cdot(1+\psi_{\alpha+1}0)=\psi_{\alpha+1}0.
\end{equation*}
Similarly, we get $\rk(L[\alpha])=\rk(L^u_{\mathbf R(\alpha)})=\psi_{\alpha+1}0$ and then
\begin{equation*}
\rk(\exists z\in L[\alpha].\,\varphi^z)=\max\{\rk(L[\alpha]),\rk(\varphi^0)+2\}=\psi_{\alpha+1}0\neq\rk(\psi).
\end{equation*}
The inequality holds by an assumption in the lemma, which thus excludes an obstructive application of clause~(iv) from Definition~\ref{def:H-controlled-deriv}.
\end{proof}

By the next result (cf.~\cite[Theorem~3.16]{buchholz-local-predicativity}), the cut rank can be reduced when no critical value~$\psi_{\alpha+1}0$ is involved. To remove this last restriction, we will later prove a collapsing result that complements cut elimination. Let us point out that $\varphi$ refers to the Veblen function from Definition~\ref{def:total-Veblen}.

\begin{proposition}[Predicative cut elimination]\label{prop:pred-cut-elim}
Consider elements $p,q\in\mathbf O$ such that $p\leq\psi_{\alpha+1}0<p+\varphi(0,q)$ fails for all~$\alpha<\nu$. We then have
\begin{equation*}
(r,a)\vdash^t_{p+\varphi(0,q)}\Gamma\quad\text{and}\quad q\in\mathcal H_r(a)\qquad\Rightarrow\qquad (r,a)\vdash^{\varphi(q,t)}_p\Gamma.
\end{equation*}
\end{proposition}
\begin{proof}
The assumption of the desired implication entails $q,t\in\mathcal H_r(a)$, due to the initial condition from Definition~\ref{def:H-controlled-deriv}. We get $\varphi(q,t)\in\mathcal H_r$ by Corollary~\ref{cor:H-nice} in conjunction with Lemma~\ref{lem:Veblen-addition-support}. Let us now argue by main induction on~$q$ and side induction on~$t$ (where $p$ may vary during the induction). In the crucial case, we are concerned with clause~(iii) of Definition~\ref{def:total-Veblen}, so that we have
\begin{equation*}
(r,a)\vdash^s_{p+\varphi(0,q)}\Gamma,\psi\qquad\text{and}\qquad (r,a)\vdash^s_{p+\varphi(0,q)}\Gamma,\neg\psi
\end{equation*}
with $\rk(\psi)<p+\varphi(0,q)$ and $s<t$. For later use we record $\supp^+(\psi)\subseteq\mathcal H_r(a)$. The side induction hypothesis yields
\begin{equation*}
(r,a)\vdash^{\varphi(q,s)}_p\Gamma,\psi\qquad\text{and}\qquad (r,a)\vdash^{\varphi(q,s)}_p\Gamma,\neg\psi.
\end{equation*}
If we have $\rk(\psi)<p$, then we can conclude by clause~(iii) of Definition~\ref{def:H-controlled-deriv}, since Proposition~\ref{prop:total-Veblen} yields~$\varphi(q,s)<\varphi(q,t)$. By the same proposition and Lemma~\ref{lem:addition}, we even have $\varphi(q,s)+\varphi(q,s)<\varphi(q,t)$. Now assume $p\leq\rk(\psi)$ and note that this entails $\rk(\neg\psi)=\rk(\psi)\notin\{\psi_{\alpha+1}0\,|\,\alpha<\nu\}$. Let us recall that $\neg\neg\psi$ is syntactically equal to $\psi$, since we treat negation as a defined operation on formulas in negation normal form. Either $\psi=\neg\neg\psi$ or $\neg\psi$ is disjunctive, as seen in~\cite[Definition~3.12]{freund-equivalence}. We can thus use reduction (Lemma~\ref{lem:reduction}) and weakening to get
\begin{equation*}
(r,a)\vdash^{\varphi(q,t)}_{\rk(\psi)}\Gamma.
\end{equation*}
Lemma~\ref{lem:addition} yields $\rk(\psi)=p+s$ for some $s\in\mathbf O=\Gamma(\mathbf K)$. By Definition~\ref{def:Gamma(X)} we may write $s=\langle s_0,\ldots,s_{n-1}\rangle$ with $s_i\in\mathsf H$ (not necessarily with $n>1$). We thus get
\begin{equation*}
p=p+s(0)\quad\text{and}\quad\rk(\psi)=p+s(n)\quad\text{for}\quad s(i):=\langle s_0,\ldots,s_{i-1}\rangle.
\end{equation*}
By an auxiliary induction from $i=n$ down to~$i=0$, we now show
\begin{equation*}
(r,a)\vdash^{\varphi(q,t)}_{p+s(i)}\Gamma.
\end{equation*}
In the induction step, we use Proposition~\ref{prop:total-Veblen} to write $s_i\in\mathsf H$ in the form $\varphi(p_i,q_i)$. Let us set $q(i):=q_i$ when $p_i=0$ and $q(i):=s_i$ when $0<p_i$. In the second case, Proposition~\ref{prop:total-Veblen} yields $\varphi(0,s_i)=s_i$. So we always get
\begin{equation*}
s(i+1)=s(i)+s_i=s(i)+\varphi(0,q(i)).
\end{equation*}
Let us observe that we have
\begin{equation*}
\supp^{\mathbf L}_{\Gamma(\mathbf X)}(q(i))=\supp^{\mathbf L}_{\Gamma(\mathbf X)}(s_i)\subseteq\supp^{\mathbf L}_{\Gamma(\mathbf X)}(\rk(\psi)).
\end{equation*}
As Lemma~\ref{lem:ranks-operators} provides $\rk(\psi)\in\mathcal H_0(\supp^+(\psi))\subseteq\mathcal H_r(a)$, we obtain $q(i)\in\mathcal H_r(a)$ by Corollary~\ref{cor:H-nice}. Furthermore, we have $q(i)<q$ due to
\begin{equation*}
p+s(i)+\varphi(0,q(i))=p+s(i+1)\leq\rk(\psi)<p+\varphi(0,q)\leq p+s(i)+\varphi(0,q).
\end{equation*}
Given the auxiliary induction hypothesis (with $i+1$ at the place of~$i$), we use the main induction hypothesis (with $p+s(i)$ and $q(i)$ at the place of~$p$ and $q$) to get
\begin{equation*}
(r,a)\vdash^{\varphi(q(i),\varphi(q,t))}_{p+s(i)}\Gamma.
\end{equation*}
From Proposition~\ref{prop:total-Veblen} we know that $q(i)<q$ entails $\varphi(q(i),\varphi(q,t))=\varphi(q,t)$. So the step of the auxiliary induction is completed. Taking $i=0$ completes the present case of the side and main induction step. The remaining cases are straightforward.
\end{proof}

So far, the notation $\varphi^a$ for relativization has been introduced for~$a\in\mathbf L^u_{\Gamma(\mathbf X)}$ only. We now use the embedding $\Gamma(I):\Gamma(\mathbf X)\to\Gamma(\mathbf K)=\mathbf O$ to overload the notation.

\begin{definition}
Given an $\mathbf L^u_{\Gamma(\mathbf X)}$-formula $\varphi$ and an element $t\in\rng(\Gamma(I))\subseteq\mathbf O$, we set $\varphi^t:=\varphi^a$ with $a:=L^u_s\in\mathbf L^u_{\Gamma(\mathbf X)}$ for the unique $s\in\Gamma(\mathbf X)$ with $\Gamma(I)(s)=t$.
\end{definition}

It is instructive to recall $\psi_{\alpha+1}0=\Gamma(I)\circ\mathbf R(\alpha)$ and $L[\alpha]=L^u_{\mathbf R(\alpha)}$, which yields
\begin{equation*}
\varphi^t=\varphi^{L[\alpha]}\quad\text{for}\quad t=\psi_{\alpha+1}0\in\mathbf O.
\end{equation*}
By Corollary~\ref{cor:normal-initial} and Definition~\ref{def:pi_X}, the range of~$\Gamma(I)$ is an initial segment of~$\mathbf O$. It follows that $\varphi^t$ is defined whenever $t\leq\psi_{\alpha+1}0$ holds for some~$\alpha<\nu$. We~will later need the following variant of inversion (cf.~Lemma~\ref{lem:inversion}).

\begin{lemma}\label{lem:inversion-variant}
Given $q\leq\psi_{\alpha+1}0$ and $q\in\mathcal H_r(a)$, we get
\begin{equation*}
(r,a)\vdash^t_s\Gamma,\forall x\in L[\alpha].\,\theta\qquad\Rightarrow\qquad (r,a)\vdash^t_s\Gamma,(\forall x.\theta)^q,
\end{equation*}
for any bounded $\mathbf L^u_{\Gamma(\mathbf X)}$-formula~$\theta=\theta(x)$.
\end{lemma}
\begin{proof}
Write $\varphi:=\forall x\in L[\alpha].\,\theta$ and $\psi:=(\forall x.\theta)^q=\forall x\in L^u_p.\,\theta$ with $\Gamma(I)(p)=q$, and note that $q\leq\psi_{\alpha+1}0$ entails $p\leq\mathbf R(\alpha)$. The initial condition of Definition~\ref{def:H-controlled-deriv} is preserved as we have $\supp^+(\psi)\subseteq\supp^+(\varphi)\cup\{q\}$. In view of~\cite[Definition~3.12]{freund-equivalence}, the formulas $\varphi$ and $\psi$ are both conjunctive, and we have $\varphi_b=\theta(b)=\psi_b$ for any
\begin{equation*}
b\in\iota(\psi)=\{a\in\mathbf L^u_{\Gamma(\mathbf X)}\,|\,\supp^{\mathbf L}_{\Gamma(\mathbf X)}(a)\subseteq_{\Gamma(\mathbf X)}p\}\subseteq\iota(\varphi).
\end{equation*}
So whenever clause~(i) of Definition~\ref{def:H-controlled-deriv} is used to derive $\varphi$, it can also derive $\psi$. Based on this observation, the claim is readily established by induction on~$t$.
\end{proof}

Let us also record how relativization interacts with our assignment of a disjunction $\varphi\simeq\bigvee_{b\in\iota(\varphi)}\varphi_b$ or conjunction $\varphi\simeq\bigwedge_{b\in\iota(\varphi)}\varphi_b$ to each formula~$\varphi$.

\begin{lemma}\label{lem:Sigma-formulas}
The following holds for any $\mathbf L^u_{\Gamma(\mathbf X)}$-formula~$\varphi$ and any $t\in\rng(\Gamma(I))$:
\begin{enumerate}[label=(\alph*),itemsep=3pt]
\item The formula $\varphi^t$ is conjunctive or disjunctive, respectively, if and only if the same holds for~$\varphi$.
\item We have $(\varphi^t)_b=(\varphi_b)^t$ for any $b\in\iota(\varphi^t)\subseteq\iota(\varphi)$.
\item For any $b\in\iota(\varphi)$ with $\supp^+(b)\subseteq_{\mathbf O}t$, we have $b\in\iota(\varphi^t)$.
\item If $\varphi$ is a $\Sigma(\alpha)$-formula, then so is $\varphi_b$ for any $b$ in the set
\begin{equation*}
\iota(\varphi^{L[\alpha]})=\{b\in\iota(\varphi)\,|\,\supp^+(b)\subseteq_{\mathbf O}\psi_{\alpha+1}0\}.
\end{equation*}
\item Assume $\varphi$ is a conjunctive $\Sigma(\alpha)$-formula. We then have~$\iota(\varphi^t)=\iota(\varphi)$. Also, there is an $s\in\supp^+(\varphi)\cup\{0\}$ with $\supp^+(b)\subseteq_{\mathbf O}s$ for all~$b\in\iota(\varphi)$.
\end{enumerate}
\end{lemma}
In part~(e), we get~$s<\psi_{\alpha+1}0$ due to the definition of~$\Sigma(\alpha)$-formulas. So when~$\varphi$ is conjunctive, part~(d) applies to any element $b\in\iota(\varphi)$.
\begin{proof}
All claims can be verified explicitly, based on~\cite[Definition~3.12]{freund-equivalence}. Details for a representative case are given in the proof of~\cite[Lemma~9.1]{freund-equivalence}. Concerning~(d), we note that $\supp^+(b)\subseteq_{\mathbf O}\psi_{\alpha+1}0$ is equivalent to $\supp^{\mathbf L}_{\Gamma(\mathbf X)}(b)\subseteq_{\Gamma(\mathbf X)}\mathbf R(\alpha)$, which relates to~Definition~\ref{def:supp-extend}. In part~(e), the crucial point is that $\varphi$ cannot begin with an unbounded quantifier.
\end{proof}

Clause~(iv) from Definition~\ref{def:H-controlled-deriv} is an obstruction to cut elimination, as we have seen in the proof of Lemma~\ref{lem:reduction}. The following result (cf.~\cite[Lemma~3.17]{buchholz-local-predicativity}) will allow us to circumvent this clause, since $\varphi^t$ with $t<\psi_{\alpha+1}0$ entails $\exists z\in L[\alpha].\,\varphi^z$.

\begin{proposition}[Boundedness]\label{prop:boundedness}
For each $\Sigma(\alpha)$-formula $\varphi$ with $\alpha<\nu$ we have
\begin{equation*}
(r,a)\vdash_q^s\Gamma,\varphi^{L[\alpha]}\text{ and }s\leq t<\psi_{\alpha+1}0\text{ with }t\in\mathcal H_r(a)\quad\Rightarrow\quad(r,a)\vdash_q^s\Gamma,\varphi^t.
\end{equation*}
\end{proposition}
\begin{proof}
First note that the antecedent of the desired implication entails
\begin{equation*}
\supp^+(\varphi^t)\subseteq\supp^+(\varphi^{L[\alpha]})\cup\{t\}\subseteq\mathcal H_r(a),
\end{equation*}
so that the initial condition from Definition~\ref{def:H-controlled-deriv} is preserved. We now argue by induction on~$s$. When the relevant clause from Definition~\ref{def:H-controlled-deriv} does not refer to~$\varphi^{L[\alpha]}$, it is straightforward to reduce to the induction hypothesis. In case clause~(i) applies to~$\varphi^{L[\alpha]}$, the latter is conjunctive and we have
\begin{equation*}
(r,a\cup\supp^+(b))\vdash^{s(b)}_q\Gamma,\varphi^{L[\alpha]},(\varphi^{L[\alpha]})_b\quad\text{with }s(b)<s\text{ for all }b\in\iota(\varphi^{L[\alpha]}).
\end{equation*}
The previous lemma ensures that $\varphi_b$ is a $\Sigma(\alpha)$-formula with $(\varphi_b)^{L[\alpha]}=(\varphi^{L[\alpha]})_b$, for any $b\in\iota(\varphi)=\iota(\varphi^{L[\alpha]})$. Thus two applications of the induction hypothesis yield
\begin{equation*}
(r,a\cup\supp^+(b))\vdash^{s(b)}_q\Gamma,\varphi^t,(\varphi_b)^t\quad\text{for all }b\in\iota(\varphi).
\end{equation*}
Using the previous lemma once again, we learn that $\varphi$ and hence~$\varphi^t$ is conjunctive with $(\varphi^t)_b=(\varphi_b)^t$ for all $b\in\iota(\varphi^t)\subseteq\iota(\varphi)$. In order to conclude the present case of the induction step, we can thus reapply clause~(i). A similar argument covers clause~(ii), as the previous lemma ensures the following: for any $b\in\iota(\varphi^{L[\alpha]})\subseteq\iota(\varphi)$ with $\supp^+(b)\subseteq_{\mathbf O}s\leq t<\psi_{\alpha+1}$, we have $b\in\iota(\varphi^t)$ and $\varphi_b$ is a $\Sigma(\alpha)$-formula. Finally, we consider an application of clause~(iv) for a $\Sigma(\beta)$-formula $\theta$ with $\beta<\nu$ and $(\exists z\in L[\beta].\,\theta^z)=\varphi^{L[\alpha]}$, where we have
\begin{equation*}
(r,a)\vdash^{s(0)}_q\Gamma,\varphi^{L[\alpha]},\theta^{L[\beta]}\quad\text{for some }s(0)<s.
\end{equation*}
If we have~$\alpha\neq\beta$, then $L[\beta]$ occurs in~$\varphi$, and the definition of~$\Sigma(\alpha)$-formulas yields
\begin{equation*}
\mathbf R(\beta)\in\supp^{\mathbf L}_{\Gamma(\mathbf X)}(L[\beta])\subseteq\supp^{\mathbf L}_{\Gamma(\mathbf X)}(\varphi)\subseteq_{\Gamma(\mathbf X)}\mathbf R(\alpha).
\end{equation*}
So in any case we have $\beta\leq\alpha$. By a similar argument, it follows that $L[\alpha]$ cannot occur in the $\Sigma(\beta)$-formula~$\theta$. In case $\beta<\alpha$ we thus get $\varphi^{L[\alpha]}=\varphi=\varphi^t$, which makes the claim trivial. Now assume $\beta=\alpha$ and note that this forces~$\varphi=\exists z.\,\theta^z$. We apply the induction hypothesis twice (once with $s(0)$ at the place of~$t$), to get
\begin{equation*}
(r,a)\vdash^{s(0)}_q\Gamma,\varphi^t,\theta^{s(0)}.
\end{equation*}
For $p\in\Gamma(\mathbf X)$ with $\Gamma(I)(p)=t$ we have
\begin{equation*}
\varphi^t=\exists z\in L^u_p.\,\theta^z\simeq\textstyle\bigvee_{b\in\iota(\varphi^t)}\theta^b\quad\text{with}\quad\iota(\varphi^t)=\{b\in\mathbf L^u_{\Gamma(\mathbf X)}\,|\,\supp^{\mathbf L}_{\Gamma(\mathbf X)}(b)\subseteq_{\Gamma(\mathbf X)}p\}.
\end{equation*}
Now $\theta^{s(0)}=\theta^b$ holds for $b:=L^u_{p(0)}$ with $\Gamma(I)(p(0))=s(0)<s\leq t$, which yields
\begin{equation*}
\supp^{\mathbf L}_{\Gamma(\mathbf X)}(b)=\{p(0)\}\subseteq_{\Gamma(\mathbf X)}p\quad\text{and}\quad\supp^+(b)=\{s(0)\}\subseteq_{\mathbf O}s.
\end{equation*}
We can thus conclude by an application of clause~(ii) from Definition~\ref{def:H-controlled-deriv}.
\end{proof}

The following definition adapts notation from~\cite[Section~4]{buchholz-local-predicativity}, which will be used for the crucial result on collapsing and impredicative cut elimination. The reader may wish to recall Definitions~\ref{def:psi-C} and~\ref{def:opertors-H} as well as the paragraph before Theorem~\ref{thm:embedding}.

\begin{definition}\label{def:overline-K}
For $\alpha<\nu$ and $r,s\in\mathbf O$ and $a\in[\mathbf O]^{<\omega}$, we abbreviate
\begin{equation*}
\mathcal A(a;r,\alpha,s)\quad :\Leftrightarrow\quad r,s\in\mathcal H_r(a)\text{ and }a\subseteq\textstyle\bigcap_{\beta\geq\alpha} C_\beta(r+1).
\end{equation*}
Let us also put $\overline K:=\{\Omega(\alpha)\,|\,\alpha\leq\nu\}$ with
\begin{alignat*}{3}
\Omega(0)&:=0,\qquad &\Omega(\alpha+1)&:=(\psi_{\alpha+1}0)+1,\\
\Omega(\nu)&:=\Gamma_{\mathbf X+\langle\rangle}\qquad &\Omega(\lambda)&:=\psi_\lambda 0\quad\text{for each limit $\lambda<\nu$}.
\end{alignat*}
\end{definition}

Note that we have $\Omega(\alpha)\in\mathcal H_0(\emptyset)$ for all $\alpha\leq\nu$, as a consequence of Theorem~\ref{thm:embedding} and Corollary~\ref{cor:operator-collapsing}. For $\alpha<\nu$, the following result characterizes $\Omega(\alpha)\in\rng(\Gamma(I))$ as a supremum.

\begin{lemma}\label{lem:Omega-alpha-sup}
For any $\alpha\leq\nu$ and $s\in\rng(\Gamma(I))\subseteq\mathbf O$ we have
\begin{equation*}
s<_{\mathbf O}\Omega(\alpha)\quad\Leftrightarrow\quad s\leq_{\mathbf O}\psi_{\beta+1}0\text{ for some }\beta<\alpha.
\end{equation*}
\end{lemma}
\begin{proof}
The claim is immediate when $\alpha$ is zero or a successor. Let us now assume that $\alpha<\nu$ is a limit. The non-trivial task is to show that $s<\psi_\alpha0$ entails $s\leq\psi_{\beta+1}0$ for some $\beta<\alpha$. Invoking Definitions~\ref{def:psi+} and~\ref{def:psi-C} as well as Lemma~\ref{lem:Gamma-normal}, we see that any $\delta<\nu$ validates
\begin{equation*}
s<_{\mathbf O}\psi_\delta 0=\gamma_{\mathbf K}\circ I(\psi^{\mathbf X}_\delta 0)\quad\Leftrightarrow\quad\supp^\Gamma_{\mathbf K}(s)\subseteq_{\mathbf K}I(\psi^{\mathbf X}_\delta 0).
\end{equation*}
Assume that these equivalent statements hold for~$\delta=\alpha$. We need to find a $\beta<\alpha$ such that they hold for~$\delta=\beta+1$ as well. Let us recall that the range of $I:\mathbf X\to\mathbf K$ is an initial segment. The maximal element of the finite set~$\supp^\Gamma_{\mathbf K}(s)$ can thus be written as $I(x)$, except in the trivial case where the support is empty. Due to Definition~\ref{def:pi_X} we get $x=\psi^{\mathbf X}_\beta t$ for some~$\beta<\nu$ and $t\in\mathbf O$. Clearly, the right side above holds for~$\delta=\beta+1$. Also, the right side for~$\delta=\alpha$ entails $I(x)<I(\psi^{\mathbf X}_\alpha 0)$ and hence $\beta<\alpha$, as desired. Now consider the case of~$\alpha=\nu$. For any~$s\in\rng(\Gamma(I))$, Definition~\ref{def:dilator} and Lemma~\ref{lem:psi-less-Gamma} yield $s<\Omega(\nu)$ and $\supp^\Gamma_{\mathbf K}(s)\subseteq\rng(I)$. Due to the latter, we can find a $\beta<\nu$ with $s<\psi_{\beta+1}0$, as in the limit case.
\end{proof}

The following transfers~\cite[Lemma~4.7]{buchholz-local-predicativity} into our setting.

\begin{lemma}\label{lem:mathcal-A}
If we have $\mathcal A(a;r,\alpha,s)$ and $t\in\mathcal H_r(a)$, then the following holds:
\begin{enumerate}[label=(\alph*)]
\item Given $\alpha\leq\beta$ and $t<\psi_{\beta+1}0$, we get $t<\psi_\beta(r+1)$.
\item For $p:=r+\varphi_0(s+t)$ we have $p\in\mathcal H_r(a)$ as well as $\psi_\alpha p\in\mathcal H_p(a)$, and $p<p'$ entails $\psi_\alpha p<\psi_\alpha p'$.
\end{enumerate}
\end{lemma}

\begin{proof}
(a) In view of Definition~\ref{def:opertors-H}, the assumptions entail $t\in C_\beta(r+1)$. Now the conclusion follows by Proposition~\ref{prop:C-psi}.

(b) In view of $r\leq p$, we can use Corollary~\ref{cor:H-nice} to get $p\in\mathcal H_r(a)\subseteq\mathcal H_{p}(a)$, which entails $\psi_\alpha p\in\mathcal H_{p}(a)$ by Corollary~\ref{cor:operator-collapsing}. Given $p<p'$, we now obtain $\psi_\alpha p\in C_\alpha(p')$, as $\mathcal A(a;r,\alpha,s)$ provides $a\subseteq C_\alpha(r+1)\subseteq C_\alpha(p')$. In order to conclude $\psi_\alpha p<\psi_\alpha p'$, it suffices to invoke Proposition~\ref{prop:C-psi} once again.
\end{proof}

Our abstract ordinal analysis culminates in the following (cf.~\cite[Theorem~4.8]{buchholz-local-predicativity}).

\begin{theorem}[Collapsing and impredicative cut elimination] For $\alpha<\nu$, assume
\begin{equation*}
(r,a)\vdash^t_s\Gamma\quad\text{with}\quad\mathcal A(a;r,\alpha,s)\quad\text{and}\quad s\in\overline K,
\end{equation*}
where all elements of $\Gamma$ have the form~$\varphi^{L[\alpha]}$ for a $\Sigma(\alpha)$-formula~$\varphi$. We then get
\begin{equation*}
(p,a)\vdash^q_q\Gamma\quad\text{with}\quad p=r+\varphi_0(s+t)\quad\text{and}\quad  q=\psi_\alpha p.
\end{equation*}
\end{theorem}
\begin{proof}
We argue by main induction on~$s$ and side induction on~$t$ (where $\alpha$ and the other parameters may vary in the induction). The previous lemma secures the initial condition from Definition~\ref{def:H-controlled-deriv}. In clause~(i) of the latter, we are concerned with a conjunctive formula~$\varphi^{L[\alpha]}\in\Gamma$ such that we have
\begin{equation*}
(r,a\cup\supp^+(b))\vdash^{t(b)}_s\Gamma,\varphi^{L[\alpha]}_b\quad\text{with }t(b)<t\text{ for all }b\in\iota(\varphi^{L[\alpha]}).
\end{equation*}
Here we write $\varphi^{L[\alpha]}_b$ for $(\varphi^{L[\alpha]})_b$, which coincides with $(\varphi_b)^{L[\alpha]}$ due to Lemma~\ref{lem:Sigma-formulas}. The latter also yields a $t'\in\supp^+(\varphi)\cup\{0\}\subseteq\mathcal H_r(a)$ with $\supp^+(b)\subseteq_{\mathbf O}t'<\psi_{\alpha+1}0$ for all $b\in\iota(\varphi^{L[\alpha]})$. To establish $\mathcal A(a\cup\supp^+(b);r,\alpha,s)$ for any such~$b$, we consider an arbitrary $\beta\geq\alpha$. By the previous lemma we get $t'<\psi_\beta(r+1)$. Let us also note that $r+1\in\mathcal H_r(a)\subseteq C_\beta(r+1)$ holds due to $\mathcal A(a;r,\alpha,s)$ and Definition~\ref{def:opertors-H}. Thus Proposition~\ref{prop:C-psi} yields $\supp^+(b)\subseteq C_\beta(r+1)$, as required. We may now use the side induction hypothesis to infer
\begin{equation*}
(p(b),a\cup\supp^+(b))\vdash^{q(b)}_{q(b)}\Gamma,\varphi^{L[\alpha]}_b\quad\text{with}\quad p(b)=r+\varphi_0(s+t(b))\quad\text{and}\quad\psi_\alpha p(b),
\end{equation*}
for any $b\in\iota(\varphi^{L[\alpha]})$. With $p$ and $q$ as in the theorem, we see that $t(b)<t$ entails $p(b)<p$ and then $q(b)<q$, by Lemma~\ref{lem:mathcal-A} with $a\cup\supp^+(b)$ at the place of~$a$. To conclude the present case of the induction step, we use weakening and reapply clause~(i) of Definition~\ref{def:H-controlled-deriv}. Now consider clause~(ii) for a disjunctive $\varphi^{L[\alpha]}\in\Gamma$~with
\begin{equation*}
(r,a)\vdash^{t(0)}_s\Gamma,\varphi^{L[\alpha]}_b\quad\text{for some}\quad t(0)<t\quad\text{and}\quad b\in\iota(\varphi^{L[\alpha]}).
\end{equation*}
As in the proof of Lemma~\ref{lem:reduction}, we may assume $\supp^+(b)\subseteq\supp^+(\varphi^{L[\alpha]}_b)$. The latter entails $\supp^+(b)\subseteq\mathcal H_r(a)$, by the initial condition from Definition~\ref{def:H-controlled-deriv}. Since we also have $\supp^+(b)\subseteq_{\mathbf O}\psi_{\alpha+1}0$ due to Lemma~\ref{lem:Sigma-formulas}, we can use Lemma~\ref{lem:mathcal-A} to get
\begin{equation*}
\supp^+(b)\subseteq_{\mathbf O}\psi_\alpha(r+1)\leq_{\mathbf O}\psi_\alpha p(0)=:q(0)\quad\text{with}\quad p(0):=r+\varphi_0(s+t(0)).
\end{equation*}
Let us recall that our version of~$\psi$ is not even weakly increasing. To secure the weak inequality above, one invokes Lemma~\ref{lem:mathcal-A} with $s=0=t$. The given bound on $\supp^+(b)$ allows us to reapply clause~(ii) after the side induction hypothesis has been used. Before we come to the crucial clause~(iii), let us consider an application of~(iv), where $\Gamma$ contains $\exists z\in L[\beta].\,\varphi^z$ for some $\Sigma(\beta)$-formula~$\varphi$. As in the proof of Proposition~\ref{prop:boundedness}, we necessarily have $\beta\leq\alpha$. To conclude by the side induction hypothesis, we need only observe that $\varphi^{L[\beta]}=\psi^{L[\alpha]}$ holds for some $\Sigma(\alpha)$-formula~$\psi$. We can take $\psi:=\varphi$ for $\beta=\alpha$ and $\psi:=\varphi^{L[\beta]}=(\varphi^{L[\beta]})^{L[\alpha]}$ for~$\beta<\alpha$. As~preparation for clause~(iv), we establish the following claim (which is adapted from the proof by Buchholz~\cite{buchholz-local-predicativity}). The quantities that appear in the theorem should be considered as fixed (for the induction step), while $p(0),q(0)$ and $\varphi$ can be arbitrary.
\begin{claim}
Assume that we have $r\leq p(0)<p$ and $p(0)\in\mathcal H_{p(0)}(a)$, and that there exists a~$\beta<\nu$ with $s(0):=\max\{q(0),\rk(\varphi)\}<\psi_{\beta+1}0\leq s$. We then get
\begin{equation*}
(p(0),a)\vdash^{q(0)}_{q(0)}\Gamma,\varphi\quad\text{and}\quad (p(0),a)\vdash^{q(0)}_{q(0)}\Gamma,\neg\varphi\qquad\Rightarrow\qquad (p,a)\vdash^q_q\Gamma.
\end{equation*}
\end{claim}
\noindent To establish the claim, we first note that clause~(iii) of Definition~\ref{def:H-controlled-deriv} yields
\begin{equation*}
(p(0),a)\vdash^{q(0)+1}_{s(0)+1}\Gamma.
\end{equation*}
For any $\beta$ as in the claim, we have $s(1):=\Omega(\beta)+\varphi_0(s(0)+1)<\psi_{\beta+1}0$, since the bound is strongly critical (cf.~the proof of Proposition~\ref{prop:C-psi}). So there is no $\gamma<\nu$ with $\Omega(\beta)\leq\psi_{\gamma+1}0<s(1)$. We can thus use predicative cut elimination to get
\begin{equation*}
(p(0),a)\vdash^{t(0)}_{\Omega(\beta)}\Gamma\quad\text{with}\quad t(0):=\varphi(s(0)+1,q(0)+1).
\end{equation*}
It is straightforward to check that we have $\mathcal A(a;p(0),\alpha,\Omega(\beta))$. We can now use the main induction hypothesis to infer
\begin{equation*}
(p(1),a)\vdash^{q(1)}_{q(1)}\Gamma\quad\text{with}\quad p(1)=p(0)+\varphi_0(\Omega(\beta)+t(0))\quad\text{and}\quad q(1)=\psi_\alpha p(1).
\end{equation*}
We have $p(0)<p=r+\varphi_0(s+t)$ by assumption, and the above yields
\begin{equation*}
\varphi_0(\Omega(\beta)+t(0))<\varphi_0(s+t)\in\mathsf H\subseteq\Gamma(\mathbf X).
\end{equation*}
Using Lemmas~\ref{lem:addition} and~\ref{lem:mathcal-A}, we obtain $p(1)<p$ and then $q(1)<q$. An application of weakening (Lemma~\ref{lem:weakening}) concludes the proof of the claim. Let us now consider an application of clause~(iii) from Definition~\ref{def:H-controlled-deriv}, where we have
\begin{equation*}
(r,a)\vdash^{t(0)}_s\Gamma,\varphi\quad\text{and}\quad(r,a)\vdash^{t(0)}_s\Gamma,\neg\varphi
\end{equation*}
for some $t(0)<t$ and some bounded $\mathbf L^u_{\Gamma(\mathbf X)}$-formula~$\varphi$ with $\rk(\varphi)<s$. First assume
\begin{equation*}
\rk(\varphi)<_{\mathbf O}\psi_{\alpha+1}0=\omega\cdot(1+\psi_{\alpha+1}0),
\end{equation*}
where the equality holds because~$\psi_{\alpha+1}0$ is strongly critical. From Lemma~\ref{lem:ranks-operators} we learn that $\varphi$ and $\neg\varphi$ are $\Sigma(\alpha)$-formulas. Given that any bounded formula~$\theta$ is equal to~$\theta^{L[\alpha]}$, the side induction hypothesis provides
\begin{equation}\label{eq:collapsing-cut}\tag{$\star$}
(p(0),a)\vdash^{q(0)}_{q(0)}\Gamma,\varphi\quad\text{and}\quad(p(0),a)\vdash^{q(0)}_{q(0)}\Gamma,\neg\varphi
\end{equation}
with $p(0)=r+\varphi_0(s+t(0))$ and $q(0)=\psi_\alpha p(0)$. Also by Lemma~\ref{lem:ranks-operators}, we have
\begin{equation*}
\rk(\varphi)\in\mathcal H_0(\supp^+(\varphi))\subseteq\mathcal H_r(a),
\end{equation*}
which entails $\rk(\varphi)<\psi_\alpha(r+1)\leq q(0)$ due to Lemma~\ref{lem:mathcal-A}. To conclude the present case of the induction step, we can thus reapply clause~(iii). Next, assume we have
\begin{equation*}
\psi_{\alpha+1}0\leq\rk(\varphi)\notin\{\psi_{\beta+1}0\,|\,\beta<\nu\}.
\end{equation*}
Due to $\rng(\Gamma(I))\ni\rk(\varphi)<s\in\overline K$, we may pick a $\beta<\nu$ with $\rk(\varphi)\leq\psi_{\beta+1}0<s$, by Lemma~\ref{lem:Omega-alpha-sup}. In the present case this upgrades to $\rk(\varphi)<\psi_{\beta+1}0$, which entails that we have $\alpha<\beta$. It follows that $\Gamma,\varphi,\neg\varphi$ consists of bounded $\Sigma(\beta)$-formulas. Indeed, for $\psi^{L[\alpha]}\in\Gamma$ with a $\Sigma(\alpha)$-formula $\psi$, we get
\begin{equation*}
\supp^{\mathbf L}_{\Gamma(X)}(\psi^{L[\alpha]})\subseteq\supp^{\mathbf L}_{\Gamma(X)}(\psi)\cup\{\mathbf R(\alpha)\}\subseteq_{\Gamma(\mathbf X)}\mathbf R(\beta).
\end{equation*}
From $\mathcal A(a;r,\alpha,s)$ and $\alpha<\beta$ we immediately get $\mathcal A(a;r,\beta,s)$. Thus the side induction hypothesis yields~(\ref{eq:collapsing-cut}), but now with $q(0)=\psi_\beta p(0)$ for the same $p(0)$. We can conclude the present case by the claim that we have established above. Finally, assume that we have $\rk(\varphi)=\psi_{\beta+1}0$ with $\alpha\leq\beta<\nu$. Recall that $\varphi$ and $\neg\neg\varphi$ are syntactically equal, due to our treatment of negation as a defined operation. We may thus assume that $\varphi$ (rather than $\neg\varphi$) is disjunctive. In view of Definition~\ref{def:ranks}, we must have $\varphi=\exists x\in L[\beta].\,\theta$ for some bounded $\mathbf L^u_{\Gamma(\mathbf X)}$-formula $\theta=\theta(x)$ that satisfies $\rk(\theta(0))<\psi_{\beta+1}0$. The latter entails that $\exists x.\,\theta$ is a $\Sigma(\beta)$-formula. Now the side induction hypothesis and boundedness (Proposition~\ref{prop:boundedness}) yield
\begin{equation*}
(p(0),a)\vdash^{q(0)}_{q(0)}\Gamma,(\exists x.\,\theta)^{q(0)}\quad\text{with } p(0)=r+\varphi_0(s+t(0))\text{ and } q(0)=\psi_\beta p(0).
\end{equation*}
From $(r,a)\vdash^{t(0)}_s\Gamma,\neg\varphi$ with $\neg\varphi=(\forall x\in L[\beta].\,\neg\theta)$ we also obtain
\begin{equation*}
(p(0),a)\vdash^{t(0)}_s\Gamma,(\forall x.\,\neg\theta)^{q(0)},
\end{equation*}
by weakening and Lemma~\ref{lem:inversion-variant}. Once readily derives $\mathcal A(a;p(0),\beta,s)$. As $(\forall x.\neg\theta)^{q(0)}$ is a bounded $\Sigma(\beta)$-formula, the side induction hypothesis provides
\begin{equation*}
(p(1),a)\vdash^{q(1)}_{q(1)}\Gamma,(\forall x.\neg\theta)^{q(0)}\quad\text{with }p(1)=p(0)+\varphi_0(s+t(0))\text{ and }q(1)=\psi_\beta p(1).
\end{equation*}
As we have $p(0)<p(1)$ and $q(0)<q(1)\in\mathcal H_{p(1)}(a)$ by Lemma~\ref{lem:mathcal-A}, the above can be weakened to
\begin{equation*}
(p(1),a)\vdash^{q(1)}_{q(1)}\Gamma,(\exists x.\,\theta)^{q(0)}.
\end{equation*}
Note that we have $p(1)<p$, due to Lemma~\ref{lem:addition}. Using Lemma~\ref{lem:ranks-operators}, we also see that $\rk(\theta(0))<\psi_{\beta+1}0=\omega\cdot(1+\psi_{\beta+1}0)$ entails
\begin{equation*}
\supp^+\left((\exists x.\,\theta)^{q(0)}\right)\subseteq\supp^+(\theta(0))\cup\{q(0)\}\subseteq_{\mathbf O}\psi_{\beta+1}0
\end{equation*} 
and hence $\rk((\exists x.\,\theta)^{q(0)})<\psi_{\beta+1}0=\rk(\varphi)<s$. We can thus conclude by the claim that was shown above (with $p(1)$ and $q(1)$ at the place of $p(0)$ and $q(0)$).
\end{proof}

One can use collapsing and boundedness to obtain quantitative information from proofs, as in~\cite[Theorem~4.9]{buchholz-local-predicativity}. For our purpose, it will be enough to have the following consistency result (recall that the empty sequent represents contradiction). Let us stress that our ordinal analysis was conditional on Assumptions~\ref{ass:u} and~\ref{ass:collapse-search}. In fact, our aim was to refute these assumptions. This aim is achieved by the following result, since it contradicts Theorem~\ref{thm:embedding} (embedding). The conclusions from this contradiction will be drawn in the next section.

\begin{corollary}[Consistency]\label{cor:consistency}
We do not have $(0,\emptyset)\vdash^t_{\Omega(\nu)}\langle\rangle$ for any~$t\in\mathbf O$.
\end{corollary}
\begin{proof}
Assume the claim is false. Then the previous theorem yields a $p\in\mathbf O$ with
\begin{equation*}
(p,\emptyset)\vdash^q_q\langle\rangle\quad\text{for}\quad q=\psi_0p.
\end{equation*}
Note that we have $q=\varphi(0,q)\leq\psi_{\alpha+1}0$ for all $\alpha<\nu$. We can thus use predicative cut elimination (Proposition~\ref{prop:pred-cut-elim}) to get
\begin{equation*}
(p,\emptyset)\vdash^{\varphi(q,q)}_0\langle\rangle.
\end{equation*}
The latter cannot hold, because no clause from Definition~\ref{def:H-controlled-deriv} applies: clauses~(i,ii) and~(iv) require a formula in~$\langle\rangle$, while clause~(iii) demands $\rk(\psi)<0$.
\end{proof}

\section{Fixed points, comprehension, and admissible sets}\label{sect:conclude}

In this section, we combine our previous work in order to prove Theorem~\ref{thm:main} and its corollaries, which were stated in the introduction. The following result provides the most difficult implication. It relies on an extensive argument that was developed in Sections~\ref{sect:search-trees} to~\ref{sect:ordinal-analysis}. More intuitive explanations of the following proof can be found in the introduction and in Section~\ref{sect:search-trees}.

\begin{theorem}\label{thm:main-crucial-implication}
For the following statements from Theorem~\ref{thm:main}, the theory $\atrs$ proves that (ii) implies~(iv) for any infinite ordinal~$\nu$:
\begin{enumerate}[label=(\roman*)]
\setcounter{enumi}{1}
\item any dilator has a well founded $\nu$-fixed point,
\setcounter{enumi}{3}
\item for any set~$u$, there is a sequence of admissible sets $\mathsf{Ad}_\alpha\ni u$ for $\alpha<\nu$, such that $\alpha<\beta<\nu$ entails $\mathsf{Ad}_\alpha\in\mathsf{Ad}_\beta$.
\end{enumerate}
\end{theorem}
\begin{proof}
As mentioned before, the restriction to infinite~$\nu$ is convenient because it allows us to reduce to the limit case. Indeed, it entails that we have $\nu\leq\mu+\omega$ for limits $\mu,\omega\leq\nu$. Given that~(ii) holds for $\nu$, it does also hold for~$\mu$ and for~$\omega$, by Corollary~\ref{cor:no-fp-monotone} in conjunction with Corollary~\ref{cor:nu-fp-unique} and Theorem~\ref{thm:psi-is-fp}. Assuming the limit case of the present theorem, we thus get~(iv) for~$\mu$ and for~$\omega$. To deduce~(iv) for~$\nu$ and a given set~$u$, we build two increasing sequences of admissibles $\mathsf{Ad}'_\alpha\ni u$ for $\alpha<\mu$ and $\mathsf{Ad}''_n\ni\bigcup_{\alpha<\mu}\mathsf{Ad}'_\alpha$ for~$n<\omega$. Note that we always have $\mathsf{Ad}'_\alpha\in\mathsf{Ad}''_n$, as admissible sets are transitive. To obtain the desired sequence of admissibles $\mathsf{Ad}_\alpha$ for $\alpha<\nu$, we set $\mathsf{Ad}_\alpha:=\mathsf{Ad}'_\alpha$ when $\alpha<\mu$ and $\mathsf{Ad}_\alpha:=\mathsf{Ad}''_n$ when $\alpha=\mu+n<\nu$. For the rest of this proof, we assume that $\nu$ is a limit such that~(ii) holds. Note that $\Pi^1_1$-comprehension becomes available by Corollary~\ref{cor:1-fixed-points}. It suffices to establish~(iv) for transitive~$u$ (replace~$u$ by the transitive closure~$u'$ of~$\{u\}$). We may also assume that the intersection $o(u)=u\cap\ord$ with the class of ordinals is a successor~$o(u)>1$ (replace~$u'$ by $u'\cup\{0,1,o(u')\}$). Since $\atrs$ contains the axiom of countability (cf.~the introduction), we can fix enumerations $u=\{u_i\,|\,i\in\mathbb N\}$ and $\nu=\{\nu_i\,|\,i\in\mathbb N\}$. By these preliminary considerations we have satisfied Assumption~\ref{ass:u}. Aiming at a contradiction, we now assume that~(iv) fails for~$\nu$ and~$u$ as fixed. By Proposition~\ref{prop:S^E-dilator}, it follows that a certain predilator~$\mathbf S_0$ is a dilator. The latter gives rise to another dilator~$\Gamma\circ\mathbf S$, due to Proposition~\ref{ref:Gamma-dilator} and Definition~\ref{def:X+S^E}. We now use statement~(ii) of the present theorem, which yields a well order~$\mathbf Y$ with a $\nu$-collapse
\begin{equation*}
\pi_{\mathbf Y}:\mathbf Y\to\nu\times(\Gamma\circ\mathbf S)(\mathbf Y).
\end{equation*}
This means that Assumption~\ref{ass:collapse-search} is satisfied as well. However, we have seen that the cited assumptions entail two incompatible results: Theorem~\ref{thm:embedding} and Corollary~\ref{cor:consistency} cannot both be valid, as we have $\Omega(\nu)=\Gamma_{\mathbf X+\langle\rangle}$ by Definition~\ref{def:overline-K}. Thus we have reached the desired contradiction.
\end{proof}

The next implication follows from~\cite[Paragraph~3]{Rathjen_PhD_1988} (see also the English translation in~\cite[Section~5]{Rathjen_PhD_2013} as well as Section~3.3.5 of the survey~\cite{pohlers98}). We provide a proof because the cited references involve the notion of inductive definition.

\begin{proposition}\label{prop:admissible-to-Pi11-rec}
Over $\atrs$, statement~(iv) from Theorem~\ref{thm:main} (or Theorem~\ref{thm:main-crucial-implication}) entails the following, for any ordinal~$\nu$:
\begin{enumerate}[label=(\roman*)]
\item $\Pi^1_1$-recursion along $\nu$ holds.
\end{enumerate}
\end{proposition}
\begin{proof}
We want to establish recursion for a given $\Pi^1_1$-formula~$\varphi(x,\alpha,X,\mathbf Z)$ with parameters $x\in\mathbb N$, $\alpha<\nu$ and $X,\mathbf Z\subseteq\mathbb N$. Recall (e.\,g.~from \cite[Lemma~V.1.4]{simpson09}) that we have a set theoretic $\Sigma$-formula~$\psi(x,\alpha,X,\mathbf Z)$ such that our base theory proves
\begin{equation*}
\text{``$A$ is admissible"}\to\forall x,\alpha,X,\mathbf Z\in A\,\left(\varphi(x,\alpha,X,\mathbf Z)\leftrightarrow\psi(x,\alpha,X,\mathbf Z)^A\right),
\end{equation*}
where the superscript denotes relativization. Since the cited reference employs inductive definitions, we recall an alternative argument: We have $\varphi(x,\alpha,X,\mathbf Z)$ precisely when a certain computable tree $T=T(x,\alpha,X,\mathbf Z)$ is well founded (see e.\,g.~\cite[Lemma~V.1.4]{simpson09}). Let $\psi(x,\alpha,X,\mathbf Z)$ assert that there is an $f:T\to\ord$ that descends along branches. Crucially, if $T\in A$ is indeed well founded, then such an $f$ exists in~$A$ (see e.\,g.~\cite[Theorem~4.6]{jaeger-admissibles}). In the following, we rely on the presentation of $\Pi^1_1$-recursion in the second paragraph after Theorem~\ref{thm:main}. Note that statement~(iv) holds for $\nu+1$ if it holds for~$\nu>0$. We may thus consider a sequence of admissibles $\mathsf{Ad}(\alpha)\in\mathsf{Ad}(\beta)$ for $\alpha<\beta\leq\nu$, such that $\mathsf{Ad}(0)$ contains given parameters~$\mathbf Z$. By primitive recursion in the sense of~\cite{jensen-karp}, we define a function $\nu\ni\alpha\mapsto Y^\alpha$ with $Y^0:=\emptyset$ and
\begin{align*}
Y^{\alpha+1}&:=Y^\alpha\cup\{\langle\alpha,x\rangle\,|\,x\in\mathbb N\text{ and }\psi(x,\alpha,Y^\alpha,\mathbf Z)^{\mathsf{Ad}(\alpha+1)}\},\\
Y^\lambda&:=\textstyle\bigcup_{\alpha<\lambda}Y^\alpha\quad\text{for limit $\lambda$}.
\end{align*}
We then set $Y:=\bigcup_{\alpha<\nu}Y^\alpha$ and observe $Y^\alpha=\{\langle\gamma,x\rangle\in Y\,|\,\gamma<\alpha\}$ for~$\alpha<\nu$, as in the presentation after Theorem~\ref{thm:main}. Our task is to establish
\begin{equation*}
\{x\in\mathbb N\,|\,\langle\alpha,x\rangle\in Y\}=\{x\in\mathbb N\,|\,\varphi(x,\alpha,Y^\alpha,\mathbf Z)\},
\end{equation*}
where the left side is commonly denoted by~$Y_\alpha$. The claim reduces to
\begin{equation*}
\varphi(x,\alpha,Y^\alpha,\mathbf Z)\leftrightarrow\psi(x,\alpha,Y^\alpha,\mathbf Z)^{\mathsf{Ad}(\alpha+1)}.
\end{equation*}
This equivalence holds by the choice of~$\psi$, once we have established $Y^\alpha\in\mathsf{Ad}(\alpha+1)$. We show the latter by induction on~$\alpha<\nu$. In the crucial case of a limit $\alpha$, we get
\begin{equation*}
\psi(x,\gamma,Y^\gamma,\mathbf Z)^{\mathsf{Ad}(\gamma+1)}\leftrightarrow\psi(x,\gamma,Y^\gamma,\mathbf Z)^{\mathsf{Ad}(\alpha)}\quad\text{for}\quad\gamma<\alpha.
\end{equation*}
Indeed, both sides are equivalent to $\varphi(x,\gamma,Y^\gamma,\mathbf Z)$, as $Y^\gamma\in\mathsf{Ad}(\gamma+1)\subseteq\mathsf{Ad}(\alpha)$~holds by induction hypothesis. So we can view $\alpha\geq\gamma\mapsto Y^\gamma$ as dependent on~$\mathsf{Ad}(\alpha)$ rather than $\alpha\ni\gamma\mapsto\mathsf{Ad}(\gamma)$. Now since $\mathsf{Ad}(\alpha+1)$ contains $\mathsf{Ad}(\alpha)$, it will also contain $Y^\alpha$, as admissible sets are closed under primitive recursive set functions.
\end{proof}

In Section~\ref{sect:fixed-points-exist} we have constructed a linear order~$\psi_\nu(D)$, relative to a given well order~$\nu$ and predilator~$D$. Besides the statements~(i,ii) and~(iv) that that have been recalled above, Theorem~\ref{thm:main} involves the following assertion:
\begin{enumerate}[label=(\roman*)]\setcounter{enumi}{2}
\item if~$D$ is a dilator (rather than just a predilator), then $\psi_\nu(D)$ is a well order.
\end{enumerate}
We now combine the previous results in order to deduce our main theorem.

\begin{proof}[\textbf{Proof of Theorem~\ref{thm:main}}] Due to Corollary~\ref{cor:nu-fp-unique} and Theorem~\ref{thm:psi-is-fp}, the order~$\psi_\nu(D)$ is the unique $\nu$-fixed point of~$D$, up to isomorphism. Together with Theorem~\ref{thm:Pi11-rec-to-wf}, it follows that we have
\begin{equation*}
(i)\quad\Rightarrow\quad (ii)\quad\Leftrightarrow\quad (iii)
\end{equation*}
for any well order~$\nu$, provably in $\rca_0$. As in the desired Theorem~\ref{thm:main}, we now assume that $\nu$ is infinite (though this could probably be avoided). From Corollary~\ref{cor:1-fixed-points} we know that~(ii) entails $\Pi^1_1$-comprehension. To show that (ii) implies~(i) over the theory~$\rca_0$, it is thus enough to prove the same implication in $\mathsf{ATR}_0$ or indeed in the conservative extension $\atrs$, a set theory due to Simpson. As stated in the introduction, our version of $\atrs$ contains the axiom of countability, which is included in~\cite{simpson09} but marked as `optional' in~\cite{simpson82}. Also recall that $\atrs$ contains axiom beta, which allows us to assume that $\nu$ is an ordinal (rather than just a well order). Over~$\atrs$, Theorem~\ref{thm:main-crucial-implication} and Proposition~\ref{prop:admissible-to-Pi11-rec} yield
\begin{equation*}
(ii)\quad\Rightarrow\quad (iv)\quad\Rightarrow\quad (i),
\end{equation*}
which closes our circle of implications.
\end{proof}

In the introduction, we have stated a corollary which asserts that (ii) and (iii) for~$\nu=\omega$ are equivalent to the following:
\begin{enumerate}[label=(\roman*')]
\item every subset of~$\mathbb N$ is contained in a countable $\beta$-model of $\Pi^1_1$-comprehension.
\end{enumerate}
This result holds by our main theorem and the following standard argument.

\begin{proof}[\textbf{Proof of Corollary~\ref{cor:beta-model-Pi11-CA}}]
We first assume~(i') and derive (ii) for~$\nu=\omega$, over~$\rca_0$. In fact we may work in $\mathsf{ATR}_0$ (e.\,g.~by~\cite[Exercise~VII.2.10]{simpson09}). Due to Theorem~\ref{thm:main}, it is enough to establish $\Pi^1_1$-recursion along~$\omega$. Given a $\Pi^1_1$-formula $\varphi(x,n,X,\mathbf Z)$ and paramters~$\mathbf Z$, we invoke~(i') to get a countable $\beta$-model $\mathcal M\ni\mathbf Z$ of $\Pi^1_1$-comprehension. Satisfaction in $\mathcal M$ is arithmetical for instances of~$\varphi$ (cf.~\cite[Definition~VII.2.1]{simpson09}). We can thus use arithmetical recursion to construct the set
\begin{equation*}
Y=\{\langle n,x\rangle\in\omega\times\mathbb N\,|\,\mathcal M\vDash\varphi(x,n,Y^n,\mathbf Z)\},
\end{equation*}
with $Y^n=\{\langle m,x\rangle\in Y\,|\,m<n\}$ as before. The given definition presumes $Y^n\in\mathcal M$, which we get by induction: in the step, $\Pi^1_1$-comprehension in~$\mathcal M$ yields
\begin{equation*}
Y^{n+1}=Y^n\cup\{\langle n,x\rangle\,|\,x\in\mathbb N\text{ and }\mathcal M\vDash\varphi(x,n,Y^n,\mathbf Z)\}\in\mathcal M.
\end{equation*}
Since $\mathcal M$ is a $\beta$-model (cf.~\cite[Lemma~VII.2.6]{simpson09}), we have
\begin{equation*}
\varphi(x,n,Y^n,\mathbf Z)\quad\leftrightarrow\quad\mathcal M\vDash\varphi(x,n,Y^n,\mathbf Z).
\end{equation*}
In the notation from the introduction we thus have $H_\varphi(Y)$, as needed to establish the given instance of $\Pi^1_1$-recursion. To show that (ii) for~$\nu=\omega$ entails (i'), we may work over $\atrs$, as in the proof of Theorem~\ref{thm:main}. By the latter, we get a hierarchy of admissible sets $\ad(m)\in\ad(n)$ for $m<n<\omega$, where we can assume that $\ad(0)$ contains a given subset of~$\mathbb N$. Let us put
\begin{equation*}
\mathcal S:=\{Z\in A\,|\,Z\subseteq\mathbb N\}\quad\text{with}\quad A:=\textstyle\bigcup_{n<\omega}\ad(n).
\end{equation*}
We shall show that $\mathcal M:=(\mathbb N,\mathcal S)$ is the $\beta$-model required by~(i'). First note that the countability of~$\mathcal S$ is for free, because~$\atrs$ includes an axiom that makes all sets countable (cf.~the previous proof). To show that $\mathcal M$ is a $\beta$-model, we consider an arbitrary $\Pi^1_1$-formula~$\varphi(x,Z)$. As in the proof of Proposition~\ref{prop:admissible-to-Pi11-rec}, we obtain a \mbox{$\Sigma$-formula}~$\psi(x,Z)$ such that $\varphi(x,Z)$ and $\psi(x,Z)^{\ad(n)}$ are equivalent for~$Z\in\ad(n)$. The indicated proof of equivalence relativizes to~$A$ (for details see \cite[Section~7]{jaeger-admissibles} or~\cite[Section~3.3.2]{pohlers98}, noting that $A\vDash\mathsf{KPl}^{\operatorname{r}}$). This means that we get
\begin{equation*}
\varphi(x,Z)\leftrightarrow\psi(x,Z)^{\ad(n)}\leftrightarrow\mathcal M\vDash\varphi(x,Z)\quad\text{when}\quad Z\in\ad(n).
\end{equation*}
As any $Z\in A$ is contained in $\ad(n)$ for some~$n\in\mathbb N$, it follows that $\mathcal M$ is a $\beta$-model. Invoking bounded separation in~$\ad(n+1)$, we also see that $Z\in\ad(n)$ entails
\begin{equation*}
\{x\in\mathbb N\,|\,\mathcal M\vDash\varphi(x,Z)\}=\{x\in\mathbb N\,|\,\psi(x,Z)^{\ad(n)}\}\in\ad(n+1)\subseteq A,
\end{equation*}
which shows that $\mathcal M$ satisfies $\Pi^1_1$-comprehension.
\end{proof}

To conclude this paper, we derive the final result that was stated in the introduction. It is concerned with the principle of $\Pi^1_1$-transfinite recursion, which asserts that statement~(i) of Theorem~\ref{thm:main} holds for every well order~$\nu$.

\begin{proof}[\textbf{Proof of Corollary~\ref{cor:Pi11-TR}}]
Consider the statements~(i) to~(iii) from Theorem~\ref{thm:main}. For each of these statements, we define the variants
\begin{enumerate}
\item[($\forall$n)] statement (n) holds for every well order~$\nu$,
\item[($\infty$n)] statement (n) holds for every infinite well order~$\nu$.
\end{enumerate}
By Theorem~\ref{thm:main}, statements ($\infty$i) and ($\infty$ii) and ($\infty$iii) are pairwise equivalent. The corollary claims that the same holds for ($\forall$i) and ($\forall$ii) and ($\forall$iii). This is true because statements ($\forall$n) and ($\infty$n) are in fact equivalent. The latter is immediate in the case of~(i). For the other statements, it follows from Corollary~\ref{cor:no-fp-monotone} (in~conjunction with Corollary~\ref{cor:nu-fp-unique} and Theorem~\ref{thm:psi-is-fp}).
\end{proof}

\bibliographystyle{amsplain}
\bibliography{Hierarchies-admissibles_Freund-Rathjen}

\end{document}